\documentclass[11pt]{article} 
    
\usepackage{a4wide,authblk,epsfig,amsmath,amssymb,makeidx,graphics,color,tikz,stmaryrd} 
           
\title{Geometric aspects of the symmetric inverse M-matrix problem}
\date{\today} 
\author{Jan Brandts and Apo Cihangir} 
\affil{\small Korteweg-de Vries Institute for Mathematics, Faculty of Science, University of Amsterdam, P.O.~Box 94248, 1090 GE, Amsterdam, Netherlands.\\
E-mail: J.H.Brandts@UvA.nl, A.Cihangir@UvA.nl\\
\bigskip
In memory of Professor Miroslav Fiedler (1926--2015)}  
                       
\begin{document}                                
               
\newtheorem{Th}{Theorem}[section]                 
\newtheorem{Le}[Th]{Lemma}      
\newtheorem{Co}[Th]{Corollary}        
\newtheorem{Pro}[Th]{Proposition}         
\newtheorem{Def}[Th]{Definition}   
\newtheorem{Obs}[Th]{Observation}   
\newtheorem{Con}[Th]{Conjecture}         
\newtheorem{rem}[Th]{Remark}        
\newcommand{\be}{\begin{equation}}        
\newcommand{\ee}{\end{equation}}          
\newcommand{\RR}{\mathbb{R}}      
\newcommand{\half}{\frac{1}{2}}    
\newcommand{\hdrie}{\hspace{3mm}} 
\newcommand{\und}{\hdrie\mbox{\rm and }\hdrie} 
\newcommand{\sth}{\hdrie | \hdrie}
\newcommand{\argmax}{\mbox{\rm argmax}}   
\newcommand{\inter}{\mbox{\rm int}} 
\def\zeros {{\bf 0}} 
\def\ones{{\bf 1}} 
\def\CC {\mathcal{C}} 
\def\DD {\mathcal{D}} 
\def\PP {\mathcal{P}}
\def\Bc {\mathcal{B}}
\def\VV {\mathcal{V}}
\def\PP {\mathcal{P}}
\def\Bn {\mathcal{B}_n} 
\def\bb {\mathbb{b}}
\def\GG {\mathcal{G}}
\def\UU {\mathcal{U}} 
\def\II {\mathbb{I}} 
\def\BB {\mathbb{B}}   
\def\MM {i\!\mathcal{M}} 
\def\HH {\mathcal{H}}
\def\PP {\mathcal{P}}  
\def\FF {\mathcal{F}}  
\def\SIS {\mathcal{S}} 
\def\TT {\mathcal{T}}
\def\sort {{\rm sort}}  
\def\Eqv {\hdrie\Leftrightarrow\hdrie}  
\def\Bnk {\BB^{n\times k}}
\def\Bnn {\BB^{n\times n}}
\def\Mnn {\MM^{n\times n}} 
\def\Mnk {\MM^{n\times k}}  
\def\mod {\,{\rm mod }\,}  
\def\conv {\mbox{\rm conv}}  
\def\vol {\mbox{\rm vol}} 
\newcommand{\supp}{{\rm supp}}   
\newcommand{\ol}{\overline}  
\newcommand{\vep}{\varepsilon}  
\newcommand{\zok}{$0/1$-$k$} 
\newcommand{\zo}{$0/1$} 
\newcommand{\tet}{\boxslash}
\newcommand{\dd}{{\rm dd}} 
\newcommand{\four}{\boxtimes}
\newcommand{\Rspd}{\mathbb{R}_{\rm spd}}
\newcommand{\Rnspd}{\mathbb{R}_{\rm spd}^{n\times n}}
\newcommand{\Rnn}{\mathbb{R}^{n\times n}}
\newcommand{\npo}{\!n\!+\!1\!}
\newcommand{\npt}{\!n\!+\!2\!}
\newcommand{\nmo}{\!n\!-\!1\!}
    
\maketitle  

\begin{abstract}
We investigate the symmetric inverse M-matrix problem from a geometric perspective. The central question in this geometric context is, which conditions on the $k$-dimensional facets of an $n$-simplex $S$ guarantee that $S$ has no obtuse dihedral angles. The simplest of such conditions is that if all triangular facets of $S$ are equilateral, then $S$ is regular and thus nonobtuse. First we study the properties of an $n$-simplex $S$ whose $k$-facets are all nonobtuse, and generalize some classical results by Fiedler~\cite{Fie}. We prove that if all $(\nmo)$-facets of an $n$-simplex $S$ are nonobtuse, each makes at most one obtuse dihedral angle with another facet. This helps to identify a special type of tetrahedron, which we will call sub-orthocentric, with the property that if all tetrahedral facets of $S$ are sub-orthocentric, then $S$ is nonobtuse. Rephrased in the language of linear algebra, this constitutes a purely geometric proof of the fact that each symmetric ultrametric matrix~\cite{VarNab} is the inverse of a weakly diagonally dominant M-matrix. The geometric proof provides valuable insights that supplement the discrete measure theoretic~\cite{MaMiSa} and the linear algebraic~\cite{NabVar} proof. 

\smallskip

The review papers~\cite{Joh,JoSm} support our believe that the linear algebraic perspective on the inverse M-matrix problem dominates the literature. The geometric perspective however connects sign properties of entries of inverses of a symmetric positive definite matrix to the dihedral angle properties of an underlying simplex, and enables an explicit visualization of how these angles and signs can be manipulated. This will serve to formulate purely geometric conditions on the $k$-facets of an $n$-simplex $S$ that may render $S$ nonobtuse also for $k>3$. For this, we generalize the class of sub-orthocentric tetrahedra that gives rise to the class of ultrametric matrices, to sub-orthocentric simplices that define symmetric positive definite matrices $A$ with special types of $k\times k$ principal submatrices for $k>3$. Each sub-orthocentric simplices is nonobtuse, and we conjecture that any simplex with sub-orthocentric facets only, is sub-orthocentric itself. 

\smallskip

Along the way, several additional new concepts will be introduced, such as vertex Gramians, simplicial matrix classes, the dual hull and the sub-orthocentric set of a nonobtuse simplex, and nonblocking matrices. These concepts may also be of use in a different linear algebraic setting, in particular in the context if completely positive matrices, for which we will prove some auxilliary results as well. 
\end{abstract} 
\section{Introduction}\label{Mintro}
An $n$-{\em simplex}\index{simplex} is the convex hull of $\npo$ affinely independent {\em vertices} $v_0,\dots,v_n$ in Euclidean $n$-space. Its geometric properties can be well studied in terms of graphs and linear algebra, as demonstrated in over six decades of work of Miroslav Fiedler~\cite{Fie,Fie2,Fie3,Fie4,FiPt}. Fiedler usually describes simplices with the $(\npo)\times(\npo)$ {\em Menger matrix} 
\index{Menger matrix} $M=(m_{ij})$ whose entry $m_{ij}$ equals the squared distance between the vertices $v_i$ and $v_j$, or a bordered $(\npt)\times(\npt)$ variation of it. We will use a different approach and describe an $n$-simplex $S$ using {\em vertex Gramians}. These are symmetric positive definite $n\times n$ matrices. They are defined as follows. Choose an orthogonal coordinate system with the origin located at a vertex $v_\ell$ of $S$. Then the coordinate vectors in $\RR^n$ of the remaining $n$ vertices of $S$ are linearly independent. Any $n\times n$ matrix $P_\ell$ having these vectors as columns fully captures the geometry of $S$. Now, observe that column permutations, left-multiplications with an orthogonal matrix, and choosing the origin at another vertex of $S$ may lead to an entirely different matrix $P_\ell$. However, since $P_\ell^\top P_\ell = R^\top R$ if and only if $R=UP_\ell$ for some orthogonal transformation $U$, we can nullify the effect of the orthogonal transformations $U$ by considering the {\em Gramians} of the matrices $P_\ell$.  
\begin{Def}[Vertex Gramian] \index{vertex Gramian} \label{Mdef-1}{\rm Let $S$ be an $n$-simplex with vertices $v_0,\dots,v_n\in\RR^n$.  A {\em vertex Gramian} $G_\ell$ of $S$ associated with vertex $v_\ell$ is the Gramian $G_\ell=P_\ell^\top P_\ell$ of $P_\ell$, where $P_\ell$ is a 
$n \times n$ matrix whose $n$ columns are the vectors $\{ v_k-v_\ell\,|\, k\not=\ell\}$.}
\end{Def} 
With this definition, any pair $G_\ell^1,G_\ell^2$ of vertex Gramians corresponding to the {\em same} vertex $v_\ell$ of $S$ are simultaneous row- and column permutations of one another. We will use the notation $G_\ell^1\overset{\pi}\sim G_\ell^2$ to denote such a {\em permutation equivalence}. Figure~\ref{Mfigure21} illustrates vertex Gramians $G_u$ and $G_v$ for the drawn triangle $\Delta$, corresponding to two {\em different} vertices $u$ and $v$.
\begin{figure}[h]
\begin{center}  
\begin{tikzpicture}[scale=0.8, every node/.style={scale=0.8}]
\draw (0,0)--(3,0)--(6,3)--cycle; 
\draw[fill=black] (0,0) circle [radius=0.07];
\draw[fill=black] (6,3) circle [radius=0.07];
\draw[gray, fill=gray] (3,0) circle [radius=0.07];
\node[scale=0.9] at (0,-0.5) {$(0,0)$};
\node[scale=0.9] at (3,-0.5) {$(1,0)$};
\node[scale=0.9] at (6.4,2.6) {$(2,1)$};
\node[scale=0.9,gray] at (2.7,0.4) {$(0,0)$};
\node[scale=0.9,gray] at (-0.5,0.4) {$(-1,0)$};
\node[scale=0.9,gray] at (5,3) {$(1,1)$};
\node at (3,1) {$\Delta$};
\node at (-0.3,0) {$u$};
\node[gray] at (3.3,0) {$v$};
\node[scale=0.9,gray] at (-1,1.4) {$P_v=\left[\begin{array}{rr} -1 & 1 \\ 0 & 1\end{array}\right]$};
\node[scale=0.9] at (-1,2.7) {$P_u=\left[\begin{array}{rr} \phantom{-}1 & 2 \\ 0 & 1\end{array}\right]$};
\node[scale=0.9] at (9.5,1.5) {$G_u=\left[\begin{array}{rr} \phantom{-}1 & \phantom{-}2 \\ 2 & 5\end{array}\right] \sim \left[\begin{array}{rr} 1 & -1 \\ -1 & 2\end{array}\right]=G_v$};
\node[scale=0.9] at (9.5,0.3) {$G_u^{-1}=\left[\begin{array}{rr} 5 & -2 \\ -2 & 1\end{array}\right] \not \sim \left[\begin{array}{rr} \phantom{-}2 & \phantom{-}1 \\ 1 & 1\end{array}\right]=G_v^{-1}$};
\end{tikzpicture}    
\end{center} 
\caption{\small{Two vertex Gramians $G_u\sim G_v$ for a triangle $\Delta$. Only $G_u$ is an inverse M-matrix.}}
\label{Mfigure21}
\end{figure}
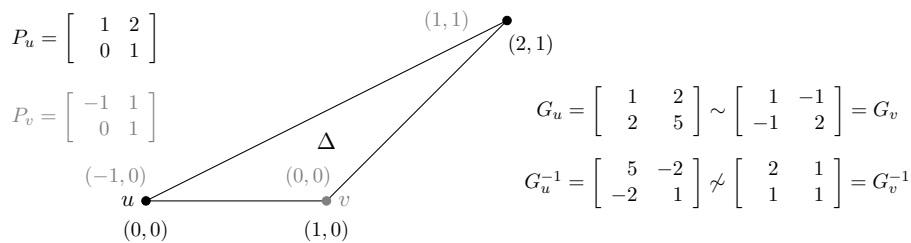 \\[2mm] 
The observation in the following definition is trivial but relevant in the context of this thesis. Write $\Rspd^{n\times n}$ for the set of $n\times n$ real symmetric positive definite matrices.   
\begin{Def}[Underlying simplex] \label{Mdef-2} {\rm Each $A\in\Rspd^{n\times n}$ is a vertex Gramian of some $n$-simplex $S$. We will call $S$ the {\em underlying simplex} \index{underlying simplex} of $A$ and denote their relation by $S=\SIS(A)$}.
\end{Def}      
Definitions~\ref{Mdef-1} and~\ref{Mdef-2} invite to define and investigate the equivalence relation $\sim$ on $\Rnspd$ given by 
\be\label{Meqv} A \sim B \hdrie \Leftrightarrow \hdrie \SIS(A) = \SIS(B). \ee
The equivalence class $\GG(A)$ of $A\in\Rnspd$ with respect to the relation~(\ref{Meqv}) consists precisely of all vertex Gramians of the underlying $n$-simplex $\SIS(A)$ and includes all $B\overset{\pi}{\sim}A$.     
\begin{Def}[Simplicial matrix class]\label{Msmc}  {\rm A set $\CC\subset\Rnspd$ is a {\em simplicial matrix class} \index{simplicial matrix class} if 
\be A\in \CC \hdrie \Rightarrow \hdrie \GG(A) \subset \CC, \ee
or in other words, if $\CC$ consists of all the vertex Gramians of some collection of $n$-simplices}.
\end{Def}   
The class $\MM$ of \index{symmetric inverse M-matrices} {\em symmetric inverse M-matrices}~\cite{Joh,JoSm,Wil}, which is the main topic of investigation in this chapter, is {\em not} simplicial; it is well-known to be closed under $\overset{\pi}{\sim}$, but it is not closed under $\sim$. Figure~\ref{Mfigure21} may serve as proof of this statement. On the other hand, the subclass $\MM_\dd\subset\MM$ of inverses of {\em weakly diagonally dominant} symmetric M-matrices is indeed simplicial. We will show this in Section~\ref{Msect-2}, where we also review~\cite{BrKoKr,BrKoKrSo} that $\MM_\dd$ equals the set of all vertex Gramians of simplices without any obtuse {\em dihedral angles} between its facets: 
\be\label{Minvm} \MM_\dd = \left\{ A\in\Rnspd \sth \mbox{\rm $\SIS(A)$ is a {\em nonobtuse} simplex}\right\}. \ee
This straightforward geometric characterization of $\MM_\dd$ does not seem to be common know\-ledge in the literature: the review papers \cite{Joh,JoSm} do not mention it. Similarly, the class $\DD$ of non-singular {\em doubly nonnegative matrices} \cite{BePl} is {\em not} simplicial. It is closed under $\overset{\pi}{\sim}$, but not under $\sim$. However, the subclass of {\em pointwise weakly diagonally dominant} matrices $\DD_\dd\subset\DD$ is again a simplical matrix class:
\be\label{Mdnn} \DD_\dd = \left\{ A\in\Rnspd \sth \mbox{\rm $\SIS(A)$ has only nonobtuse {\em triangular} facets}\right\}. \ee
Since any nonobtuse simplex has only nonobtuse facets, $\MM_\dd$ is a {\em simplicial subclass} of $\DD_\dd$, and both are extremes in the chain of simplical matrix classes,
\be\label{Mchain} \MM_\dd=\VV_n^n\subset \VV_{n-1}^n \subset \cdots  \subset  \VV_3^n \subset \VV_2^n = \DD_\dd \subset \VV_1^n = \Rnspd, \ee
where $\VV_k^n$ is the class of all vertex Gramians of $n$-simplices whose $k$-facets are all nonobtuse.\\[2mm]
To further illustrate the equivalence $\sim$, let $e_1,\dots,e_n$ be the standard basis vectors of $\RR^n$, and $S$ the $4$-simplex with vertices $v_0,v_1=v_0+e_1,v_2=v_1+e_2,v_3=v_2+e_3$, and $v_4=v_3+e_4$. Because it has a path of four mutually orthogonal edges, $S$ is called a {\em path simplex}\index{path simplex} \cite{Cox2,Schl}. Although $S$ has five vertices, due to a reflection symmetry in $S$ each vertex Gramian of $S$ is permutation equivalent to one of the vertex Gramians $G_0,G_1,G_2$ associated with the vertices $v_0,v_1,v_2$,
\be\label{Mexample}\small G_0=\left[\begin{array}{cccc} 1 & 1 & 1 & 1\\ 1 & 2 & 2 & 2\\ 1 & 2 & 3 & 3\\ 1 & 2 & 3 & 4\end{array}\right], \hdrie 
G_1=\left[\begin{array}{cccc} 1 & 0 & 0 & 0\\ 0 & 1 & 1 & 1\\ 0 & 1 & 2 & 2\\ 0 & 1 & 2 & 3\end{array}\right], \und G_2 =\left[\begin{array}{cccc} 2 & 1 & 0 & 0\\ 1 & 1 & 0 & 0\\ 0 & 0 & 1 & 1\\ 0 & 0 & 1 & 2\end{array}\right]. \ee
The matrices $G_0\sim G_1\sim G_2$ are considered of {\em different} types in the literature: $G_0$ is of {\em type-D} \index{type-D} as defined by Markham \cite{Mar}, whereas $G_1$ and $G_2$ are examples of {\em symmetric ultrametric matrices} \cite{MaMiSa,NabVar,VarNab}. We will show that the class $\UU$ of symmetric ultrametric matrices is simplicial, which establishes ultrametricity of $A$ as a {\em geometric} property of its underlying simplex $\SIS(A)$. Moreover, by showing that $\UU\subset\MM_\dd$, we provide a {\em geometric proof} of the fact that each symmetric ultrametric matrix is the inverse of a weakly diagonally dominant M-matrix, supplementing the proofs from discrete measure theory \cite{MaMiSa} and linear algebra~\cite{NabVar}. 
\begin{rem}{\rm It is an interesting question which other known matrix classes are simplicial, and also to investigate which matrix properties are invariant under $\sim$. An obvious example is that if $A\sim B$ then $|\det(A)| = |\det(B)|$. Also, invariants associated with simplices may have matrix equivalents that are of interest. For instance, writing $e=e_1+\cdots+e_n$ for the {\em all-ones} vector, the radii $r_{\rm{i}}(A)$ and $r_{\rm{c}}(A)$ of the {\em inscribed} and {\em circumscribed} ball about the underlying simplex $\SIS(A)$ of $(a_{ij})=A\in\Rnspd$, are given by
\be\label{Mx-3} \frac{1}{r_{\rm{i}}(A)} = \frac{1}{\sqrt{e_{\phantom{j}}^\top A^{-1}e}}+\sum_{j=1}^n \frac{1}{\sqrt{e_j^\top A^{-1}e_j}} \und r_{\rm{c}}(A) = \half\sqrt{v^\top A^{-1}v}, \hdrie \mbox{where} \hdrie v = \sum_{j=1}^n a_{jj}e_j.\ee 
Consequently, we have that $r_{\rm{i}}(A)=r_{\rm{i}}(B)$ and $r_{\rm{c}}=r_{\rm{c}}(B)$ whenever $A\sim B$. The formulas given in~(\ref{Mx-3}) can be derived using the 
linear algebraic description of simplices in Section~\ref{Msect-2}}.
\end{rem}          
\subsection{Detailed outline, motivation, and main results}
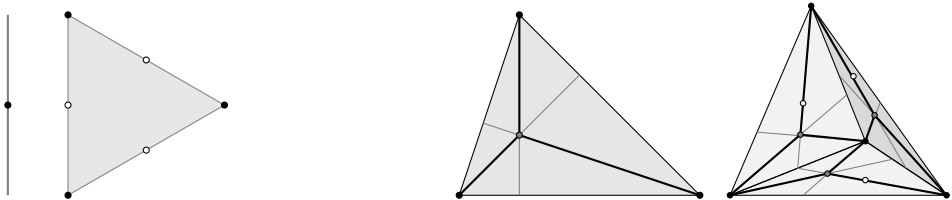
\begin{figure}[h]
\begin{center}  
\begin{tikzpicture}[scale=0.8, every node/.style={scale=0.8}]
\draw[gray,thick] (-1,0)--(-1,3);
\draw[fill=black] (-1,1.5) circle [radius=0.05];
\draw[gray,fill=white!80!gray] (0,0)--(0,3)--(2.6,1.5)--cycle;
\draw[fill=black] (0,0)     circle [radius=0.05];
\draw[fill=black] (0,3)     circle [radius=0.05];
\draw[fill=black] (2.6,1.5) circle [radius=0.05];
\draw[fill=white] (0,1.5)    circle [radius=0.05];
\draw[fill=white] (1.3,2.25) circle [radius=0.05];
\draw[fill=white] (1.3,0.75) circle [radius=0.05];
\begin{scope}[shift={(6.5,0)}]
\draw[white!80!gray, fill=white!80!gray] (0,0)--(4,0)--(1,3)--cycle;
\draw (0,0)--(4,0)--(1,3)--cycle;
\draw[fill=black] (0,0) circle [radius=0.05];
\draw[fill=black] (4,0) circle [radius=0.05];
\draw[fill=black] (1,3) circle [radius=0.05];
\draw[thick]      (1,3)--(1,1);
\draw[gray] (1,1)--(1,0);
\draw[thick]          (4,0)--(1,1);
\draw[gray] (1,1)--(0.4,1.2);
\draw[thick]          (0,0)--(1,1);
\draw[gray] (1,1)--(2,2);
\draw[fill=gray] (1,1) circle [radius=0.05];
\end{scope}
\begin{scope}[shift={(11,0)},scale=0.9]
\draw[white!90!gray,fill=white!90!gray] (0,0)--(4,0)--(1.5,3.5)--cycle;
\draw[white!70!gray,fill=white!70!gray] (4,0)--(2.5,1)--(1.5,3.5)--cycle;
\draw (0,0)--(4,0)--(2.5,1)--cycle;
\draw (0,0)--(2.5,1)--(1.5,3.5)--cycle;
\draw (1.5,3.5)--(4,0)--(2.5,1);
\draw[gray] (1.5,3.5)--(1.25,0.5)--(4,0);
\draw[gray] (2.5,1)--(0.5,7/6);
\draw[gray] (0,0)--(2.15,1.85);
\draw[gray] (0,0)--(3,2/3);
\draw[gray] (4,0)--(2.03,2.2);
\draw[gray] (1.5,3.5)--(3.25,0.5);
\draw[gray] (2.5,1)--(1.35,0);
\draw[gray] (2.5,1)--(2.78,1.73);
\draw[fill=black] (0,0) circle [radius=0.05];
\draw[fill=black] (4,0) circle [radius=0.05]; 
\draw[fill=black] (2.5,1) circle [radius=0.05];
\draw[fill=black] (1.5,3.5) circle [radius=0.05];
\draw[thick] (0,0)--(1.3,1.12);
\draw[thick] (1.5,3.5)--(1.3,1.12);
\draw[thick] (2.5,1)--(1.3,1.12);
\draw[fill=white] (1.35,1.7) circle [radius=0.05];
\draw[thick] (0,0)--(1.8,0.4);
\draw[thick] (2.5,1)--(1.8,0.4);
\draw[thick] (4,0)--(1.8,0.4);
\draw[fill=white] (2.5,0.28) circle [radius=0.05];
\draw[thick] (1.5,3.5)--(2.67,1.48);
\draw[thick] (4,0)--(2.67,1.48);
\draw[thick] (2.5,1)--(2.67,1.48);
\draw[fill=gray] (1.3,1.12) circle [radius=0.05];
\draw[fill=gray] (1.8,0.4) circle [radius=0.05];
\draw[fill=gray] (2.67,1.48) circle [radius=0.05];
\draw[fill=white] (2.28,2.2) circle [radius=0.05];
\end{scope}
\end{tikzpicture}    
\end{center} 
\caption{\small{ A {\em sub-orthocentric} simplex is one in which each vertex projects in the {\em sub-orthocentric set} of its opposite nonobtuse facet. The sub-orthocentric set of a {\em line segment} is its midpoint, and of a nonobtuse {\em triangle} it is the union of the three line segments between a vertex and the orthocenter. The projections of vertices are depicted as white bullets.}}
\label{Mfigure23}
\end{figure}

The newly defined concepts of vertex Gramians and simplicial matrix classes will serve as a framework to systematically study the {\em symmetric inverse M-matrix problem} \cite{DeMaMa,FiPt,Joh,JoSm,Mar,MaMiSa,NabVar,NabVar2,VarNab,Wil} from the geometric perspective provided by (\ref{Minvm}): we associate the symmetric inverse M-matrices to the simplices without obtuse dihedral angles, the so-called {\em nonobtuse} simplices \cite{BrKoKr,BrKoKrSo,Fie,Fie2,Fie3}. Their definition and properties will be reviewed in Section~\ref{Msect-2}, where we will also explain the subtle interrelation between the Stieltjes property and weak diagonal dominance. The geometric problem of most interest is, under which {\em additional} conditions a simplex $S$, all whose $(\nmo)$-dimensional {\em facets} are nonobtuse, is nonobtuse itself. This relates to the linear algebraic question of finding {\em additional} conditions on {\em submatrices} of a symmetric nonnegative matrix under which it is the inverse of a (possibly weakly diagonally dominant) M-matrix. The general setting for this will be provided in Section~\ref{Msect-3}, where we start by studying the properties that simplices whose $k$-facets are all nonobtuse already possess. To the best of our knowledge, such simplices have not been studied before, even though it has been frequently observed that a tetrahedron $\tet$ with nonobtuse triangular facets $\Delta$ need not be nonobtuse. On the other hand, it seems {\em not} to have been observed that such a $\tet$ cannot have {\em three} obtuse dihedral angles, like an arbitrary tetrahedron may have, but merely {\em two}. By introducing the {\em dual hull} $S^\ast$ of a nonobtuse simplex $S$ in Section~\ref{Msect-3.1}, derived from the concept of the dual of a simplicial cone \cite{BoVa}, we are able to prove in Section~\ref{Msect-3.2} that a facet of a simplex with nonobtuse facets makes at most one obtuse angle with another facet. Hence, as $n$ increases, simplices with nonobtuse facets can only have a very limited number of obtuse dihedral angles, in comparison with arbitrary simplices. See the table in Figure~\ref{Mfigure22}, and also Theorem~\ref{Mth-2} in Section~\ref{Msect-3}. 
\begin{figure}[h]
\begin{center}  
\begin{tikzpicture}[scale=0.8, every node/.style={scale=0.8}]
\node at (0,1.8) {$\begin{array}{|c||r|r|r|r|r|r|r|r|}
\hline  
n & 1 & 2 & 3 & 4 & 5 & 6 & 7 & 8\\
\hline
\hline
\omega(n) & 0 & 1 & 2 & 2 & 3 & 3 & 4 & 4 \\
\hline
\mu(n)     &  0 & 1 & 3 & 6 & 10 & 15 & 21 & 28 \\
\hline 
\end{array}$ };
\begin{scope}[shift={(3,0)}]
\node[scale=0.9] at (4.2,1.3) {$ 
\left[\begin{array}{cc|cc|cc}  
+ & + & \leq & \cdots &\cdots &\leq \\ 
+ & + & \leq &  &  & \vdots \\ 
\hline 
\leq & \leq &\ddots &\ddots &  & \vdots  \\
 \vdots  &  &\ddots &\ddots&\leq &\leq \\ 
\hline 
\vdots  &  & & \leq &+ &+ \\ 
\leq  & \cdots & \cdots & \leq &+ & + 
\end{array}\right]$};
\node at (-0.7,0.4) {$D\geq 0,\,\, C\geq 0, \hspace{3mm} A^{-1} \overset{\pi}{\sim}D-C =$};
\end{scope}
\end{tikzpicture}    
\end{center} 
\caption{\small{Table of the maximal number $\omega(n)$ of obtuse dihedral angles in an $n$-simplex with nonobtuse dihedral facets versus the corresponding maximal number $\mu(n)$ in {\em any} $n$-simplex. Right: block-$2\times 2$ {\em sign pattern} of the inverse of a vertex Gramian $A$ associated with such a simplex. Here, the $\leq$ sign stands for any nonpositive number and $+$ for a positive number.}}
\label{Mfigure22}
\end{figure} \\[2mm] 
The small number of obtuse dihedral angles in such simplices relates to a small number of positive entries in the inverse of a corresponding vertex Gramian $A\in\Rnspd$. In fact, $A^{-1}\overset{\pi}{\sim}D-C$ with $D\geq 0$ and $C\geq 0$ and $D$ a block $2\times 2$ diagonal matrix. Thus, even though $A^{-1}$ is generally not in $\MM_\dd$, it is quite close to one: at most $n/2$ entries on the first upper diagonal of $A^{-1}$ may be positive, and additionally, at most two of its row sums may be negative. See Theorem~\ref{Mth-10} and the impression in the right of Figure~\ref{Mfigure22}. Hence, intuitively, the additional conditions that are needed on the $(\nmo)$-facets of an $n$-simplex $S$ in order for $S$ to be nonobtuse, should become {\em weaker} as $n$ increases, as there are relatively less potential obtuse angles to conquer.\\[2mm]
In Section~\ref{Msect-3.3} we prove that if the simplex $S$ with only nonobtuse facets has indeed an obtuse angle, it possesses a vertex Gramian $G$ having a column $g_j$, none of whose off-diagonal entries is minimal in its row in $G$. We will call such a column $g_j$ a {\em blocking column} of $G$.  See Theorem~\ref{Mth-13}. This implies that $G$ cannot satisfy the conditions on its $3\times 3$ submatrices that render it a symmetric ultrametric matrix \cite{VarNab}. An ultrametric matrix is an example of a {\em nonblocking matrix}, a matrix with no blocking columns. This observation immediately leads to a geometric proof that the inverse of a symmetric ultrametric matrix is a weakly diagonally dominant Stieltjes matrix \cite{MaMiSa,NabVar}. To be precise, in \cite{MaMiSa,NabVar} only {\em strictly} ultrametric matrices were considered. Here, we consider the larger class of matrices introduced in \cite{VarNab}. In the remainder of Section~\ref{Msect-3.3} we prove a final property of $n$-simplices $S$ with nonobtuse facets only, which is Theorem~\ref{Mth-4}: all their vertex Gramians are {\em completely positive} \cite{BeSh} with cp-rank $n$. This is equivalent to stating that such $S$ can be isometrically embedded in the nonnegative orthant of $\RR^n$ with any of its vertices placed at the origin, even though $S$ may have obtuse dihedral angles. Note that this concerns the inverse $A$ of the matrix $A^{-1}$ displayed in Figure~\ref{Mfigure22}. In Section~\ref{Msect-3.4} we consider, dually to $n$-simplices with nonobtuse $(\nmo)$-facets, simplices with nonobtuse $2$-facets. Their vertex Gramians $G$ are nonnegative, and as inverses of finite element {\em stiffness matrices} in the finite element method \cite{Bra} to approximate solutions of PDEs they would imply discrete maximum and comparison principles \cite{BrKoKr}. We also state a conjecture about simplices with nonobtuse $k$-facets for other values of $k$, extending classical results by Fiedler \cite{Fie} in combination with the new results in this chapter.\\[2mm] 
Sections~\ref{Msect-4},~\ref{Msect-5},~\ref{Msect-6}, and~\ref{Msect-7} build on the framework provided in Section~\ref{Msect-3} by studying simplices with special nonobtuse $2,3,4$ and general $k$-facets, respectively, that render $S$ nonobtuse. First, in Section~\ref{Msect-4} we study simplices $S$ with those nonobtuse triangular facets for which $S$ nonobtuse. We consider regular simplices \cite{Cox}, then we move via Schl\"afli's \index{Schlafli's orthoschemes}  {\em orthoschemes} \cite{Cox2,Schl}, also called {\em path-simplices} in \cite{BrKoKr,BrKoKrSo} to simplices whose vertices form a {\em finite ultrametric space} with respect to the usual Euclidean metric as studied by Fiedler in \cite{Fie2}. Their corresponding simplicial matrix classes correspond to the simplest types of symmetric ultrametric matrices, such as the {\em type-D matrices} of Markham \cite{Mar}. The vertex Gramians of simplices whose vertices form a finite ultrametric space form only a modest subclass of the strictly ultrametric matrices in the sense of \cite{MaMiSa,NabVar}. To describe the underlying simplices of symmetric ultrametric matrices \cite{VarNab}, in Section~\ref{Msect-4.2} we introduce the new notion of {\em sub-orthocentric tetrahedron}. The clear figures in \cite{HaWe} stood model for our figures in Section~\ref{Msect-5}. We call a tetrahedron sub-orthocentric if each of its vertices projects in the {\em sub-orthocentric set} $\Delta_\ast$ of its opposite nonobtuse triangular facet $\Delta$. This set $\Delta_\ast$ is the union of the line segments between the vertices and the orthocenter of $\Delta$, as depicted in Figure~\ref{Mfigure23}.  If we adopt the convention that the midpoint of a {\em line segment} is its sub-orthocentric set, this turns the sub-orthocentric tetrahedron into a three-dimensional generalization of the equilateral triangle in the sense that simplices whose tetrahedral facets are all sub-orthocentric, are nonobtuse. In Theorem~\ref{Mth-5} we prove that the class $\UU$ of symmetric ultrametric matrices equals the set of all vertex Gramians of simplices having only sub-orthocentric tetrahedral facets and that $\UU\subset\MM_\dd$. In Section~\ref{Msect-5.3} we briefly discuss the subclass of {\em orthogonal simplices}. 
 \\[2mm] 
As a logical next step, in Section~\ref{Msect-6} we discuss the existence of conditions on the $4$-facets of a simplex $S$ that would render $S$ nonobtuse, or in other words, for conditions on the principal $4\times 4$ submatrices of a matrix $A\in\Rnspd$ that would guarantee that $A\in\MM_\dd$. As a first attempt we look for conditions on {\em individual} entries of the $4\times 4$ submatrices of $A$, in a similar way as ultrametricity of a matrix can be read from the entries of its $3\times 3$ submatrices. This corresponds to considering the projections of vertices of the underlying simplex $\SIS(A)$ on its opposite {\em edges}. Although this appeared sufficient in three dimensions, we will prove a negative result in Theorem~\ref{Mth-12}: in four dimensions, the conditions derived from the corresponding concept of {\em blocking columns} are unfortunately so restrictive that only $4\times 4$ submatrices, all whose $3\times 3$ submatrices are ultrametric, remain: only the $4\times 4$ ultrametric submatrices.\\[2mm] 
The geometrical considerations in Section~\ref{Msect-6.3} seem to provide a way out. Instead of considering the projections of each vertex of a $4$-simplex on its {\em six} opposite edges, we will consider their projections on the {\em four} opposite {\em triangular} facets instead. This is discussed in detail in Section~\ref{Msect-6.4}, where we generalize the concept of sub-orthocentric set of a {\em line segment} and a nonobtuse {\em triangle} to a nonobtuse {\em tetrahedron}. In Theorem~\ref{Mth-14} we show that if each vertex of a $4$-simplex $\four$ projects in the sub-orthocentric set of each of its opposite {\em triangular} facets, it also projects in the sub-orthocentric set of its opposite {\em tetrahedral} facet. In other words: a $4$-simplex $\four$ whose tetrahedral facets $\tet$ are all sub-orthocentric, is sub-orthocentric itself.\\[2mm]
This last observation finally gives rise to a general definition of {\em sub-orthocentric set} $S_\ast$ of a nonobtuse $n$-simplex $S$ in Section~\ref{Msect-7}. Correspondingly, we call a simplex with nonobtuse facets {\em sub-orthocentric} if each of its vertices projects in the sub-orthocentric set of its opposite facet. Encouraged by the results for tetrahedra and $4$-simplices, we state the following conjecture.\\[2mm]
{\bf Conjecture. } A simplex with sub-orthocentric facets is sub-orthocentric.\\[2mm]
The main consequence of the validity of this conjecture is, that an $n$-simplex $S$ whose $k$-facets with $2\leq k$ are all sub-orthocentric, is itself sub-orthocentric. As each sub-orthocentric simplices is nonobtuse, this yields exactly the general result that we aimed for in this chapter. Moreover, the sub-orthocentric set $S_\ast$ of an $n$-simplex $S$ seems to be a larger subset of $S$ as $n$ grows: for a triangle $\Delta$ the set $\Delta_\ast$ has dimension one, whereas for tetrahedra $\tet$ that are not orthocentric, the set $\tet_\ast$ has dimension three. This reflects our belief that as $k$ increases, conditions on the $k$-facets of $S$ that render $S$ nonobtuse, become weaker. 
\section{Preliminaries regarding simplices}\label{Msect-2}
The purpose of this section is twofold. In Section~\ref{Msect-2.1} we will familiarize the reader with the linear algebraic description of simplices and their dihedral angles. Most of the material and notation stems from \cite{BrKoKr}, though some of the results date back to work by Fiedler \cite{Fie,Fie2,Fie3,Fie4,FiPt}. In Section~\ref{Msect-2.2} we elaborate on the intimate connection between nonobtuse simplices, Stieltjes matrices and weak diagonal dominance, and show the use of simplicial matrix classes. In a similar manner, the concepts of nonnegativity and weak pointwise diagonal dominance will be related to one another by means of simplices with nonobtuse triangular facets. Finally, we comment on complete positivity \cite{BeSh} and its relation to simplices with a certain shape.
\subsection{Simplices and the dihedral angles}\label{Msect-2.1}
Let $S$ be an $n$-simplex with vertices $v_0,\dots,v_n$. The convex hull of $\{v_0,\dots,v_n\}\setminus\{v_j\}$ is the {\em facet} \index{facet} $F_j$ of $S$ opposite $v$. Let $P_0=(p_1|\dots|p_n)$ be the matrix with columns $p_j=v_j-v_0$. Then, in accordance with Definition~\ref{Mdef-1}, $G_0=P_0^\top P_0$ is a vertex Gramian of $S$ associated with $v_0$. The inverse of $G_0$ equals the matrix $Q_0^\top Q_0$, where $Q_0^\top P_0=I$. Since the $j$-th column $q_j$ of $Q_0$ is orthogonal to all but the $j$-th column of $P_0$, it is a {\em normal vector} to $F_j$, and due to $q_j^\top p_j=1>0$, it points {\em into} $S$. Also due to $q_j^\top p_j=1$, the {\em height} $h_j$ of $p_j$ above $F_j$, the length of the orthogonal projection of $p_j$ on $q_j$, reduces to 
\be h_j = \left\| \frac{q_jq_j^\top}{q_j^\top q_j}p_j\right\| = \frac{1}{\|q_j\|}. \ee
The above definitions do not include an inward pointing normal $q_0$ to the facet $F_0$ {\em opposite} $v_0$. It can be found either by considering a vertex Gramian of $S$ associated with a vertex other than $v_0$, or to interpret each affine functional $q_j^\top\!:S\rightarrow\RR$ as a {\em barycentric coordinate} of $S$, taking value zero on $F_j$ and increasing linearly to value one at $p_j$. Both approaches show that
\be\label{Mq0} q_0=-(q_1+\cdots+q_n)=-Qe , \hdrie\mbox{\rm where }\hdrie  e=e_1+\cdots+e_n\in\RR^n. \ee 
Each pair of facets $F_i$ and $F_j$ of $S$ makes a so-called {\em dihedral} angle \index{dihedral angle} $\alpha_{ij}$, which equals the {\em supplement} $\pi-\gamma_{ij}$ of the angle $\gamma_{ij}$ between their respective inward pointing normals $q_i,q_j$. If $S$ is a triangle, the concept of dihedral angle reduces to the usual concept of angle.
\begin{figure}[h]
\begin{center}  
\begin{tikzpicture}[scale=0.8, every node/.style={scale=0.8}]
\draw (0,0)--(4,0)--(1.2,4)--(0,0);
\draw[dashed] (0,0)--(2,1.4);
\draw[dashed] (1.2,4)--(2,1.4)--(4,0);
\draw[fill] (4,0) circle [radius=0.07];
\draw[fill] (0,0) circle [radius=0.07];
\draw[fill] (1.2,4) circle [radius=0.07];
\draw[fill,gray] (2,1.4) circle [radius=0.07];
\draw[dotted] (1.2,4)--(1.2,0.5);
\draw[dotted] (1.4,1.4)--(4,0);
\draw[fill] (1.2,0.5) circle [radius=0.04];
\draw[fill] (1.4,1.4) circle [radius=0.04];
\node[scale=0.9] at (4.4,0) {$v_1$};
\node[scale=0.9] at (1,4.3) {$v_2$};
\node[scale=0.9] at (-0.4,0) {$v_0$};
\node[gray,scale=0.9] at (2.3,1.7) {$v_3$};
\draw[->] (0,0)--(0,3);
\draw[->] (0,0)--(2.6,-1.4);
\node[scale=0.9] at (0,3.2) {$q_2$};
\node[scale=0.9] at (2.6,-1) {$q_1$};
\node[scale=0.9] at (0.4,2.3) {$F_1$};
\node[scale=0.9] at (2.5,-0.3) {$F_2$};
\node[scale=0.9] at (0.9,1.2) {$h_1$};
\node[scale=0.9] at (2,0.8) {$h_2$};
\begin{scope}[shift={(8,2)},scale=0.8]
\draw (0,0) circle [radius=3];
\draw (-3,0) arc (180:360:3cm and 0.7cm);
\draw[dashed] (-3,0) arc (180:0:3cm and 0.7cm);
\draw[rotate=-55] (-3,0) arc (180:360:3cm and 1cm);
\draw[rotate=-55,dashed] (-3,0) arc (180:0:3cm and 1cm);
\draw[fill] (0,-3) circle [radius=0.07];
\draw[fill] (2,1.5) circle [radius=0.07];
\draw[dotted] (0,0)--(-0.574*3,0.819*3);
\draw[dotted] (-2.27,0.45)--(0,0);
\draw[<->,dotted] (-2.27,0.45) to[out=92,in=230 ] (-0.574*3,0.819*3);
\node[scale=0.9] at (-1.3,-0.4) {$\alpha_{ij}$};
\node[scale=0.9] at (-1.9,0.8) {$\alpha_{ij}$};
\node[scale=0.9] at (3,-3) {\fbox{$\alpha_{ij}+\gamma_{ij}=\pi$}};
\draw[<->] (-1.7,-0.55) to[out=90,in=205 ] (-1.32,0.175);
\draw[dotted] (-0.65,-0.62)--(0,0);
\draw[dotted] (0,0)--(2,1.5);
\draw[->]     (2,1.5)--(3,2.25); 
\draw[dotted] (0,0)--(0,-3);
\draw[->]     (0,-3)--(0,-3.5);
\node[scale=0.9] at (0.3,-3.5) {$q_i$};
\node[scale=0.9] at (3.2,2.5) {$q_j$};
\node[scale=0.9] at (-1,-1.5) {$\mathbb{S}^2$};
\draw[<->] (0.1,-3) arc (273:388:2.3cm and 3cm);
\draw[<->,dotted] (0,-1) arc (270:390:0.8cm and 1cm);
\node[scale=0.9] at (2.5,-0.9) {$\gamma_{ij}$};
\node[scale=0.9] at (0.4,-0.2) {$\gamma_{ij}$};
\draw[fill,gray] (0,0) circle [radius=0.07];
\end{scope}
\end{tikzpicture}    
\end{center} 
\caption{\small{Illustrating normals, heights, and dihedral angle of a tetrahedron.}}
\label{Mfigure24}
\end{figure}
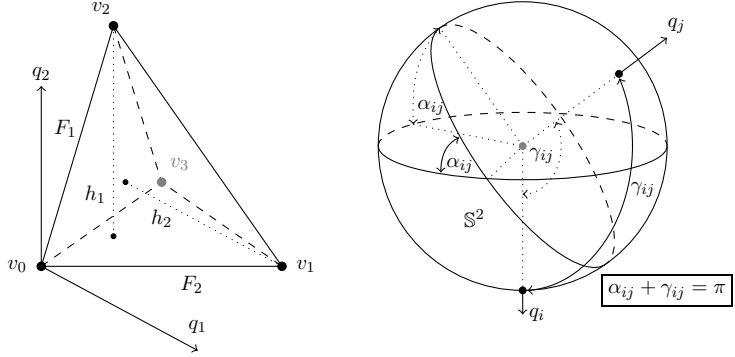 \\[2mm] 
In Figure~\ref{Mfigure24} we illustrate the concepts just introduced for a tetrahedron. Now, a dihedral angle $\alpha_{ij}$ that exceeds $\frac{\pi}{2}$ is called {\em obtuse} \index{obtuse angle}, and this is the case if and only if $q_i^\top q_j$ is positive. If none of its dihedral angles is obtuse, then $S$ is called a {\em nonobtuse} simplex \index{nonobtuse simplex}.
\begin{rem}\label{Mrem-x}{\rm In the context of the finite element method \cite{Bra} to approximate solutions of PDEs, triangulations using only nonobtuse simplices play an important role \cite{BrKoKrSo}, both to ensure {\em superconvergence} \cite{BrKr1} and {\em discrete maximum principles} \cite{BrKoKr}. In fact, the geometric properties of the simplices in the triangulation result in finite element system matrices in $\MM_\dd$.}
\end{rem}
From here on, for given integers $k \leq \ell$ we will write $I_k^\ell=\{k,\dots,\ell\}$.
\subsection{Simplicial matrix classes}\label{Msect-2.2}
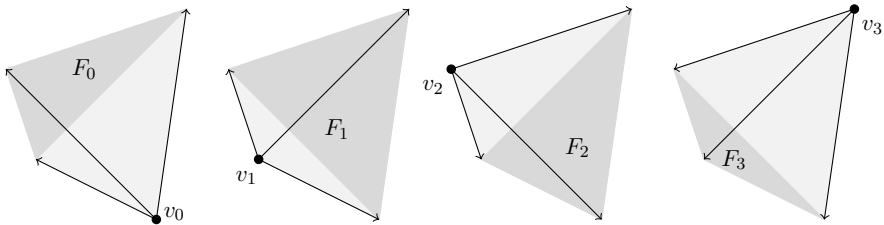
\begin{figure}[h]
\begin{center}  
\begin{tikzpicture}[scale=0.8, every node/.style={scale=0.8}]
\begin{scope}[shift={(0,0)}]
\draw[gray!10!white, fill=gray!10!white] (2.5,0)--(0.5,1)--(0,2.5)--(3,3.5)--cycle;
\draw[gray!30!white, fill=gray!30!white] (0.5,1)--(0,2.5)--(3,3.5)--cycle;
\draw[->] (2.5,0)--(0.5,1);
\draw[->] (2.5,0)--(0,2.5);
\draw[->] (2.5,0)--(3,3.5);
\draw[fill] (2.5,0) circle [radius=0.07];
\node at (2.8,0.1) {$v_0$};
\node at (1.3,2.5) {$F_0$};
\end{scope}
\begin{scope}[shift={(3.7,0)}]
\draw[gray!10!white, fill=gray!10!white] (2.5,0)--(0.5,1)--(0,2.5)--(3,3.5)--cycle;
\draw[gray!30!white, fill=gray!30!white] (2.5,0)--(0,2.5)--(3,3.5)--cycle;
\draw[->] (0.5,1)--(2.5,0);
\draw[->] (0.5,1)--(0,2.5);
\draw[->] (0.5,1)--(3,3.5);
\draw[fill] (0.5,1) circle [radius=0.07];
\node at (0.3,0.7) {$v_1$};
\node at (1.8,1.5) {$F_1$};
\end{scope}
\begin{scope}[shift={(7.4,0)}]
\draw[gray!10!white, fill=gray!10!white] (2.5,0)--(0.5,1)--(0,2.5)--(3,3.5)--cycle;\
\draw[gray!30!white, fill=gray!30!white] (0.5,1)--(2.5,0)--(3,3.5)--cycle;
\draw[->] (0,2.5)--(2.5,0);
\draw[->] (0,2.5)--(0.5,1);
\draw[->] (0,2.5)--(3,3.5);
\draw[fill] (0,2.5) circle [radius=0.07];
\node at (-0.3,2.2) {$v_2$};
\node at (2.1,1.2) {$F_2$};
\end{scope}
\begin{scope}[shift={(11.1,0)}]
\draw[gray!10!white, fill=gray!10!white] (2.5,0)--(0.5,1)--(0,2.5)--(3,3.5)--cycle;
\draw[gray!30!white, fill=gray!30!white] (0.5,1)--(0,2.5)--(2.5,0)--cycle;
\draw[->] (3,3.5)--(2.5,0);
\draw[->] (3,3.5)--(0.5,1);
\draw[->] (3,3.5)--(0,2.5);
\draw[fill] (3,3.5) circle [radius=0.07];
\node at (3.3,3.2) {$v_3$};
\node at (1,1) {$F_3$};
\end{scope}
\end{tikzpicture}    
\end{center} 
\caption{\small{ Proposition~\ref{Mpro-2}: the {\em Stieltjes property} of $G_j^{-1}$ concerns nonobtuse dihedral angles between facets meeting at a vertex $v_j$; the {\em weak diagonal dominance} of $G_j^{-1}$ is about the remaining nonobtuse dihedral angles at the facet $F_j$ opposite $v_j$. Proposition~\ref{Mpro-15}: the {\em nonnegativity} of $G_j$ concerns nonobtuse angles of edges of $S$ meeting at $v_j$; the {\em weak pointwise diagonal dominance} of $G_j$ is about nonobtuse angles between edges meeting at a vertex of~$F_j$.}}
\label{Mfigure25}
\end{figure}

The nonobtusity of a simplex $S$ is of course fully determined by its vertex Gramians \cite{BrKoKr}. Note that the following proposition actually proves the claim in (\ref{Minvm}).
\begin{Pro}\label{Mpro-2} Let $S$ be a simplex and $G$ any vertex Gramian of $S$. Then $S$ is a nonobtuse simplex if and only if $G\in\MM_\dd$, or in other words, if both\\[2mm]
$(1)$ $G^{-1}$ is a Stieltjes matrix: $e_i^\top G^{-1}e_j \leq 0$ for all $i,j\in I_1^n, (i\not=j)$;\\[2mm]
$(2)$ $G^{-1}$ is weakly diagonally dominant: $G^{-1}e\geq 0$, where $e=e_1+\cdots+e_n$.
\end{Pro}  
{\bf Proof. } Let $i,j\in I_1^n$. The entry $e_i^\top G^{-1}e_j$ of $G^{-1}$ is the inner product between inward pointing normals $q_i$ and $q_j$ to the facets $F_i$ and $F_j$ of $S$. Thus $e_i^\top G^{-1}e_j\leq 0$ if and only if $F_i$ and $F_j$ make a nonobtuse dihedral angle. Due to (\ref{Mq0}), the $i$-th entry of $G^{-1}e$ equals the inner product between the {\em outward} pointing normal $-q_0$ to the facet $F_0$ of $S$, and the {\em inward} normal $q_i$ to $F_i$. It is nonnegative if and only if the dihedral angle between $F_0$ and $F_i$ is nonobtuse.~\hfill $\Box$
\begin{rem}\label{Mrem-1}{\rm  Proposition~\ref{Mpro-2} expresses a subtle geometrical duality between the property (1) of $G^{-1}$ being a Stieltjes matrix, which relates to the $\half n(\nmo)$ dihedral angles between $n$ facets of $S$ that {\em meet at a given vertex} $v$ of $S$, and the weak diagonal dominance (2) of $G^{-1}$, which relates to the $n$ dihedral angles of $S$ {\em at the facet opposite} $v$. See also Figure~\ref{Mfigure25}}.
\end{rem}
Consequently, the relation between a nonobtuse simplex and its vertex Gramians is as follows.
\begin{Th}\label{Mth-1} Let $S$ be an $n$-simplex. The following five statements are equivalent:\\[1mm]
$\bullet$ $S$ is a nonobtuse simplex;\\[1mm]
$\bullet$ $S$ has a vertex Gramian whose inverse is a weakly diagonally dominant Stieltjes;\\[1mm]
$\bullet$ the inverse of each vertex Gramian of $S$ is a Stieltjes matrix;\\[1mm]
$\bullet$ the inverse of each vertex Gramian of $S$ is weakly diagonally dominant;\\[1mm]
$\bullet$ each vertex $v$ of $S$ projects onto its opposite facet $F$.
\end{Th}    
{\bf Proof. } The equivalence of the first four items can be concluded from the proof of Proposition~\ref{Mpro-2} above and Remark~\ref{Mrem-1}. Since the fifth item expresses that the dihedral angle between $F$ and another facet of $S$ is nonobtuse, also the fifth equivalence follows.~\hfill $\Box$\\[2mm]
See Figure~\ref{Mfigure25} for an illustration of Theorem~\ref{Mth-1} for $n=3$. The Stieltjes property of $G_j^{-1}$ relates to the three nonobtuse dihedral angles between the three facets {\em meeting} at the vertex $v_j$, whereas the weak diagonal dominance of $G_j^{-1}$ relates to the three nonobtuse dihedral angles that the facet $F_j$ {\em opposite} $v_j$ (depicted in darker gray) makes with the facets meeting at $v_j$. The dihedral angles of the tetrahedron $\tet$ are all nonobtuse if {\em one} of both properties is shared by {\em all} vertex Gramians of $\tet$. Or, to be more precise, by {\em all but one} of them.\\[2mm]
In (\ref{Mdnn}) we claimed that the vertex Gramians of simplices with {\em nonobtuse triangular facets} form a simplicial matrix class equal to $\DD_\dd$. It is proved in the next counterpart of Proposition~\ref{Mpro-2}.
\begin{Pro}\label{Mpro-15} Let $S$ be an $n$-simplex and $G=(g_{ij})$ any vertex Gramian of $S$. Then all triangular facets of $S$ are nonobtuse if and only if $G\in\DD_\dd$, or in other words, if both\\[2mm]
(1) $G$ is nonsingular doubly nonnegative: $G\in\Rnspd$ and $g_{ij} \geq 0$ for all $i,j\in I_1^n$;\\[2mm]
(2) $G$ is pointwise weakly diagonally dominant: $g_{ii}\geq g_{ij}$ for all $i,j\in I_1^n$. 
\end{Pro} 
{\bf Proof. } Each entry $g_{ij}$ of the vertex Gramian $G$ associated with a vertex $v_\ell$ of $S$ equals the inner product between vectors constituting edges of $S$ that meet at $v_\ell$. The difference $g_{ii}-g_{ij}$ 
for $i \neq j$ equals the inner product between two edges meeting at $v_i$.~\hfill $\Box$\\[2mm]
The counterpart of Remark~\ref{Mrem-1} in this context is that the {\em nonnegativity} of $G$ concerns the triangular angles meeting at a given vertex $v$ of $S$, whereas the {\em pointwise weak diagonal dominance} relates to angles in triangular facets oposite $v$. See again Figure~\ref{Mfigure25}. 
Now, the following theorem is the counterpart of Theorem~\ref{Mth-1}.
\begin{Th}\label{Mth-1b} Let $S$ be an $n$-simplex. The following five statements are equivalent:\\[1mm]
$\bullet$ $S$ is a simplex with nonobtuse triangular facets only;\\[1mm]
$\bullet$ $S$ has a nonnegative pointwise weakly diagonally dominant vertex Gramian;\\[1mm]
$\bullet$ each vertex Gramian of $S$ is nonnegative;\\[1mm]
$\bullet$ each vertex Gramian of $S$ is pointwise weakly diagonally dominant;\\[1mm]
$\bullet$ each vertex $v$ of $S$ projects into each opposite edge.
\end{Th}    
{\bf Proof. } The equivalence of the first four items can be concluded from the proof of Proposition~\ref{Mpro-15} above and Remark~\ref{Mrem-1}. The fifth item is trivial.~\hfill $\Box$\\[2mm]
The relevance of the observations in Theorem~\ref{Mth-1} and~\ref{Mth-1b} become more apparent in the context of the simplicial matrix classes, as defined in Definition~\ref{Msmc}.
\begin{Co}\label{Mcor-5} Let $\CC$ be a simplicial matrix class. Then equivalent are:\\[2mm]
$\bullet$ $Ax=e$ has a solution $x\geq 0$ for all $A\in\CC$;\\[2mm]
$\bullet$ $A^{-1}$ is a Stieltjes matrix for all $A\in\CC$.\\[2mm]
Moreover, also the following are equivalent:\\[2mm]
$\bullet$ $A\geq 0$ for all $A\in\CC$;\\[2mm]
$\bullet$ $A$ is pointwise weakly diagonally dominant for all $A\in\CC$.
\end{Co} 
{\bf Proof. } Let $A\in\CC$ be given and $S=\SIS(A)$ the simplex underlying $A$. Since $\CC$ is a simplicial matrix class, all other vertex Gramians of $S$ are also in $\CC$. The assumption that $Ax=e$ has solution $x\geq 0$ for all $A\in\CC$ shows that their inverses are all weakly diagonally dominant. Theorem~\ref{Mth-1} yields that this is equivalent with all their inverses being Stieltjes matrices.~\hfill $\Box$\\[2mm]
To motivate the use of Corollary~\ref{Mcor-5}, we mention again that the set $\UU$ of symmetric ultrametric matrices is a simplicial matrix class. This will be proved in detail in Section~\ref{Msect-5}. Therefore, to prove that their inverses are diagonally dominant Stieltjes matrices \cite{MaMiSa,VarNab} it is sufficient to show {\em either} the Stieltjes property of an arbitrary ultrametric matrix, {\em or} the existence of a nonnegative solution $x$ of $Ax=e$. The latter was, in fact, already proved as the very first lemma in \cite{MaMiSa}, where $x$ is called an {\em equilibrium potential} \index{equilibrium potential}.\\[2mm]
Although the simplicial matrix classes $\MM_\dd$ and $\DD_\dd$ consist of well-known and well-studied matrix types, far less seems to be know about the classes $\VV_k^n$ from (\ref{Mchain}) for $2<k<n$.
\subsection{Nonobtusity by facets \& complete positivity}
The fact that the {\em converse} of the following proposition does {\em not} hold lies at the basis of the symmetric {\em inverse M-matrix problem}. We will discuss this aspect in detail in Section~\ref{Msect-3}.
\begin{Pro}[\cite{Fie}]\label{Mpro-3} Each facet of a nonobtuse $n$-simplex is nonobtuse.
\end{Pro}
The nonobtusity of all $(\nmo)$-facets of an $n$-simplex $S$ is therefore a {\em necessary condition} for $S$ itself to be nonobtuse. As such, it will be of interest to study properties of $n$-simplices $S$ having only nonobtuse $(\nmo)$-facets. This will be done in Section~\ref{Msect-3}. As far as we know, such simplices have not been studied before.\\[2mm]
Before that, we will prove that each vertex Gramian $G$ of a nonobtuse $n$-simplex is {\em completely positive with cp-rank} \index{completely positive with cp-rank} $n$ \cite{BeSh}, which means that $G$ is the Gramian $P^\top P$ of a {\em nonnegative} $n\times n$ matrix $P$. This is a much stronger property than the nonnegativity of $G$ itself. 
\begin{Co}\label{Mcor-2} Let $G$ be a vertex Gramian of a nonobtuse simplex $S$. Then $G$ is completely positive with cp-rank $n$. In other words, $G=P^\top P$ for a nonnegative $n\times n$ matrix $P$.
\end{Co}
{\bf Proof. } The fact that the factor $P$ of the Gramian $G$ is nonnegative is equivalent with the posibility to rigidly transform $S$ into the nonnegative orthant of $\RR^n$ with the vertex of $S$ associated with $G$ located at the origin. Now, assume inductively that a facet $F$ of $S$, which is nonobtuse due to Proposition~\ref{Mpro-3}, lies in the nonnegative orthant of $\RR^{n-1}$ with one of its vertices at the origin. Since by Theorem~\ref{Mth-1}, the remaining vertex of $S$ projects into $F$, it has at least $n-1$ nonnegative coordinates. In case it has one negative coordinate, reflection of $S$ in its facet $F$ places it in the nonnegative orthant of $\RR^n$.~\hfill $\Box$\\[2mm]
The inductive structure of the proof shows that in fact $G=U^\top U$ for a nonnegative {\em upper triangular} $n\times n$ matrix and constitutes a geometric supplement to Theorem~2.27 in \cite{BeSh}.
\begin{rem}{\rm The class $\PP$ of all completely positive $n\times n$ matrices with cp-rank $n$ is {\em not} a simplicial matrix class. This can be seen from the simple example in 
Figure~\ref{Mfigure21}, in which $G_u=P_u^\top P_u$ with $P_u\geq 0$, whereas trivially, $G_v$ cannot be factorized at $G_v=R^\top R$ with $R\geq 0$. Now, write $\PP_\dd$ for the largest simplicial subclass of $\PP$.
In Section~\ref{Msect-3} we will prove that, surprisingly, for all $n\geq 3$, each vertex Gramian of an $n$-simplex with only nonobtuse facets is completely positive with cp-rank $n$, in spite of the fact that $S$ may have obtuse dihedral angles. For $n=2$, this is trivially false (see again Figure~\ref{Mfigure21}). In terms of the chain (\ref{Mchain}) of simplicial matrix classes, this shows that for $n\geq 3$,
\be \MM_\dd = \VV_n^n \subset \VV_n^{n-1} \subset \PP_\dd \subset \DD_\dd = \VV_n^2. \ee
It is well known that $\PP=\DD$ for $n\leq 4$ {\em only} \cite{BeSh}. It is therefore an interesting question what are the optimal values of $k$ and $\ell$ such that $V_n^k\subset \PP_\dd\subset \VV_n^\ell$, as well as their dependence on $n$.} 
\end{rem} 
 \section{Simplices whose $k$-facets are all nonobtuse}\label{Msect-3} 
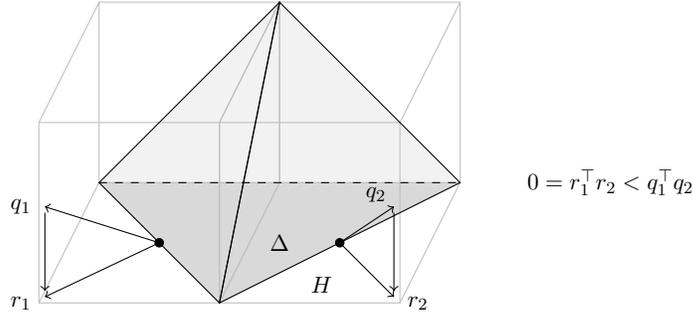
\begin{figure}[h]
\begin{center}  
\begin{tikzpicture}[scale=0.8, every node/.style={scale=0.8}]
\draw[gray!10!white, fill=gray!10!white] (3,0)--(1,2)--(4,5)--(7,2)--cycle;
\draw[gray!30!white, fill=gray!30!white] (3,0)--(1,2)--(7,2)--cycle;
\draw[gray!50!white] (0,0)--(3,0)--(3,3)--(0,3)--cycle;
\draw[gray!50!white] (1,2)--(4,2)--(4,5)--(1,5)--cycle;
\draw[gray!50!white] (0,0)--(1,2);
\draw[gray!50!white] (3,0)--(4,2);
\draw[gray!50!white] (3,3)--(4,5);
\draw[gray!50!white] (0,3)--(1,5);
\begin{scope}[shift={(3,0)}]
\draw[gray!50!white] (0,0)--(3,0)--(3,3)--(0,3)--cycle;
\draw[gray!50!white] (1,2)--(4,2)--(4,5)--(1,5)--cycle;
\draw[gray!50!white] (0,0)--(1,2);
\draw[gray!50!white] (3,0)--(4,2);
\draw[gray!50!white] (3,3)--(4,5);
\draw[gray!50!white] (0,3)--(1,5);
\end{scope}
\draw (3,0)--(1,2)--(4,5)--cycle;
\draw (4,5)--(7,2)--(3,0)--cycle;
\draw[dashed] (1,2)--(7,2);
\draw[->] (2,1)--(0.1,0.1);
\draw[->] (2,1)--(0.1,1.6);
\node at (-0.3,1.6) {$q_1$};
\node at (-0.3,0) {$r_1$};
\draw[fill=black] (2,1) circle [radius=0.07];
\draw[->] (5,1)--(5.9,0.1);
\draw[->] (5,1)--(5.9,1.6);
\node at (5.6,1.8) {$q_2$};
\node at (6.3,0) {$r_2$};
\node at (9.5,2) {$0 = r_1^\top r_2 < q_1^\top q_2$};
\node at (4,1) {$\Delta$};
\node at (4.7,0.3) {$H$};
\draw[fill=black] (5,1) circle [radius=0.07];
\draw[->] (0.1,1.5)--(0.1,0.2);
\draw[->] (5.9,1.5)--(5.9,0.2);
\end{tikzpicture}    
\end{center} 
\caption{\small{All facets of the tetrahedron $\tet$ are nonobtuse, but $\tet$ has an obtuse dihedral angle. The projections $r_1,r_2$ of normals $q_1,q_2$ of $\tet$ are normals to facets of the bottom facet~$\Delta$.}}
\label{Mfigure26}
\end{figure}

Recall from Proposition~\ref{Mpro-3} that each $k$-facet of a nonobtuse simplex, seen as a $k$-simplex in its ambient $k$-space, is again nonobtuse \cite{Fie}. Reason for this is, that the normals to {\em the facets of a facet} $F$ of an $n$-simplex $S$ are the projections onto the affine hyperspace $H$ containing $F$ of the normals to the $n$ facets of $S$ other than $F$. If these $n$ normals lie, together with $S$, in the same half-space separated by $H$, projecting them on $H$ decreases their inner product. See Figure~\ref{Mfigure26}, where we depict two normals $q_1,q_2$ to facets of a tetrahedron $\tet$ and their projections $r_1,r_2$ onto the plane $H$ containing a facet $\Delta$ of $\tet$. Then $r_1,r_2$ are normals to facets of $\Delta$, and their zero inner product is smaller than the inner product between $q_1$ and $q_2$ themselves. Thus $\tet$ has an obtuse dihedral angle, even though $\tet$ consists of two tetrahedral corners of adjacent cubes glued together: as such, $\tet$ has two {\em right} and two {\em equilateral} triangular faces. Due to Theorem~\ref{Mth-1}, this means that $\tet$ has a vertex $v$ that does not project onto its opposite facet $\Delta$, even though $v$ projects onto {\em all facets} of $\Delta$, which are the three edges of $\Delta$.
\\[2mm]
Given an $n$-simplex $S$ with nonobtuse $k$-facets, it is both a logical and an interesting question which {\em additional} properties these $k$-facets should have for $S$ itself to be nonobtuse, and possibly to have the same additional properties. One may start with $k=2$ and observe (see also Section~\ref{Msect-4}) that having only {\em equilateral} triangular faces can be rightfully considered as the {\em strongest} such property, as a simplex with only  equilateral triangular facets is regular and hence nonobtuse. The challenge lies in finding weaker conditions. As argued in the Introduction, this {\em geometric} question is intimately related to the {\em algebraic} question which symmetric nonnegative matrices are {\em inverse M-matrices} \cite{Joh,JoSm,Wil}. As nonobtusity of the facets of an $n$-simplex is a necessary condition for $S$ to be nonobtuse, in this section we will first study general properties of simplices with nonobtuse facets.
\subsection{The dual hull of a nonobtuse simplex and its properties} \label{Msect-3.1}
Let $S\subset \RR^n$ be an $n$-simplex with vertices $v_0,\dots,v_n$ and respective opposite facets $F_0,\dots,F_n$ and inward pointing normals $q_0,\dots,q_n$. For each $j\in I_0^n$, let $H_j$ be the affine hyperplane \index{affine hyperplane} with $F_j\subset H_j$, and $H_j^+$ the closed half-space \index{closed half-space} separated by $H_j$ with $S\subset H_j^+$. Then obviously
\be S=\bigcap_{j=0}^n H_j^+.\ee
Also the other regions defined by intersections of some of these half-spaces are of interest.
\begin{Def}\label{MDef-2}{\rm Given the half-spaces $H_0^+,\dots,H_n^+$ in $\RR^n$ intersecting as an $n$-simplex $S$, let
\be \nu_S:\RR^n\rightarrow \mathbb{N}: \hdrie x \mapsto \chi(H_0^+)+\cdots+\chi(H_n^+), \ee
where $\chi(X)$ is the characteristic function \index{characteristic function} of the set $X\subset \RR^n$}.
\end{Def} 
Observe that $\nu_S(x)$ simply counts the number of half-spaces $H_j^+$ in which $x$ lies. It can also be used to count the number of nonobtuse angles of an $(\npo)$-simplex having $S$ as a facet.
\begin{rem}\label{MRem-7}{\rm Given the simplex $S$, each $x\in\RR^n$ lies in $1\leq \nu_S(x)\leq n+1$ of the half-spaces $H_j^+$. Suppose now that $S$ is supplemented with a vertex $v$ as to form an $(\npo)$-simplex $\hat{S}$. Then the number of nonobtuse angles that the facet $S$ of $\hat{S}$ makes with the other facets equals $\nu_S(y)$, where $y$ is the projection of $v$ onto the ambient space of the facet $S$ of $\hat{S}$. Compare this with the fifth bullet in Theorem~\ref{Mth-1}, implicitly using the case $\nu_F(y)=n$ for a facet $F$}.
\end{rem} 
For each facet $F_j$ of $S$, let $C_j$ be the cylinder \index{cylinder} defined by $F_j$ and $q_j$.  Thus, $C_j$ consists of precisely those points that project orthogonally onto $F_j$, and consequently,
\be \label{Meq-15-1} S^\ast = \bigcap_{j=0}^n C_j \ee
is the set of points that project onto {\em each} of the facets of $S$. Observe that $S^\ast$ is the intersection of convex sets, and thus $S^\ast$ itself is convex.
\begin{Pro} The simplex $S$ is nonobtuse if and only if $S\subset S^\ast$.
\end{Pro}
{\bf Proof. } Due to Theorem~\ref{Mth-1}, each vertex $v_j$ projects onto its opposite facet $F_j$ if and only if $S$ is nonobtuse. This is both necessary and sufficient for $S$ to be contained in $C_j$ for all $j$.~\hfill $\Box$\\[2mm]
As the above proposition shows that if $S$ has an obtuse angle, $S$ is not any more contained in the set $S^\ast$, the  following definition is now well-motivated.
\begin{Def}{\rm  Whenever $S$ is nonobtuse, we will call the set $S^\ast$ in (\ref{Meq-15-1}) the {\em dual hull} of $S$} \index{dual hull}.
\end{Def}
\begin{rem}{\rm The name {\em dual hull} stems from the fact that $S^\ast$ can equivalently be defined as the intersection of the {\em dual cones} \cite{BoVa} of the simplicial cones emerging from each vertex of $S$. To be explicit, let
\be\label{Mcones} K_j = \bigcap_{i\not=j} H_i^+. \ee
Then $K_j$ is the pointed simplicial cone with origin $v_j$ spanned by the vectors $v_i-v_j$, $i\not=j$. The {\em dual cone} \index{dual cone} $K_j^\ast$ of $K_j$ is then defined as
\be K_j^\ast = \{ x\in\RR^n \sth x^\top y \geq 0 \hdrie\mbox{\rm for all } y\in K_j \}. \ee
Taking intersections of dual cones originating from different positions in $\RR^n$ is technically cumbersome, and therefore we prefer to work with the cyclinders}.
\end{rem}
To study the dual hull $S^\ast$ of a nonobtuse simplex $S$ in more detail, we split each cylinder $C_j$ in two closed half-cylinders as follows,
\be\label{Meq-15-2} C_j = C_j^+ \cup C_j^- \hdrie F_j = C_j^+ \cap C_j^-, \hdrie S\subset C_j^+.\ee
We will refer to $C_j^+$ as the {\em interior} cyclinder \index{interior cylinder} and to $C_j^-$ as the {\em exterior} cylinder \index{exterior cylinder} of $F_j$. 
\begin{Pro}\label{MProp-4} Any pair of distinct exterior cylinders $C_i^-,C_j^-$ have disjoint interior.
\end{Pro}
{\bf Proof. } Any point outside a simplex projects in the interior of at most one of its facets.~\hfill $\Box$\\[3mm]
Application of De Morgan's law to (\ref{Meq-15-1}) with each $C_j$ split up in two parts as in (\ref{Meq-15-2}) results in the union of $2^{n+1}$ intersections of $n+1$ sets. Due to Proposition~\ref{MProp-4}, only those intersections containing at most one exterior cylinder may be nonempty. Therefore, since the intersection of all interior cylinders equals $S$ itself, we have that
\be S^\ast = \bigcap_{j=0}^n (C_j^+\cup C_j^-) = S \cup S^\ast_0 \cup \cdots \cup S^\ast_n,\ee
where    
\be\label{MDandE} S_j^\ast = C_j^- \cap \bigcap_{i\not=j} C_i^+. \ee 
Figure~\ref{Mfigure26} illustrates the sets $\Delta_0^\ast,\Delta_1^\ast,\Delta_2^\ast$ for a nonobtuse triangle $\Delta$.
\begin{Pro}\label{Mpro-6} Let $S$ be a nonobtuse $n$-simplex. The $n+2$ sets $S,S_0^\ast,\dots,S_n^\ast$ have pairwise disjoint interiors. In particular, $S\cap S_j^\ast=F_j$ for all $j\in I_0^n$. 
\end{Pro}
{\bf Proof. } This follows because $S_j^\ast\subset C_j^-$ for all $j\in I_0^n$ together with Proposition~\ref{MProp-4}, whereas the interior of $S$ does not lie in any exterior cylinder.~\hfill $\Box$
\subsection{$n$-Simplices with all $(\nmo)$-facets nonobtuse}\label{Msect-3.2}
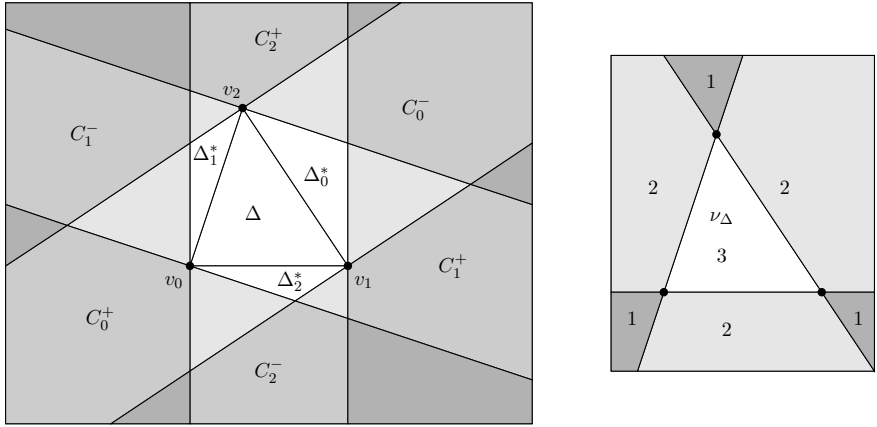
\begin{figure}[h]
\begin{center}  
\begin{tikzpicture}[scale=0.7, every node/.style={scale=0.7}]
\draw (0,0)--(3,0)--(1,3)--cycle;
\draw (0,-3)--(0,5);
\draw (3,-3)--(3,5);  
\draw (-7/2,7/6)--(13/2,-13/6);
\draw (-7/2,9/2)--(13/2,7/6);
\draw (-7/2,0)--(4,5);
\draw (-3/2,-3)--(13/2,7/3);
\draw (-7/2,-3)--(13/2,-3)--(13/2,5)--(-7/2,5)--cycle;
\draw[fill=gray!0!white] (0,0)--(2,-2/3)--(3,0)--cycle;
\draw[fill=gray!0!white] (0,0)--(0,7/3)--(1,3)--cycle;
\draw[fill=gray!0!white] (3,0)--(3,7/3)--(1,3)--cycle;
\draw[fill=gray!20!white] (3,7/3)--(3,13/3)--(1,3)--cycle;
\draw[fill=gray!20!white] (3,7/3)--(3,0)--(16/3,9.3/6)--cycle;
\draw[fill=gray!20!white] (3,0)--(3,-1)--(2,-2/3)--cycle;
\draw[fill=gray!20!white] (0,0)--(2,-2/3)--(0,-2)--cycle;
\draw[fill=gray!20!white] (0,0)--(0,7/3)--(-7/3,0.77)--cycle;
\draw[fill=gray!20!white] (0,7/3)--(1,3)--(0,10/3)--cycle;
\draw[fill=gray!40!white] (3,7/3)--(3,13/3)--(4,5)--(13/2,5)--(13/2,7/3)--(16/3,9.3/6)--cycle;
\draw[fill=gray!40!white] (0,0)--(0,-2)--(-3/2,-3)--(-7/2,-3)--(-7/2,0)--(-7/3,0.77)--cycle;
\draw[fill=gray!40!white] (0,-2)--(0,-3)--(3,-3)--(3,-1)--(2,-2/3)--cycle;
\draw[fill=gray!40!white] (0,10/3)--(0,5)--(3,5)--(3,13/3)--(1,3)--cycle;
\draw[fill=gray!40!white] (-7/2,7/6)--(-7/2,9/2)--(0,10/3)--(0,7/3)--(-7/3,0.77)--cycle;
\draw[fill=gray!40!white] (3,0)--(3,-1)--(13/2,-13/6)--(13/2,7/6)--(16/3,9.3/6)--cycle;
\draw[fill=gray!60!white] (0,10/3)--(0,5)--(-7/2,5)--(-7/2,9/2)--cycle;
\draw[fill=gray!60!white] (-7/2,0)--(-7/2,7/6)--(-7/3,0.77)--cycle;
\draw[fill=gray!60!white] (0,-2)--(0,-3)--(-3/2,-3)--cycle;
\draw[fill=gray!60!white] (3,13/3)--(3,5)--(4,5)--cycle;
\draw[fill=gray!60!white] (13/2,7/6)--(13/2,7/3)--(16/3,9.3/6)--cycle;
\draw[fill=gray!60!white] (3,-1)--(3,-3)--(13/2,-3)--(13/2,-13/6)--cycle;
\draw[fill=black] (0,0) circle [radius=0.07];
\draw[fill=black] (3,0) circle [radius=0.07];
\draw[fill=black] (1,3) circle [radius=0.07];
\node at (-0.3,-0.3) {$v_0$};
\node at (3.3,-0.3) {$v_1$};
\node at (0.8,3.3) {$v_2$};
\node at (1.2,1) {$\Delta$};
\node at (4.3,3) {$C_0^-$};
\node at (-1.7,-1) {$C_0^+$};
 \node at (5,0) {$C_1^+$};
\node at (-2,2.5) {$C_1^-$};
\node at (1.5,-2) {$C_2^-$};
\node at (1.5,4.3) {$C_2^+$};
\node at (2.4,1.7) {$\Delta_0^\ast$}; 
\node at (0.3,2.1) {$\Delta_1^\ast$}; 
\node at (1.9,-0.3) {$\Delta_2^\ast$}; 
\begin{scope}[shift={(9,-0.5)}]
\draw (-1,-1.5)--(4,-1.5)--(4,4.5)--(-1,4.5)--cycle;
\draw (0,4.5)--(4,-1.5);
\draw (-1,0)--(4,0);
\draw (-0.5,-1.5)--(1.5,4.5);
\draw[fill=gray!20!white] (0,0)--(-1,0)--(-1,4.5)--(0,4.5)--(1,3)--cycle;
\draw[fill=gray!20!white] (3,0)--(1,3)--(1.5,4.5)--(4,4.5)--(4,0)--cycle;
\draw[fill=gray!20!white] (0,0)--(-0.5,-1.5)--(4,-1.5)--(3,0)--cycle;
\draw[fill=gray!60!white] (0,0)--(-0.5,-1.5)--(-1,-1.5)--(-1,0)--cycle;
\draw[fill=gray!60!white] (3,0)--(4,0)--(4,-1.5)--cycle;
\draw[fill=gray!60!white] (0,4.5)--(1.5,4.5)--(1,3)--cycle;
\draw[fill=black] (0,0) circle [radius=0.07];
\draw[fill=black] (3,0) circle [radius=0.07];
\draw[fill=black] (1,3) circle [radius=0.07];
\node at (1.1,1.4) {$\nu_\Delta$};
\node at (1.1,0.7) {$3$};
\node at (-0.2,2) {$2$};
\node at (2.3,2) {$2$};
\node at (1.2,-0.7) {$2$};
\node at (-0.6,-0.5) {$1$};
\node at (3.7,-0.5) {$1$};
\node at (0.9,4) {$1$};
\end{scope}
\end{tikzpicture}    
\end{center} 
\caption{\small{Left: The {\em cylinders} $C_0,C_1,C_2$ associated with the edges of a nonobtuse triangle $\Delta$ divide space into regions belonging to $0,1,2$ or $3$ of the cylinders. Right: the {\em half-spaces} $H_0^+,H_1^+,H_2^+$ divide space into regions belonging to $\nu_\Delta(x)\in\{1,2,3\}$ of the half-spaces.}}
\label{Mfigure27}
\end{figure}

The dual hull $S^\ast$ of a nonobtuse simplex $S$ can be used to characterize simplices with nonobtuse facets. We will first illustrate this by studying tetrahedra with nonobtuse triangular facets. Consider the left picture in Figure~\ref{Mfigure27}. The white hexagon is the dual hull $\Delta^\ast$ of the nonobtuse triangle $\Delta$. It consists of $\Delta$ and the triangles $\Delta_j^\ast$ defined in (\ref{MDandE}), each sharing an edge with $\Delta$. By construction, the angles of $\Delta$ and $\Delta_j^\ast$ at their shared edge add pairwise, at each vertex, to a right angle. In the right picture, space is divided into regions belonging to $3,2$ or $1$ of the half-spaces $H_0^+,H_1^+,H_2^+$; in other words, we show the function $\nu_\Delta$ from Definition~\ref{MDef-2}. Superimposing the left picture on top of the right, we see that each triangle $\Delta_j^\ast$ lies in all half-spaces but $H_j^+$. This yields that $\nu_\Delta(x)=2$ for all $x\in E_j(\Delta)$.
\\[2mm]
The latter observation shows, in combination with Remark~\ref{MRem-7}, that if $\Delta$ is supplemented with a fourth vertex $v$ as to form a tetrahedron $\tet$, then if $v$ projects onto $\Delta^\ast$, the facet $\Delta$ of $\tet$ makes {\em at most one} obtuse dihedral angle with the other three facets of $\tet$. As a consequence, assuming that all triangular facets of $\tet$ are nonobtuse, each of them makes at most one obtuse angle with an other, and thus, $\tet$ has at most {\em two} obtuse dihedral angles, whereas an arbitrary tetrahedron may have {\em three} obtuse angles. To conclude, choose the origin of $\RR^2$ at any of the vertices of $\Delta$. Consider the union of $\Delta$ with any of the three triangles $\Delta_j^\ast$. Obviously, this union can be isometrically embedded into the nonnegative quadrant of $\RR^2$. Since $v$ projects in at most one of the sets $\Delta_j^\ast$, this surprisingly shows that all vertex Gramians of $\tet$ are {\em completely positive} with cp-rank $3$, even though $\tet$ may have two obtuse dihedral angles.\\[2mm]
We will now proceed to prove similar conclusions for higher dimensional simplices.
\begin{Le}\label{Mlem-6} Let $S$ be a nonobtuse $n$-simplex. Then:\\[2mm]
$\bullet$ for each $j\in I_0^n$, the set $S_j^\ast$ is an $n$-simplex sharing the facet $F_j$ with $S$;\\[2mm]
$\bullet$ $S_j^\ast$ lies in the simplicial cone $K_j$ defined in (\ref{Mcones});\\[2mm]
$\bullet$ the $n$ pairs of dihedral angles of $S$ and $S_j^\ast$ at each facet of $F_j$ add up to right.
\end{Le} 
{\bf Proof. } Assuming that $S$ is nonobtuse, $F_j$ makes a nonobtuse angle with each facet of $K_j$ and thus $S_j^\ast\subset C_j^- \subset K_j$. By definition, each facet of $S_j^\ast$ is orthogonal to a facet of $S$ and $F_j$ bisects this angle. ~\hfill  $\Box$\\[2mm]
An immediate corollary of the above is as follows.
\begin{Co}\label{Mcor-4} Let $S$ be a nonobtuse $n$-simplex with dual hull $S^\ast$. Then $\nu_S(x)=n+1$ for all $x\in S$ and 
$\nu_S(x)=n$ for all $x\in S^\ast\setminus S$.
\end{Co} 
{\bf Proof. } By Lemma~\ref{Mlem-6}, $S_j^\ast\subset K_j$. Now, due to (\ref{Mcones}) we have $\nu_S(x)=n$ for all $x\in K_j\setminus S$ and trivially $\nu_S(x)=n+1$ for all $x\in S$.~\hfill $\Box$
 \begin{Th}\label{Mth-2} Let $S$ be an $n$-simplex. Then the following two statements are equivalent:\\[2mm]
$\bullet$ each facet of $S$ is a nonobtuse $(\nmo)$-simplex;\\[2mm]
$\bullet$ each vertex $v$ of $S$ projects onto the dual hull $F^\ast$ of its opposite facet $F$;\\[2mm]
and both of them imply:\\[2mm]
$\bullet$  each facet of $S$ makes at most one obtuse dihedral angle with the remaining $n$ facets.\\[2mm]
Thus, an $n$-simplex $S$ with nonobtuse facets has at most $\lfloor\half(\npo)\rfloor$ obtuse dihedral angles.\\[2mm]
See also the table in Figure~\ref{Mfigure22} in introduction of this chapter for an illustration.
\end{Th}   
{\bf Proof. } Let $v$ be a vertex of $S$ opposite the facet $F$ of $S$. If $v$ projects onto $F^\ast$, it also projects onto each of the facets of $F$. If this happens for all vertices of $S$, then all facets of $S$ are nonobtuse. Conversely, if a vertex $v$ of $S$ does not project onto $F^\ast$, it also does not project into at least one facet of $F$. This proves the equivalence of the first two items. Suppose next that $v$ projects into the dual hull $F^\ast$ of its opposite facet $F$. Then Corollary~\ref{Mcor-4} proves that $F$ makes at least $\nmo$ nonobtuse angles.  Since any obtuse angle involves normals to a pair of facets, this proves also that $S$ cannot have more than  $\lfloor\half(\npo)\rfloor$ obtuse dihedral angles.~\hfill $\Box$\\[3mm] 
The following theorem rephrases some of the above results in terms of linear algebra.
 \begin{Th}\label{Mth-10} Let $S$ be an $n$-simplex with nonobtuse facets, and $G$ a vertex Gramian of $S$. Consider the solutions $x$ and $y_j$ of
\be Gx = e \und Gy_j = e_j. \ee
Then $x$ has at most one negative entry and each $y_j$ has at most two positive entries. Thus, the inverse of $G$ has the form
\be\label{Mx-1} G^{-1} \overset{\pi}{\sim} D-C \ee
with $C\geq 0$, and where $D\geq 0$ is block-diagonal with blocks of size $1\times 1$ and $2\times 2$ only.\\[2mm]
See also Figure~\ref{Mfigure22} in introduction of this chapter for an illustration.
\end{Th} 
{\bf Proof. } Combine Theorem~\ref{Mth-2} with Theorem~\ref{Mth-1}.~\hfill $\Box$\\[2mm]
The standard proof that the inverse $A$ of an M-matrix $A^{-1}$ is a nonnegative matrix, is to split the matrix $A^{-1}$ as $A^{-1}=D-C=D(I-D^{-1}C)$ with $D\geq 0$ diagonal and invertible, and $C\geq 0$, to prove that the spectral radius of $D^{-1}C$ is less than one, and to use the subsequent convergence of the Neumann series to conclude that
\be\label{Mx-2} A = \left(I-D^{-1}C\right)^{-1}D^{-1} = \sum_{j=0}^\infty \left(D^{-1}C\right)^j D^{-1} \geq 0. \ee
Even though the vertex Gramian $G$ in (\ref{Mx-1}) is known to be nonnegative (see also Lemma~\ref{Mpro-15} below), we have not been able to prove the nonnegativity of $G$ (nor its complete positivity, see Theorem~\ref{Mth-4} below) in a manner similar as in (\ref{Mx-2}) from the properties of $C$ and the block $2\times 2$ diagonal $D$. We leave this algebraic question as challenge to the reader.  
\begin{rem}\label{Mrem-cp}{\rm If an $n$-simplex $S$ has a vertex Gramian $G$ corresponding to a vertex $v$ of $S$ all whose $(\nmo)\times(\nmo)$ principal submatrices are inverses of weakly diagonally dominant Stieltjes matrices, then each facet $F$ of $S$ with $v\in F$ is nonobtuse. This does {\em not} imply that the facet $F_v$ of $S$ opposite $v$ is nonobtuse. But it does imply that the {\em facets} of $F_v$ are nonobtuse, and thus, in the light of Theorem~\ref{Mth-2}, that $F_v$ is {\em close} to being nonobtuse itself}.
\end{rem}
An interesting consequence of the difference between odd and even $n$ is the following.
\begin{Co}\label{Mcor-7} Let $n\geq 2$ be even, and $S$ an $n$-simplex with nonobtuse facets. Then $S$ has a facet $F$ making nonobtuse angles with all other facets. Consequently, for the vertex Gramian $G$ associated with the vertex $v$ opposite $F$, the equation $Gx=e$ has a solution $x\geq 0$, or equivalently, $G$ is weakly diagonally dominant.
\end{Co}
{\bf Proof. } Each pairs of facets of $S$ can make at most one obtuse angle, and thus if the number of facets of $S$ is odd, there must be a facet making no obtuse angle with another facet.~\hfill $\Box$\\[2mm]
In the next section we will study simplices with only nonobtuse facets, that have an obtuse dihedral angle, in order to try to distinguish them from nonobtuse simplices.
\subsection{Blocking columns and nonblocking matrices}\label{Msect-3.3}
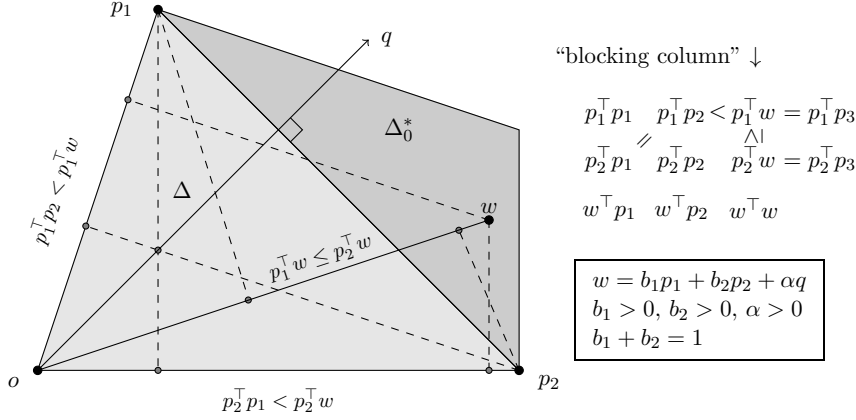
\begin{figure}[h]
\begin{center}  
\begin{tikzpicture}[scale=0.8, every node/.style={scale=0.8}]
\draw[white!80!gray, fill=white!80!gray] (0,0)--(8,0)--(2,6)--cycle;
\draw[white!60!gray, fill=white!60!gray] (8,0)--(8,4)--(2,6)--cycle;
\draw (0,0)--(8,0)--(2,6)--cycle;
\draw (8,0)--(8,4)--(2,6)--cycle;
\draw[fill=black] (0,0) circle [radius=0.07];
\draw[fill=black] (8,0) circle [radius=0.07];
\draw[fill=black] (2,6) circle [radius=0.07];
\draw[fill=black] (7.5,2.5) circle [radius=0.07];
\draw[dashed] (2,6)--(2,0);
\draw[dashed] (7.5,2.5)--(7.5,0);
\draw[dashed] (7.5,2.5)--(1.5,4.5);
\draw[dashed] (8,0)--(0.8,2.4);
\draw[dashed] (2,6)--(3.5,1.2);
\draw[dashed] (8,0)--(7,2.33);
\draw[fill=gray] (1.5,4.5) circle [radius=0.05];
\draw[fill=gray] (0.8,2.4) circle [radius=0.05];
\draw[fill=gray] (2,0) circle [radius=0.05];
\draw[fill=gray] (7.5,0) circle [radius=0.05];
\draw[fill=gray] (3.5,1.17) circle [radius=0.05];
\draw[fill=gray] (7,2.33) circle [radius=0.05];
\draw[fill=gray] (2,2) circle [radius=0.05];
\node at (8.5,-0.2) {$p_2$};
\node at (1.4,6) {$p_1$};
\node at (7.5,2.7) {$w$};
\node at (-0.4,-0.2) {$o$};
\node at (2.4,3) {$\Delta$};
\node at (6,4) {$\Delta^\ast_0$};
\node[scale=0.9] at (4,-0.5) {$p_2^\top p_1 < p_2^\top w$};
\node[rotate=18,scale=0.9] at (4.7,1.9) {$p_1^\top w \leq  p_2^\top w$};
\node[rotate=72,scale=0.9] at (0.3,3) {$p_1^\top p_2 <  p_1^\top w$};
\draw (0,0)--(7.5,2.5);
\draw[->] (0,0)--(5.5,5.5);
\draw (4.2,3.8)--(4.4,4)--(4.2,4.2);
\node at (5.8,5.5) {$q$};
\node at (11,1) {\fbox{$\begin{array}{l} w=b_1p_1+b_2p_2+\alpha q\\ b_1>0,\, b_2>0,\, \alpha>0 \\ b_1+b_2=1\end{array}$}};
\begin{scope}[shift={(9.5,2.7)}]
\node at (0,0) {$w^\top p_1$};
\node at (1.2,0) {$w^\top p_2$};
\node at (2.4,0) {$w^\top w$};
\node at (0,0.8) {$p_2^\top p_1$};
\node at (1.2,0.8) {$p_2^\top p_2$};
\node at (2.4,0.8) {$p_2^\top w$};
\node at (0,1.6) {$p_1^\top p_1$};
\node at (1.2,1.6) {$p_1^\top p_2$};
\node at (2.4,1.6) {$p_1^\top w$};
\node at (3.5,1.6) {$=p_1^\top p_3$};
\node at (3.5,0.8) {$=p_2^\top p_3$};
\node[rotate=90] at (2.4,1.2) {$\geq$};
\node at (1.8,1.6) {$<$};
\node[rotate=45] at (0.6,1.1) {$=$};
\node at (0.8,2.5) {``blocking column'' $\downarrow$};
\end{scope}
\end{tikzpicture}    
\end{center} 
\caption{\small{Illustration of the proofs of Theorem~\ref{Mth-13} and Proposition~\ref{Mprop-10} for $n=3$.}}
\label{Mfigure28}
\end{figure}

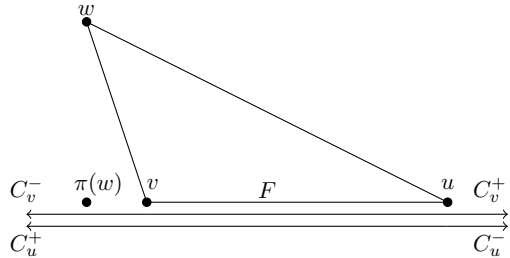
\begin{figure}[h]
\begin{center}  
\begin{tikzpicture}[scale=0.8, every node/.style={scale=0.8}]
\draw (0,0)--(5,0)--(-1,3)--cycle;
\draw[fill=black] (0,0) circle [radius=0.07];
\draw[fill=black] (5,0) circle [radius=0.07];
\draw[fill=black] (-1,3) circle [radius=0.07];
\draw[fill=black] (-1,0) circle [radius=0.07];
\node at (-1,3.2) {$w$};
\node at (0.1,0.3) {$v$};
\node at (5,0.3) {$u$};
\node at (2,0.2) {$F$};
\node at (-0.8,0.3) {$\pi(w)$};
\node at (-2,-0.7) {$C^+_u$};
\node at (5.7,0.2) {$C^+_v$};
\draw[->] (0,-0.2)--(6,-0.2);
\draw[->] (5,-0.4)--(-2,-0.4);
\draw[->] (0,-0.2)--(-2,-0.2);
\draw[->] (5,-0.4)--(6,-0.4);
\node at (-2,0.2) {$C^-_v$};
\node at (5.7,-0.7) {$C^-_u$};
\end{tikzpicture}    
\end{center} 
\caption{\small{Facet $F$ has $n$ interior cylinders, $v$ is a vertex of $\nmo$ of them, and $\pi(w)$ needs to project in $\nmo$ of them. Only if $n=2$, it may be the one of which $v$ is not a vertex.}}
\label{Mfigure29}
\end{figure}

We will now study properties of simplices with nonobtuse facets that do indeed have obtuse dihedral angles. Again, we  first consider the case $n=3$. Let $\tet$ be a tetrahedron with nonobtuse triangular facets and assume that $\tet$ is {\em not} nonobtuse. Then Theorem~\ref{Mth-2} shows that $\tet$ has a vertex that projects onto the dual hull $\Delta^\ast$ of its opposite facet $\Delta$ but not onto $\Delta$ itself. This is depicted in Figure~\ref{Mfigure28}, in which we did not draw the triangles $\Delta_1^\ast$ and $\Delta_2^\ast$ that together with $\Delta$ and $\Delta_0^\ast$ make up the dual hull of $\Delta$, and where $w$ is the projection of the fourth vertex $p_3$ of $\tet$ on the ambient plane of $\Delta$. Now, write $G_o$ for the vertex Gramian of $\tet$ corresponding to its vertex $o$. Thus, $G_o$ is the Gramian of the {\em position vectors} $p_1,p_2,p_3$  of the other vertices of $\tet$ with respect $o$, which is considered as the origin of $\RR^3$.\\[2mm] 
It is clear from Figure~\ref{Mfigure28} that if the projection $w$ of $p_3$ lies in $\Delta^\ast_0$, then both $p_1^\top p_3 = p_1^\top w > p_1^\top p_2$ and $p_2^\top p_3 = p_2^\top w > p_2^\top p_1$. This shows that the first and second entry in the third column of $G_o$ are not minimal in their respective rows. Thus, conversely, if either one of those two entries in the third column of $G_o$ is minimal in its row, then $p_3$ does not project in  $\Delta^\ast_0$. In that case, the dihedral angle between $\Delta$ and the facet of $\tet$ opposite $o$ is nonobtuse, even if $p_3$ projects elsewhere in $\Delta^\ast$. Consequently, if {\em each} column of $G_o$ has an {\em off-diagonal} entry that is minimal in its row, then the facet of $\tet$ opposite $o$ does not make an obtuse dihedral angle, and $G_o^{-1}$ is weakly diagonally dominant. Now, Theorem~\ref{Mth-1} shows the following.  
\begin{Co}\label{Mco-10} If in each vertex Gramian $G$ of a tetrahedron $\tet$ with nonobtuse facets, each column has an off-diagonal entry that is minimal in its row, then $\tet$ is nonobtuse.
\end{Co} 
The converse implication does not hold: there are points $x\in\Delta$ with $p_1^\top x > p_1^\top p_2$ and $p_2^\top x > p_2^\top p_1$. Choosing $p_3$ high enough above such a point $x$ results of course in a nonobtuse tetrahedron.
\\[2mm]
Before generalizing the above results to dimensions $n\geq 4$, we first define some terminology inspired by the case $n=3$ and make some observations.
\begin{Def}[Blocking column, nonblocking matrix]{\rm Let $A\in\Rnspd$, $A=(a_{ij})\geq 0$. A column $a_j=Ae_j$ of $A$ with the property that no entry $a_{ij}, i\not=j$, of $a_j$ is minimal in its row in $A$, is called a {\em blocking column} \index{blocking column}. If $A$ has no blocking columns, which in terms of its entries means that 
\be \forall j\in I_1^n:\, \exists i \in I_1^n\setminus\{j\}:\, \forall k\in I_1^n: \hdrie a_{ij} \leq a_{kj}, \ee
then $A$ will be called a {\em nonblocking matrix}}. \index{nonblocking matrix}
\end{Def}
The matrix property of having a blocking column is trivially invariant under $\overset{\pi}{\sim}$. Invariance under $\sim$ is a much more subtle issue. In Section~\ref{Msect-5} we will prove that for $n=3$ it is indeed invariant under $\sim$. However, for $n \geq 4$ it is not. See for instance the two vertex Gramians $G_0,G_2$ corresponding to the nonobtuse $4$-simplex whose vertices $v_0,\dots,v_4$ are the origin and the four columns of $P$,
\be\label{MG0} \small{P=\left[\begin{array}{rrrr} 0 & 4 & 4 & 3 \\ 4 & 0 & 4 & 2\\ 4 & 4 & 0 & 2 \\ 0 & 0 & 0 & 4\end{array}\right], \und G_0=\left[\begin{array}{rrrr} 32  &  16  &  16  &  16\\
    16  &  32  &  16 &   20\\ 
    16  &  16  &  32 &   20\\ 
    16  &  20  &  20 &   33 \end{array}\right] \hdrie \sim \hdrie \ \left[\begin{array}{rrrr}   32  &  16  &  16  &  12\\
    16  &  32  &  16  &  12\\
    16  &  16  &  32  &  16\\
    12  &  12  &  16  &  25\end{array}\right] = G_2.}\ee
We see that $G_0$ is nonblocking, but that $G_2$ has a blocking (third) column. Since the $3\times 3$ bottom right {\em submatrix} of $G_0$ also has a blocking column, we conclude that the property of being nonblocking is generally not inherited by principal submatrices. Conversely, $G_2$ has nonblocking principal submatrices of size $3\times 3$, but also one having a blocking column.\\[2mm]
The latter observation holds more generally, in the following sense.
\begin{Pro}\label{Mprop-10} Let $A\in\Rnspd$ with $n\geq 3$. If $A$ has a blocking column, then $A$ has a $3\times 3$ principal submatrix that possesses a blocking column.
\end{Pro}
{\bf Proof. } Without loss of generality, we may assume that the last column of $A$ is blocking and that its entries are non-decreasing from top to bottom: $a_{1n} \leq\cdots\leq a_{nn}$. Since $A$ is nonblocking by assumption, there exists a $j\in I_2^{n-1}$ such that $a_{1j}<a_{1n}$. But then due to $a_{1n} \leq a_{jn}$ and $a_{j1}=a_{1j}$, the $3\times 3$ principal submatrix of $A$ consisting of the entries $a_{\ell k}$ with $\ell,k\in\{1,j,n\}$ has a blocking third column.~\hfill $\Box$\\[2mm]
The reader may be tempted to believe that if $A$ has a blocking column, then for {\em any} $k<n$ it has a principal submatrix of size $k\times k$ with a blocking column. This is however not true, and here we present a counter example for $n=5$ and $k=4$,
\be\label{Mexam-1}  \small A = \left[\begin{array}{cccc|c} 2.00  &  1.20 &   \fbox{1.00}  &  1.20 &   1.10\\
    1.20  &  2.00  &  1.20  &  \fbox{1.05}  &  1.10\\
   \fbox{1.00}  &  1.20 &  2.00  &  1.20  &  1.10\\
    1.20  &  \fbox{1.05}  &  1.20  &  2.00  &  1.10\\ \hline
    1.10  &  1.10  &  1.10  &  1.10  &  2.00\end{array}\right], \hdrie A_{\{1,3,5\}} = \left[\begin{array}{lll} 2.00 & 1.00 & 1.10 \\  1.00 & 2.00 & 1.10 \\ 1.10 & 1.10 & 2.00\end{array}\right] \ee
The fifth column of $A$ is blocking, but each of its five $4\times 4$ principal submatrices is nonblocking. It can also be verified that $A$ is positive definite and, in fact, an element of $\MM_\dd$. However, following the proof of Proposition~\ref{Mprop-10} above, the $3\times 3$ principal submatrix $A_{\{1,3,5\}} = (a_{ij})$ with $i,j\in\{1,3,5\}$ of $A$ has indeed a blocking column. This subtlety is expressed as follows.
\begin{Co}\label{Mcor-x} If all $3\times 3$ principal submatrices of $A\in\Rnspd$ with $A\geq 0$ are nonblocking, then so are all $k\times k$ principal submatrices with $4\leq k \leq n$, including $A$. The same conclusion does not hold if merely all $4\times 4$ principal submatrices are nonblocking.
\end{Co}
{\bf Proof. } This follows from Proposition~\ref{Mprop-10} and the example in (\ref{Mexam-1}).~\hfill $\Box$\\[2mm]
Postponing the consequences and using the new terminology, we first return to the generalization to arbitrary dimension of the earlier observations in this section for the case $n=3$. To better understand the proof of the next theorem, see also Figure~\ref{Mfigure28}. 
\begin{Th}\label{Mth-13} Let $n\geq 3$. Let $S$ be an $n$-simplex with only nonobtuse facets. Assume that $S$ has an obtuse dihedral angle. Then $S$ has a vertex Gramian $G$ with a blocking column.
\end{Th}   
{\bf Proof. } Denote the vertices of $S$ by $v_0,\dots,v_n$ and their opposite facets by $F_0,\dots,F_n$, respectively. Because $S$ has an obtuse dihedral angle, it has a vertex that projects into the dual hull of its opposite facet but not in that facet itself. Using Theorem~\ref{Mth-2} and Proposition~\ref{Mpro-6} and without loss of generality, assume that $v_n$ projects in the interior of the part $(F_n)_0^\ast$ of $F_n^\ast$. Denote this projection by $w$. Consider the vertex $v_0$ of $F_n$ opposite its facet $H=F_n\cap (F_n)_0^\ast$ as the origin of $\RR^n$. Let $P_0$ be the matrix with the position vectors $v_1,\dots,v_{n-1}$ as columns. Then $G_0=P_0^\top P_0$ is the vertex Gramian of $F_n$ associated with $v_0$, and 
\be\label{Meq-101} G = \left[ \begin{array}{c|c} G_0 & g \\ \hline g^\top & \gamma \end{array}\right] \hdrie\mbox{\rm with}\hdrie g=P_0^\top w=P_0^\top v_n \und \gamma=v_n^\top v_n\ee
a vertex Gramian of $S$ itself. Now, any point in $H$ is a convex combination of the columns of $P_0$. Because $w$ lies further away from $v_0$ than $H$, we can write   
\be w=P_0b+\alpha q, \hdrie \mbox{\rm where $b>0,\,\, e^\top b=1$, and $\alpha>0$}, \ee 
where $q$ is the {\em outward} pointing normal to $H$ defined by $P_0^\top q = e$. See (\ref{Mq0}) for the {\em inward} normal. Finally, let $\ell\in I_1^{\nmo}$. Then, as $v_\ell$ is a column of $P_0$, we have $v_\ell^\top q = 1$ and
\be v_\ell^\top w = v_\ell^\top(P_0b+\alpha q) = e_\ell^\top G_0b+\alpha. \ee
Due to $b>0, e^\top b=1$, the expression $(e_\ell^\top G_0)b$ is a convex combination of the entries in the $\ell$-th row of $G$. Thus, as $\alpha>0$, $v_\ell^\top w$ is strictly larger than at least one entry in that row, proving that the rightmost column of $G$ in (\ref{Meq-101}) is a blocking column.~\hfill $\Box$\\[2mm]
Note that the above theorem shows that if the dihedral angle between $F_i$ and $F_j$ is obtuse, the vertex Gramian corresponding to not only $v_i$ but also to $v_j$ contains a blocking column. Thus, if $n$ of the $n+1$ vertex Gramians of $S$ are nonblocking, also its remaining vertex Gramian must be nonblocking, as an obtuse angle is shared between two facets.\\[2mm]
We can now formulate two theorems and a corollary, which are the culmination of the new concepts introduced. In Sections~\ref{Msect-4} to~\ref{Msect-6} they will be of central importance.
\begin{Th}\label{Mth-17} Let $S$ be a simplex with nonobtuse facets, and $G$ a vertex Gramian of $S$.\\[2mm]
$\bullet$ If $G$ is a nonblocking vertex Gramian of $S$, then $Gx=e$ has a solution $x\geq 0$.\\[2mm]
$\bullet$ If all vertex Gramians of $S$ are nonblocking, then $S$ is nonobtuse.
\end{Th}  
{\bf Proof. } The proof of Theorem~\ref{Mth-13} shows that if the vertex Gramian $G$ corresponding to a vertex $v$ of $S$  has no blocking column, the facet $F$ opposite $v$ makes no obtuse dihedral angle. This proves the first bullet. The second bullet follows from Theorem~\ref{Mth-1}.~\hfill $\Box$ 
\begin{Th}\label{Mth-8} Let $k\geq 2$, and $S$ an $n$-simplex with nonobtuse $k$-facets. If all $m\times m$ principal submatrices of all vertex Gramians of $S$ with $m>k$ are nonblocking, then $S$ is nonobtuse. 
\end{Th}  
{\bf Proof. } Let $F$ be  $(k+1)$-facet of $S$. By assumtion, it has nonobtuse $k$-facets only, and all its vertex Gramians are nonblocking. Thus, $F$ is nonobtuse by Theorem~\ref{Mth-17}. The statement now follows by induction with respect to $k$.~\hfill $\Box$
\begin{Co}\label{Mco-9} Let $S$ be an $n$-simplex with nonobtuse triangular facets. If each $3\times 3$ principal submatrix of each vertex Gramian of $S$ is nonblocking, then $S$ is nonobtuse. 
\end{Co}  
{\bf Proof. } Combine Theorem~\ref{Mth-8} with Corollary~\ref{Mcor-x}.~\hfill $\Box$\\[2mm]
The final theorem in this section may come as a surprise, as it obviously does not hold for dimension $n=2$. However, in spite of possible obtuse dihedral angles, simplices with nonobtuse facets can be isometrically embedded in the nonnegative orthant of $\RR^n$ with any of its vertices located at the origin. The theorem deserves its place for the following reason. In the finite element method~\cite{Bra} to approximate solutions of elliptic PDEs, the approximations satisfy {\em discete maximum principles} under the {\em sufficient} condition that the triangulation of the domain into tetrahedra has no obtuse angles \cite{BrKoKr}. This condition leads to M-matrices, whose inverses are completely positive due to Corollary~\ref{Mcor-2}. It may be possible to relax this condition and aim for matrices as in (\ref{Mx-1}), whose inverses are still completely positive due to Theorem~\ref{Mth-4} below. This will be a topic for another paper. See also Remark~\ref{Mrem-x}. 
\begin{Th}\label{Mth-4} Let $n\geq 3$. Then each vertex Gramian $G$ of an $n$-simplex $S$ with nonobtuse facets is completely positive with cp-rank $n$.
\end{Th} 
{\bf Proof. } Choose the origin of $\RR^n$ in a vertex $v$ of $S$. Let $F$ be a facet of $S$ with $v\in F$, and let $w$ be the vertex of $S$ opposite $F$. Theorem~\ref{Mth-2} shows that $w$ projects onto the dual hull $F^\ast$ of $F$. Denote this projection by $\pi(w)$. Then according to (\ref{MDandE}), $\pi(w)$ lies in at least $\nmo$ of the $n$ interior cyclinders associated with the $n$ facets of $F$. Also, $v$ is a vertex of exactly $\nmo$ of these cylinders. Thus, by the pigeon hole principle there is one cylinder $C^+$ associated with a facet $G$ of $F$, such that $v$ is vertex of $C^+$ and $\pi(w)\in C^+$ if and only if $n\geq 3$. Now, since $G$ is nonobtuse, Corollary~\ref{Mcor-2} yields that it can be isometrically embedded into $\RR^{n-2}_{\geq 0}$. Its interior cylinder $C^+$, which contains $F$ and $\pi(w)$, can then be isometrically embedded into $\RR^{\nmo}_{\geq 0}$. Since $w$ projects into it, a possible change of sign in the last coordinate establishes the isometrical embedding of $S$ into $\RR^n_{\geq 0}$. Thus, the vertex Gramian associated with $v$ is completely positive with cp-rank $n$.~\hfill$\Box$\\[2mm]
In terms of the chain (\ref{Mchain}) of simplicial matric classes, Theorem~\ref{Mth-4} shows that $\VV_n^{n-1}\subset \PP_\dd$.    
\begin{rem}\label{Mrem-5}{\rm Any $1$-simplex $S$ is formally nonobtuse, as the normals to its two facets make an angle $\pi$ and thus the dihedral angle between them equals zero. Thus, any triangle $\Delta$ has only nonobtuse facets. Theorem~\ref{Mth-2} confirms that $\Delta$ has at most one obtuse angle at each facet, and Corollary~\ref{Mcor-7} that $\Delta$ has a facet making nonobtuse angles with both other facets. Theorem~\ref{Mth-4} obviously does not hold true for a nonobtuse triangle $\Delta$. To see where the {\em proof} of Theorem~\ref{Mth-4} fails, see Figure~\ref{Mfigure29}}.
\end{rem}  
Note that due to Remark~\ref{Mrem-cp} we have that for $G$ to be an element of $\VV_n^{n-1}$ it is necessary but {\em not sufficient} that all $(n-1)\times (n-1)$ principal submatrices of $G$ are in $\VV_{n-1}^{n-1}$, as this does {\em not} force the facet of $\SIS(G)$ opposite $v$ to be nonobtuse. Nonetheless, for a given matrix $A\in\VV_n^2$ it is a simple, finite computation to verify if all $(n-1)$-facets of its underlying simplex $\SIS(A)$ are nonobtuse. As a first question for further research on the position of $\PP_\dd$ in the chain (\ref{Mchain}), it would be of interest to verify if the factual $\VV_5^4\subset\PP_\dd\subset \VV_5^2$ could be improved to $\VV_5^3\subset\PP_\dd\subset \VV_5^2$, or perhaps even the highly speculative $\VV_5^3=\PP_\dd$, but it will be not considered in this chapter. 
\subsection{$n$-Simplices with nonobtuse triangular facets}\label{Msect-3.4}
The previous sections were on $n$-simplices with nonobtuse $(\nmo)$-facets. Here we prove some dual results about $n$-simplices with nonobtuse {\em triangular} facets. For tetrahedra, these results coincide, as the $(\nmo)$-facets of a tetrahedron are its triangular facets.\\[2mm]
First recall Proposition~\ref{Mpro-15} and Theorem~\ref{Mth-1b}. There we proved that the set of vertex Gramians of simplices whose triangular facets are all nonobtuse, equals the simplicial matrix class $\DD_\dd$ of nonsingular, doubly nonnegative, pointwise weakly diagonally dominant matrices. \\[2mm]
The following theorem is the counterpart of Theorem~\ref{Mth-2}.  
\begin{Th}\label{Mth-3} Let $S$ be an $n$-simplex with nonobtuse triangular facets. Then:\\[2mm]
$\bullet$ each facet of $S$ makes an obtuse angle with at most $n-2$ of the remaining $n$ facets;\\[2mm]
$\bullet$ $S$ has at most $\lfloor \half(\npo)(n-2)\rfloor$ obtuse dihedral angles.
\end{Th}  
{\bf Proof}. Write $v_0,\dots,v_n$ for the vertices of $S$, and $F_0,\dots,F_n$ for the respective opposite facets with normals $q_0,\dots,q_n$. It suffices to prove the statement for one facet, say $F_1$. Consider the vertex Gramian $G_0$ of $S$ corresponding to $v_0$. The solution $x$ of $G_0x=e_1$ is the first column of $G_0^{-1}$, which equals the $n$-vector of inner products $q_1^\top q_1,q_2^\top q_1,\dots,q_n^\top q_1$. Now, $G_0\geq 0$ due to Proposition~\ref{Mpro-15}, and the diagonal entries of $G_0$ are positive, because $G_0$ is positive definite. This implies that $x>0$ contradicts $Gx=e_1$. Thus, $x$ has at least one nonpositive entry. This proves that $F_1$ makes a nonobtuse angle with at least one of $F_2,\dots,F_n$, say $F_i$. Consider next the vertex Gramian $G_i$ corresponding to the vertex $v_i$. One of the columns of $G_i^{-1}$ contains the inner products between $q_1$ and all normals, except $q_i$, because $q_i$ is the normal to the facet opposite $v_i$ Thus, for the same reason as above, $F_1$ also makes a nonobtuse angle with a facet $F_\ell$ with $1\not=\ell\not=i$. This proves the statement.~\hfill $\Box$\\[3mm]
Theorem~\ref{Mth-3} may seem of marginal value, especially when compared to Theorem~\ref{Mth-2}. Therefore, we will elaborate on its interest. First of all, Corollary~\ref{Mco-9} is about such simplices, and in Section~\ref{Msect-5} this will result in theorems about ultrametric matrices. Secondly, the result fits into Conjecture~\ref{Mcon-1}. For $k=1$, this conjecture encompasses the result by Fiedler \cite{Fie} that each facet of an arbitrary $n$-simplex $S$, all whose $1$-facets are nonobtuse due to Remark~\ref{Mrem-5}, makes at most $\nmo$ obtuse dihedral angles with the remaining facets. If $k=n$, the simplex $S$ is in fact assumed to be nonobtuse, hence the conjecture holds trivially. Theorem~\ref{Mth-2} and Theorem~\ref{Mth-3} prove the conjecture for $k=\nmo$ and $k=2$, respectively.   
\begin{Con}\label{Mcon-1} Let $S$ be an $n$-simplex whose $k$-facets are all nonobtuse. Then:\\[2mm]
$\bullet$ each normal of $S$ makes an obtuse angle with at most $n-k$ of the other $n$ normals;\\[2mm]
$\bullet$ $S$ has at most $\lfloor \half(\npo)(n-k)\rfloor$ obtuse dihedral angles.\\[2mm]
If $G$ is a vertex Gramian of $S$, consider the solutions $x$ and $y_j$ of
\be Gx = e \und Gy_j = e_j. \ee
Then $x$ has at most $n-k$ negative entries and each $y_j$ has at most $k+1$ positive entries. 
\end{Con}   
A third reason for the interest in Theorem~\ref{Mth-3} but also in the above conjecture is, that in applications, it is not always relevant to know whether a given symmetric positive definite nonnegative matrix is an inverse M-matrix (which corresponds to the symmetric inverse M-matrix problem). More interesting is to know if a given matrix $A$ of a linear system $Ax=b$ has a nonnegative inverse, such that if $b\geq 0$, then $x=A^{-1}b \geq 0$. As noted already before, this is for instance relevant in the context of discrete maximum and comparison principles in finite element methods \cite{Bra,BrKoKr}. If the simplex $S=\SIS(A)$ underlying $A$ is nonobtuse, then $A$ is obviously nonnegative and in fact, even completely positive (see Corollary~\ref{Mcor-2}). However, the property $A\geq 0$ already holds if $S$ has merely nonobtuse {\em triangular} facets, which is a much weaker condition, and which allows $A^{-1}$ to have many positive off-diagonal entries.\\[2mm] 
A final reason for the interest in simplices with nonobtuse triangular facets comes from the context of $0/1$-matrices. Let $P$ be an $n\times n$ non-singular $0/1$-matrix. The origin together with the columns of $P$ constitute an $n$-simplex. Now, observe that any triple of vertices of the unit $n$-cube is a nonobtuse triangle. As such, $0/1$-simplices and more general $0/1$-polytopes \cite{KaZi} have restricted dihedral angle properties, in comparison with general simplices and polytopes.
\begin{rem}{\rm Note that $0/1$-matrices also play an important role in the context of ultrametric matrices. In fact, any symmetric ultrametric matrix with rational entries is, apart from a multiplicative factor, the Gramian of a $0/1$-matrix. See \cite{NabVar}}.
\end{rem} 
Having studied simplices with nonobtuse triangular facets in general, in Section~\ref{Msect-4} we continue with studying simplices with {\em special} nonobtuse triangular facets.
\section{Simplices with special triangular facets}\label{Msect-4}
There are three well-known types of triangles such that each $n$-simplex $S$ having only such triangular facets is nonobtuse. They are the {\em equilateral} triangle, the {\em right} triangle, and the {\em strongly isosceles} triangle, which is a triangle having no unique longest edge. 
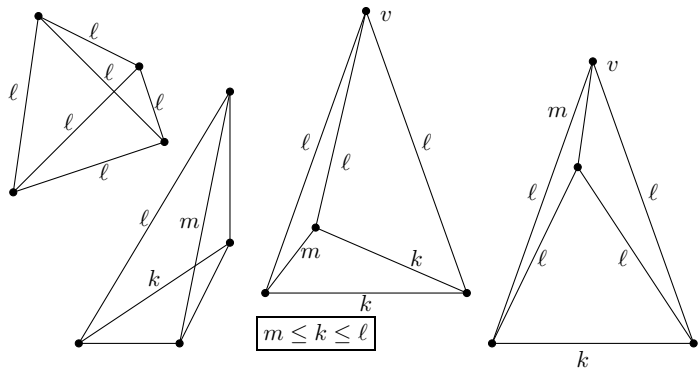
\begin{figure}[h]
\begin{center}  
\begin{tikzpicture}[scale=0.67, every node/.style={scale=0.8}]
\begin{scope}[shift={(0,2)}]
\draw (0,0)--(3,1)--(0.5,3.5)--cycle;
\draw (0.5,3.5)--(2.5,2.5)--(3,1);
\draw (0,0)--(2.5,2.5);
\node at (1.8,0.4) {$\ell$};
\node at (1.1,1.4) {$\ell$};
\node at (0,2) {$\ell$};
\node at (1.6,3.2) {$\ell$};
\node at (1.9,2.4) {$\ell$};
\node at (2.9,1.8) {$\ell$};
\draw[fill=black] (3,1) circle [radius=0.07];
\draw[fill=black] (0,0) circle [radius=0.07];
\draw[fill=black] (0.5,3.5) circle [radius=0.07];
\draw[fill=black] (2.5,2.5) circle [radius=0.07];
\end{scope}
\begin{scope}[shift={(1.3,-1)}]
\draw (0,0)--(2,0)--(3,2)--(3,5)--cycle;
\draw (0,0)--(3,2);
\draw (2,0)--(3,5);
\draw[fill=black] (0,0) circle [radius=0.07];
\draw[fill=black] (2,0) circle [radius=0.07];
\draw[fill=black] (3,2) circle [radius=0.07];
\draw[fill=black] (3,5) circle [radius=0.07];
\node at (1.3,2.5) {$\ell$};
\node at (1.5,1.3) {$k$};
\node at (2.2,2.4) {$m$}; 
\end{scope}
\begin{scope}[shift={(5,0)}]
\draw[fill=black] (0,0) circle [radius=0.07];
\draw[fill=black] (4,0) circle [radius=0.07];
\draw[fill=black] (2,5.6) circle [radius=0.07];
\draw[fill=black] (1,1.3) circle [radius=0.07];
\draw (0,0)--(2,5.6)--(4,0)--cycle;
\draw (0,0)--(1,1.3)--(4,0);
\draw (1,1.3)--(2,5.6);
\node at (0.8,3) {$\ell$};
\node at (3.2,3) {$\ell$};
\node at (2,-0.2) {$k$};
\node at (3,0.7) {$k$};
\node at (1.6,2.6) {$\ell$};
\node at (0.9,0.8) {$m$};
\node at (2.4,5.5) {$v$};
\node at (1,-0.8) {\fbox{$m\leq k \leq \ell$}};
\end{scope}
\begin{scope}[shift={(9.5,-1)}]
\node at (2.4,5.5) {$v$};
\draw[fill=black] (0,0) circle [radius=0.07];
\draw[fill=black] (4,0) circle [radius=0.07];
\draw[fill=black] (2,5.6) circle [radius=0.07];
\draw[fill=black] (1.7,3.5) circle [radius=0.07];
\draw (0,0)--(2,5.6)--(4,0)--cycle;
\draw (0,0)--(1.7,3.5)--(4,0);
\draw (1.7,3.5)--(2,5.6);
\node at (1.8,-0.3) {$k$};
\node at (2.6,1.7) {$\ell$};
\node at (0.8,3) {$\ell$};
\node at (3.2,3) {$\ell$};
\node at (1.3,4.6) {$m$};
\node at (1,1.7) {$\ell$};
\end{scope}
\end{tikzpicture}    
\end{center} 
\caption{\small{Tetrahedra with only equilateral, only right, and only strongly isosceles facets.}}
\label{Mfigure210}
\end{figure}\\[2mm]
Simplices with equilateral triangular facets are the {\em regular simplices}; the ones with only right triangular facets are the {\em Schl\"afli orthoschemes} \cite{Cox,Cox2,Fie,Schl} or {\em path-simplices} \cite{BrKoKr,BrKoKrSo}, because they have a path of mutually orthogonal edges. Simplices with strongly isosceles triangular facets are the least known of the three: they are precisely the simplices whose vertex sets are {\em ultrametric spaces} with respect to their Euclidean distances. These simplices were proved to be nonobtuse (in fact, even acute) by Fiedler in \cite{Fie2}. For completeness, in this section we study the simplicial matrix classes that each of the three types of simplices induce.\\[2mm]
In Figure~\ref{Mfigure210} we display the only type of tetrahedron having only equilateral facets, which is the regular tetrahedron, and the only type of tetrahedron having only right triangular facets, which is the path tetrahedron. Also displayed are the only two distinct types of tetrahedra having only strongly isosceles triangular facets. To see that these are the only types possible, let $\ell$ be the length of the longest edge of a tetrahedron $\tet$ with only strongly isosceles triangular faces. Then $\tet$ has a triangular face $\Delta$ with two edges of lengths $\ell$ meeting at a vertex $v$ of $\Delta$. If the third edge of $\tet$ that meets at $v$ also has length $\ell$, the facet opposite $v$ is strongly isosceles with edge lengths $m$ and $k$ with $m\leq k\leq\ell$. We call this a Type I tetrahedron. If the third edge that meets at $v$ is shorter than $\ell$, then both the two triangular faces of $\tet$ that share this edge have two edges of length $\ell$. This results in a tetrahedron of Type II, in which without loss of generality, $m\leq k$. The intersection of both types consists of tetrahedra with at least five longest edges of length $\ell$, and the regular tetrahedron in particular.\\[2mm]
The following observations are simple but instructive. Let $\tet$ be a tetrahedron and $v$ a vertex of $\tet$. If the three triangular faces $\Delta$ of $\tet$ with $v\in \Delta$ are equilateral, then the face of $\tet$ opposite $v$ is also equilateral. A  similar conclusion does not hold for tetrahedra with three right or three strongly isosceles facets. To control the shape of their triangular facets opposite a chosen vertex $v$ of an $n$-simplex $S$, it is necessary to impose conditions on all three-dimensional facets $\tet$ of $S$ containing $v$. This is not a deep observation, but we present the full details, because in Theorem~\ref{Mth-6} we use the same idea in a rather more complicated setting.
\begin{Pro}\label{Mpro-8} Let $v$ be a vertex of an $n$-simplex $S$. Then:\\[2mm]
$\bullet$ If all tetrahedral faces $\tet$ of $S$ with $v\in \tet$ are of Type I or II, then all tetrahedral faces are.\\[2mm]
$\bullet$ If all tetrahedral faces $\tet$ of $S$ with $v\in \tet$ are path tetrahedra, then all tetrahedral faces are.
\end{Pro} 
{\bf Proof. } If all tetrahedral faces $\tet$ of $S$ with $v\in \tet$ are of Type I or II, each triangular face $\Delta$ of $S$ with $v\in \Delta$ is strongly isosceles. If $\Delta$ is a triangular face with $v\not\in\Delta$, the convex hull of $v$ and $\Delta$ is a tetrahedron $\tet$ of Type I or II with $v\in \tet$. Since $\tet$ has strongly isosceles faces, $\Delta$ is strongly isosceles. Hence, all triangular faces of $S$ are strongly isosceles, and thus all tetrahedral faces are of Type I or II. The second statement can be proved similarly.~\hfill $\Box$\\[2mm]
This proposition shows that by merely inspecting all $3\times 3$ principal submatrices of one arbitrarily given vertex Gramian $G_v$ of a simplex $S$, it is possible to decide whether a simplex $S$ has only strongly isosceles triangular facets or only right triangular facets. For isosceles triangular facets this can be done by inspecting all $2\times 2$ submatrices of one vertex Gramian.
\begin{Pro}\label{Mpro-4} Let $S$ be an $n$-simplex. Then the following are equivalent:\\[2mm]
$\bullet$ $S$ is a regular simplex;\\[2mm]
$\bullet$ $S$ has a vertex $v$ such that all vertex Gramians associated with $v$ of each triangular face $\Delta$ of $S$ with $v\in\Delta$ are equal, and equal to a positive multiple of the $2\times 2$ matrix $I+ee^\top$;\\[2mm]
$\bullet$ $S$ has one vertex Gramian $G$ that is a positive multiple of $I+ee^\top$;\\[2mm]
$\bullet$ all vertex Gramians of $S$ are equal, and equal to a positive multiple of $I+ee^\top$.
\end{Pro} 
{\bf Proof. } The lengths of all edges of a regular simplex $S$ are equal, and each vertex of $S$ projects on the midpoint of any opposite edge. This proves that the off-diagonal entries of a vertex Gramian are half the diagonal values, which implies all the claims.~\hfill $\Box$\\[2mm]
Thus, the vertex Gramians of regular simplices form a small simplicial matrix class $\CC_r$, in the sense of Definition~\ref{Msmc}. It is a subclass of $\MM_\dd$ because of Corollary~\ref{Mcor-5} and the fact that $(I+ee^\top)e=(\npo)e$, which shows that $(I+ee^\top)y=e$ has a solution $y\geq 0$. Alternatively, one may say that each $3\times 3$ principal submatrix of each vertex Gramian is trivially nonblocking, and invoke Corollary~\ref{Mco-9}. \\[2mm]
Next, consider the path simplex $S$. Let $d_1,\dots,d_n$ be the lengths of the edges in its path of orthogonal edges, in consecutive order. Then
\be\label{Mdef-T} G_v  \overset{\pi}{\sim} (DT)^\top DT, \hdrie \mbox{where}\hdrie D=\left[\begin{array}{ccc} d_1 & & \\ & \ddots & \\ & & d_{n}\end{array}\right] \und T=\left[\begin{array}{ccc} 1 & \dots & 1\\ & \ddots  & \vdots \\  & & 1\end{array}\right] \ee
is the vertex Gramian of $S$ associated with a vertex at the longest edge of $S$. Indeed, the edges of respective lengths $d_1,d_2,\dots,d_n$ can be put consecutively in the coordinate axes directions, such that its vertices are the origin and the columns of $P=DT$. Now, observe that
\be\label{Meq-16} G_v \overset{\pi}{\sim} \left[\begin{array}{cccc} \alpha_1 & \dots & \dots & \alpha_1 \\ \vdots & \alpha_2 & \dots & \alpha_2\\ \vdots & \vdots & \ddots & \\ \alpha_1 & \alpha_2 & & \alpha_n\end{array}\right],\hdrie\mbox{\rm where }\hdrie \alpha_j = \sum_{i=1}^j d_i^2, \ee
which is a {\em type D-matrix} introduced by Markham \cite{Mar}.
\begin{Def}{\rm An $n\times n$ matrix $A=(a_{ij})$ is a {\em type-D matrix}~\cite{Mar} \index{type-D} if there exist numbers $\alpha_n>\alpha_{\nmo}>\dots>\alpha_1$ such that}
\be\label{Mtype-D} a_{ij} = \alpha_i \hdrie\mbox{\rm if $i\leq j$} \und a_{ij} = \alpha_j \hdrie\mbox{\rm if $i>j$}. \ee
\end{Def}  
Thus, the simplex underlying a symmetric positive definite type-D matrix is a path simplex.\\[2mm]
The vertex Gramians of all path simplices form a simplicial matrix class $\CC_p$ in the sense of Definition~\ref{Msmc}. This class contains the type-D matrices and their simultaneous row-column permutations, but not exclusively so: choosing the origin at an arbitrary vertex $v$ of $S$, we can follow the path of orthogonal edges from $v$ to one of the end points and choose the first, say $\ell$ coordinate directions along those edges. If the remaining coordinate directions are laid along the remaining $n-\ell$ orthogonal edges, starting at $v$ and following the path, we find that
\be\label{Meq-17} G_v  \overset{\pi}{\sim} (D\hat{T})^\top D \hat{T} \hdrie\mbox{\rm with }\hdrie \hat{T} =\left[\begin{array}{rr} T_\ell\\ & T_{n-\ell} \end{array}\right], \ee
where $T_\ell$ is the $\ell\times \ell$ upper triangular matrix as in (\ref{Mdef-T}) and hence $G_v$ is a block diagonal matrix with one or two diagonal blocks, both of type-D. Formally $G_v$ in (\ref{Meq-17}) is not of type-D, but symmetric ultrametric. See also the Introduction for an explicit example for $n=4$.
\begin{Pro}\label{Mpro-5} Let $S$ be an $n$-simplex. Then the following are equivalent:\\[2mm]
$\bullet$ $S$ is a path simplex;\\[2mm]
$\bullet$ each vertex Gramian $G$ of $S$ is $\overset{\pi}{\sim}$ to block diagonal with one or two diagonal type-D blocks;\\[2mm]
$\bullet$ $S$ has a vertex Gramian that is block diagonal with one or two type-D blocks;\\[2mm]
$\bullet$ $S$ has a vertex $v$ such that each vertex Gramian associated with $v$ of each tetrahedral face $\tet$ of $S$ with $v\in\tet$ is $\overset{\pi}{\sim}$ to $3\times 3$ block-diagonal type-D matrices with one or two diagonal blocks.
\end{Pro}  
The class $\CC_p$ is a subclass of $\MM_\dd$. To prove this, by Corollary~\ref{Mcor-5} it suffices to show that there exists $y\geq 0$ satisfying $Gy=e$ for each $G\in\CC_p$. Again, this is a triviality. Indeed, if $G$ is type-D as in (\ref{Meq-16}), then $Ge_1=\alpha_1 e$. In the case $G$ has two type-D blocks as in (\ref{Meq-17}), then $G(\alpha_{\ell+1} e_1+\alpha_1 e_{\ell+1})=\alpha_1\alpha_{\ell+1} e$. Alternatively, all $3\times 3$ submatrices of all the vertex Gramians of $S$ are nonblocking, hence also Corollary~\ref{Mco-9} proves the statement.
\begin{rem}{\rm The reader may observe that the set of vertex Gramians of all simplices underlying matrices with {\em any} number of type-D diagonal blocks might form a simplicial matrix class $\CC_o$ with $C_p\subset C_o$. This is true, but the underlying simplices can no longer be characterized in terms of their {\em triangular} facets only}.
\end{rem}   
To conclude this section on simplices that are nonobtuse due to their special triangular facets, we consider simplices with strongly isosceles facets only. As demonstrated by Fiedler in \cite{Fie2}, the vertex sets $X$ of such $n$-simplices are the only subsets of $\RR^n$ that, endowed with the Euclidean metric $d$, form an {\em ultrametric space} \index{ultrametric space} $(X,d)$ in the sense that 
\be d(x,y) \leq \max\{ \,d(x,z), d(z,y)\,\} \hdrie\mbox{\rm for all}\hdrie x,y,z\in X.\ee
Fiedler also proved that such $n$-simplices have only acute dihedral angles. Here, we will merely show their nonobtusity, using the concepts and insights developed in Section~\ref{Msect-3}.\\[2mm] 
Consider the tetrahedra of Type I and Type II, in the right in Figure~\ref{Mfigure210}. Of the four vertex Gramians of those of Type I, two are the same due to symmetry. The remaining three are, with $2\ell>\beta\geq\alpha\geq \ell,\ 2k> \gamma\geq k$
\be\label{Meq-28} \left[\begin{array}{rrr} \ell & \beta & \alpha \\ \beta & \ell & \alpha \\ \alpha & \alpha & \ell \end{array}\right], \hdrie  \half\left[\begin{array}{rrr} 2\ell & k & m \\ k & 2k &  m \\ m &  m & 2m \end{array}\right], \hdrie \half\left[\begin{array}{rrr} 2\ell & k &  k \\ k & 2k & \gamma \\ k & \gamma & 2k \end{array}\right], \ee
and with $\ell\geq k\geq m$ just as we assumed in Figure~\ref{Mfigure210}, whereas of the four vertex Gramians of a tetrahedron of Type II, only two are distinct, being
\be \half\left[\begin{array}{rrr} 2\ell & \alpha & k \\ \alpha & 2\ell &  k \\ k &  k & 2k \end{array}\right] \und \half\left[\begin{array}{rrr} 2\ell & \beta & m \\ \beta & 2\ell & m\\ m & m & 2m \end{array}\right], \ 2\ell>2\alpha\geq \ell, \ 2\ell>2\beta\geq \ell. \ee
As already mentioned in Proposition~\ref{Mpro-8}, these vertex Gramians cannot be completely characterized by their $2\times 2$ submatrices, although four of the five can: in case two diagonal entries are not equal, the corresponding off-diagonal entry is half times the smallest of the two, and in case two diagonal entries are equal, the corresponding off-diagonal entry is at least a half times this value. However, in the left matrix in (\ref{Meq-28}), this leaves open the option that $\beta<\alpha$. But then it is precisely the vertex Gramian of a tetrahedron $\tet$ with a triangular facet opposite the origin that, though isosceles, is not strongly isosceles. See Proposition~\ref{Mpro-8}.\\[2mm]
The above shows that any vertex Gramian $G$ of an $n$-simplex $S$ having only strongly isosceles triangular facets is permutation equivalent to a matrix
\be G \overset{\pi}{\sim} \left[\begin{array}{c|c|c} A_{11} & \dots & A_{1k} \\ \hline \vdots & & \vdots\\ \hline  A_{k1} & \dots & A_{kk}\end{array}\right], \ee 
where each of the the diagonal blocks $A_{jj}$ has constant diagonal with value $\lambda_j$, and $\lambda_1\geq\cdots\geq \lambda_k$ and where   
\be A_{ij} = \half \lambda_j ee^\top \hdrie \mbox{\rm if $j>i$}.\ee
Finally, the above list of possible $3\times 3$ principal submatrices of each vertex Gramian contains only nonblocking matrices. Therefore, also Corollary~\ref{Mco-9} proves the nonobtusity of the underlying simplices.\\[2mm]
The above vertex Gramians $G$ are all strictly ultrametric matrices. It may seem no surprise that the vertex Gramians of simplices whose vertex set is a finite ultrametric space with respect to the Euclidean metric, are strictly ultrametric matrices. But one should not be deceived by the convenient linguistics. Indeed, as argued in \cite{MaMiSa}, a matrix $A\geq 0$ is called strictly ultrametric if there exists an ultrametric $d$ on $I_n$ such that $a_{ij} \geq a_{ik}$ if and only if $d(i,j) \leq d(i,k)$ for all $i,j,k\in I_n$. The fact that the vertex Gramian is strictly ultrametric expresses that the mutual {\em inner products} between the edges of a simplex, whose {\em lengths} satisfy the ultrametric inequality, also satisfy a reversed version of the ultrametric inequality. This is a less trivial observation than the linguistics suggest, based on the fact that the three angles in a strongly isosceles triangle have no unique maximum. In the following section we show that the class of ultrametric matrices is much larger than the vertex Gramians of simplices with strongly isosceles triangular facets, and that it includes the {\em orthogonal} simplices.
\section{Tetrahedral facets of simplices}\label{Msect-4.2}\label{Msect-5}
We will now consider more general tetrahedral facets of a simplex $S$ that render $S$ itself nonobtuse. For this, we introduce in Section~\ref{Msect-5.1} the novel notion of {\em sub-orthocentric} tetrahedron in the classical context of metric geometry. Then in Section~\ref{Msect-5.2} we consider vertex Gramians of simplices with sub-orthocentric tetrahedral facets. These vertex Gramians turn out to form the simplicial matrix class of symmetric ultrametric matrices.
\subsection{The sub-orthocentric tetrahedron}\label{Msect-5.1}
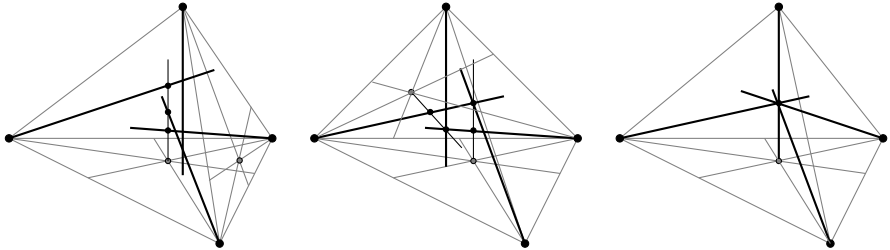
\begin{figure}[h]
\begin{center}  
\begin{tikzpicture}[scale=0.7, every node/.style={scale=0.7}]
\draw[gray] (0,2)--(5,2)--(4,0)--cycle;
\draw[gray] (0,2)--(4.667,1.333);
\draw[gray] (5,2)--(1.5,1.25);
\draw[gray] (4,0)--(2.75,2);
\draw[gray] (4,0)--(3.3,4.5)--(5,2);
\draw[gray] (3.3,4.5)--(4.55,1.1);
\draw[gray] (5,2)--(3.8,1.2);
\draw[gray] (4,0)--(4.6,2.6);
\draw[gray] (0,2)--(3.3,4.5);
\draw[fill=gray] (3.02,1.57) circle [radius=0.05];
\draw[fill=gray] (4.38,1.58) circle [radius=0.05];
\draw[fill=black] (0,2) circle [radius=0.07];
\draw[fill=black] (5,2) circle [radius=0.07];
\draw[fill=black] (4,0) circle [radius=0.07];
\draw (3.02,1.57)--(3.02,3.5);
\draw[thick] (5,2)--(2.3,2.2);
\draw[thick] (4,0)--(2.9,2.8);
\draw[thick] (0,2)--(3.9,3.3);
\draw[fill=black] (3.02,2.15) circle [radius=0.05];
\draw[fill=black] (3.02,2.5) circle [radius=0.05];
\draw[fill=black] (3.02,3) circle [radius=0.05];
\draw[fill=black] (3.3,4.5) circle [radius=0.07];
\draw[thick] (3.3,4.5)--(3.3,1.3);
\begin{scope}[shift={(5.8,0)}]
\draw[gray] (0,2)--(5,2)--(4,0)--cycle;
\draw[gray] (0,2)--(4.667,1.333);
\draw[gray] (5,2)--(1.5,1.25);
\draw[gray] (4,0)--(2.75,2);
\draw[gray] (4,0)--(2.5,4.5)--(5,2);
\draw[gray] (0,2)--(2.5,4.5);
\draw (3.02,1.57)--(3.02,3.5);
\draw[thick] (5,2)--(2.1,2.2);
\draw[thick] (4,0)--(2.77,3.33);
\draw[thick] (0,2)--(3.6,2.8);
\draw[fill=black] (3.02,2.15) circle [radius=0.05];
\draw[fill=black] (3.02,2.67) circle [radius=0.05];
\draw[thick] (2.5,4.5)--(2.5,1.46);
\draw[fill=black] (2.5,2.17) circle [radius=0.05];
\draw[fill=black] (2.2,2.5) circle [radius=0.05];
\draw[fill=gray] (3.02,1.57) circle [radius=0.05];
\draw[fill=gray] (1.84,2.88) circle [radius=0.05];
\draw (2.8,1.82)--(1.84,2.88);
\draw[gray] (0,2)--(3.4,3.6);
\draw[gray] (5,2)--(1.1,3.07);
\draw[gray] (2.5,4.5)--(1.5,2);
\draw[fill=black] (0,2) circle [radius=0.07];
\draw[fill=black] (5,2) circle [radius=0.07];
\draw[fill=black] (4,0) circle [radius=0.07];
\draw[fill=black] (2.5,4.5) circle [radius=0.07];
\end{scope}
\begin{scope}[shift={(11.6,0)}]
\draw[gray] (0,2)--(5,2)--(4,0)--cycle;
\draw[gray] (0,2)--(4.667,1.333);
\draw[gray] (5,2)--(1.5,1.25);
\draw[gray] (4,0)--(3.02,4.5);
\draw[gray] (3.02,4.5)--(5,2);
\draw[gray] (0,2)--(3.02,4.5);
\draw[fill=black] (0,2) circle [radius=0.07];
\draw[fill=black] (5,2) circle [radius=0.07];
\draw[fill=black] (4,0) circle [radius=0.07];
\draw[thick] (5,2)--(2.3,2.9);
\draw[thick] (4,0)--(2.9,2.93);
\draw[thick] (0,2)--(3.6,2.8);
\draw[fill=black] (3.02,2.67) circle [radius=0.05];
\draw[fill=black] (3.02,4.5) circle [radius=0.07];
\draw[thick] (3.02,4.5)--(3.02,1.57);
\draw[fill=gray] (3.02,1.57) circle [radius=0.05];
\draw[gray] (4,0)--(2.75,2);
\end{scope}
\end{tikzpicture}    
\end{center} 
\caption{\small{Left: a tetrahedron with four disjoint altitudes, each three of them intersecting an orthocentric perpendicular, and projecting onto the altitudes of a triangular facet. Middle: a semi-orthocentric tetrahedron with two pairs of intersecting altitudes. Right: an orthocentric tetrahedron with all altitudes intersecting in one point.}}
\label{Mfigure211}
\end{figure}

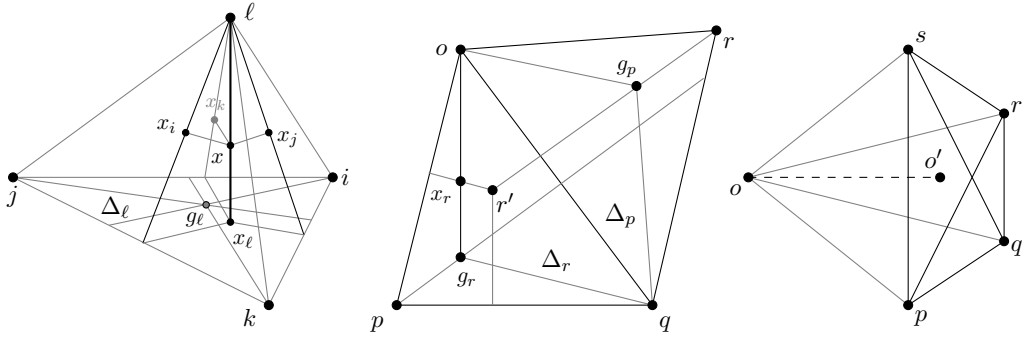
\begin{figure}[h]
\begin{center}  
\begin{tikzpicture}[scale=0.85, every node/.style={scale=0.85}]
\draw[gray] (0,2)--(5,2)--(4,0)--cycle;
\draw[gray] (4,0)--(3.4,4.5); 
\draw[gray] (3.4,4.5)--(5,2);
\draw[gray] (0,2)--(3.4,4.5);
\draw[gray] (5,2)--(1.5,1.25);
\draw[gray] (4,0)--(2.75,2);
\draw[gray] (0,2)--(4.667,1.333);
\draw[gray] (3.4,2.5)--(2.7,2.7);
\draw[gray] (3.4,2.5)--(4,2.7);
\draw (3.4,4.5)--(4.55,1.1);
\draw[gray] (4.55,1.1)--(3.4,1.3);
\draw (3.4,4.5)--(2.03,0.97);
\draw[gray] (2.03,0.97)--(3.4,1.3);
\draw[gray] (3.4,2.5)--(3.15,2.9);
\draw[gray] (3.4,4.5)--(3,2);
\draw[gray] (3,2)--(3.4,1.3);
\draw[fill=black] (0,2) circle [radius=0.07];
\draw[fill=black] (5,2) circle [radius=0.07];
\draw[fill=black] (4,0) circle [radius=0.07];
\draw[fill=black] (3.4,2.5) circle [radius=0.05];
\draw[fill=gray] (3.02,1.57) circle [radius=0.05];
\draw[fill=black] (3.4,1.3) circle [radius=0.05];
\draw[fill=black] (3.4,4.5) circle [radius=0.07];
\draw[fill=black] (2.7,2.7) circle [radius=0.05];
\draw[gray,fill=gray] (3.15,2.9) circle [radius=0.05];
\draw[fill=black] (4,2.7) circle [radius=0.05];
\draw[thick] (3.4,4.5)--(3.4,1.3);
\node at (3.7,4.6) {$\ell$};
\node at (0,1.7) {$j$};
\node at (3.7,-0.2) {$k$};
\node at (5.2,2) {$i$};
\node[scale=0.9] at (2.85,1.35) {$g_\ell$};
\node[scale=0.9] at (3.6,1) {$x_\ell$};
\node at (1.6,1.53) {$\Delta_\ell$};
\node[scale=0.9] at (3.2,2.3) {$x$};
\node[scale=0.9] at (2.4,2.8) {$x_i$};
\node[scale=0.9] at (4.3,2.6) {$x_j$};
\node[scale=0.8,gray] at (3.17,3.15) {$x_k$};
\begin{scope}[shift={(3,0)}]
\draw (3,0)--(7,0)--(4,4)--cycle;
\draw (7,0)--(8,4.3)--(4,4);
\draw[gray] (4,4)--(6.75,3.43);
\draw[gray] (8,4.3)--(4.5,1.8);
\draw[gray] (4.5,1.8)--(3.5,2.07);
\draw (4,4)--(4,0.75);
\draw[gray] (3,0)--(7.8,3.55);
\draw[gray] (7,0)--(6.75,3.43);
\draw[gray] (7,0)--(4,0.75);
\draw[gray] (4.5,1.8)--(4.5,0);
\node at (2.7,-0.3) {$p$};
\node at (7.2,-0.3) {$q$};
\node at (3.7,4) {$o$};
\node at (4.7,1.6) {$r'$};
\node at (8.2,4.1) {$r$};
\node[scale=0.9] at (4.1,0.4) {$g_r$};
\node[scale=0.9] at (3.7,1.7) {$x_r$};
\node[scale=0.9] at (6.6,3.7) {$g_p$};
\node at (5.5,0.7) {$\Delta_r$};
\node at (6.5,1.4) {$\Delta_p$};
\draw[fill=black] (3,0) circle [radius=0.07];
\draw[fill=black] (7,0) circle [radius=0.07];
\draw[fill=black] (4,4) circle [radius=0.07];
\draw[fill=black] (4.5,1.8) circle [radius=0.07];
\draw[fill=black] (8,4.3) circle [radius=0.07];
\draw[fill=black] (4,0.75) circle [radius=0.07];
\draw[fill=black] (4,1.94) circle [radius=0.07];
\draw[fill=black] (6.75,3.43) circle [radius=0.07];
\end{scope}
\begin{scope}[shift={(11.5,2)}]
\draw[gray] (0,0)--(2.5,-2);
\draw[gray] (0,0)--(2.5,2);
\draw[gray] (0,0)--(4,-1);
\draw[gray] (0,0)--(4,1);
\draw (2.5,-2)--(4,-1)--(4,1)--(2.5,2)--cycle;
\draw (4,-1)--(2.5,2);
\draw (4,1)--(2.5,-2);
\draw[fill=black] (0,0) circle [radius=0.07];
\draw[fill=black] (2.5,-2) circle [radius=0.07];
\draw[fill=black] (2.5,2) circle [radius=0.07];
\draw[fill=black] (4,-1) circle [radius=0.07];
\draw[fill=black] (4,1) circle [radius=0.07];
\node at (-0.2,-0.2) {$o$};
\node at (2.7,-2.2) {$p$};
\node at (2.7, 2.2) {$s$};
\node at (4.2,-1.1) {$q$};
\node at (4.2, 1.1) {$r$};
\node at (2.9,0.3) {$o'$};
\draw[dashed] (0,0)--(3,0);
\draw[fill=black] (3,0) circle [radius=0.07];
\end{scope}
\end{tikzpicture}    
\end{center} 
\caption{\small{Illustrations for the proofs of Lemma~\ref{Mlem-11} and Theorem~\ref{Mth-6}.}}
\label{Mfigure213}
\end{figure}

An {\em altitude} \index{altitude} of a simplex $S$ is a line segment falling from a vertex $v$ perpendicularly onto the hyperplane containing the facet $F$ opposite $v$. The three altitudes of a triangle $\Delta$ intersect in its {\em orthocenter} \index{orthocenter} $o$, which lies in $\Delta$ if and only if $\Delta$ is nonobtuse. Although the four altitudes of a tetrahedron $\tet$ generally do not intersect, they do satisfy some other long known properties, recently summarized and graphically well-illustrated in \cite{HaWe}. As in \cite{HaWe}, we call the line through the orthocenter of a triangular facet $\Delta$ of $\tet$ orthogonal to $\Delta$ its {\em orthocentric perpendicular}.
\begin{Pro}\label{Mpro-9} Let $\tet$ be a tetrahedron. Then:\\[2mm]
$\bullet$ the altitudes of the three vertices of a facet $\Delta$ of $\tet$ intersect its orthocentric perpendicular;\\[1mm]
$\bullet$ the projections of these three altitudes on $\Delta$ yield the (triangular) altitudes of $\Delta$;\\[1mm]
$\bullet$ if a pair of altitudes of $\tet$ have a common point, also the other pair has a common point;\\[1mm]
$\bullet$ such a common point lies on the orthocentric perpendicular of two triangular facets of $\tet$;\\[1mm]
$\bullet$ if three altitudes have a common point, then all four of them have this point in common.
\end{Pro}  
{\bf Proof. } See Figure~\ref{Mfigure211}, or \cite{HaWe} and the references therein.~\hfill $\Box$\\[2mm]
 A tetrahedron $\tet$ is known as {\em semi-orthocentric}~\cite{Cou2} \index{semi-orthocentric} if two altitudes of $\tet$ have a common point, and  {\em orthocentric} \cite{Cou1,Fie} if {\em all} its altitudes intersect in one point. Apart from the right tetrahedron in Figure~\ref{Mfigure211}, also the tetrahedra of type I with $m=k\leq \ell$ in Figure~\ref{Mfigure210} are orthocentric. The remaining tetrahedra of types I and II and the path tetrahedra, together with the middle tetrahedron in Figure~\ref{Mfigure211}, are only semi-orthocentric.\\[2mm] 
Whereas Figure~\ref{Mfigure26} refutes the possible conjecture that each semi-orthocentric tetrahedron with nonobtuse triangular facets is nonobtuse, there is a subclass of the semi-orthocentric tetrahedra with nonobtuse facets that does indeed have this property.
\begin{Def}[sub-orthocentric point] {\rm Any convex combination of a vertex and the orthocenter of a nonobtuse triangle $\Delta$ is called a {\em sub-ortho\-centric point} \index{sub-ortho\-centric point} 
of $\Delta$. We will write $\Delta_\ast$ for the set of sub-orthocentric points of~$\Delta$.}
\end{Def} 
\begin{Def}[sub-orthocentric tetrahedron] {\rm  A tetrahedron $\tet$ with nonobtuse triangular facets will be called {\em sub-orthocentric} \index{sub-orthocentric} if each altitude of $\tet$ meets a sub-orthocentic point on its opposite facet}.
\end{Def} 
\begin{rem}\label{Mrem-7}{\rm Any sub-orthocentric tetrahedron is semi-orthocentric. Moreover, since a sub-ortho\-centric tetrahedron $\tet$ has nonobtuse triangular facets by definition, each facet contains its sub-orthocentric points. In particular, due to Theorem~\ref{Mth-1}, $\tet$ is nonobtuse.}
\end{rem}
Now, the following proposition is an easy consequence of the third bullet in Proposition~\ref{Mpro-9}.
 \begin{Pro}\label{Mpro-11} If one altitude of a tetrahedron $\tet$ meets a sub-orthocentric point on its opposite facet, then each altitude of $\tet$ meets a sub-orthocentric point on its opposite facet. 
\end{Pro} 
{\bf Proof. } In Figure~\ref{Mfigure212}, the triangular facets $\Delta_r$ and $\Delta_q$ opposite $r$ and $q$ of a tetrahedron $\tet$ are unfolded onto the plane. For $\ell\in\{r,q\}$ we indicate by $\ell'$ the projection of $\ell$ onto $\Delta_\ell$. Obviously, $\ol{or'}\perp\ol{pq}$ if and only if $\ol{or}\perp\ol{pq'}$. Since $\ol{rr'}$ and $\ol{qq'}$ are both perpendicular to $\ol{op}$, we see that $r'$ is a sub-orthocentric point of $\Delta_r$ if and only if $q'$ is a sub-orthocentric point of $\Delta_q$. The corresponding result for the remaining two facets follows similarly.~\hfill $\Box$
\begin{figure}[h]
\begin{center}  
\begin{tikzpicture}[scale=0.8, every node/.style={scale=0.8}]
\draw (3,0)--(7,0)--(4,4)--cycle;
\draw[gray] (3,0)--(0,3)--(4,4);
\draw (4,4)--(4,0);
\draw[gray] (3,0)--(5.55,1.95);
\draw[gray] (0,3)--(4,1.8);
\draw[gray] (7,0)--(2.65,1.1);
\draw[gray] (4,4)--(1.5,1.5);
\draw (3,0)--(2,3.5);
\node at (2.7,-0.3) {$p$};
\node at (7.2,-0.3) {$q$};
\node at (4.3,3.9) {$o$};
\node at (4.3,2) {$r'$};
\node at (0,3.3) {$r$};
\node at (2.9,1.3) {$q'$};
\draw[fill=black] (3,0) circle [radius=0.07];
\draw[fill=black] (7,0) circle [radius=0.07];
\draw[fill=black] (4,4) circle [radius=0.07];
\draw[fill=black] (4,1.8) circle [radius=0.07];
\draw[fill=black] (2.67,1.1) circle [radius=0.07];
\draw[fill=black] (0,3) circle [radius=0.07];
\end{tikzpicture}    
\end{center} 
\caption{\small{If $r$ projects on a sub-orthocentric point of its opposite facet, so does $q$.}}
\label{Mfigure212}
\end{figure}\\[2mm]
We are now able to prove a remarkable lemma, that will result in the observation in Theorem~\ref{Mth-6} that if a $4$-simplex $\four$ has four sub-orthocentric tetrahedral facets, also its fifth facet will be sub-orthocentric. This generalizes the corresponding property of equilateral triangles to the newly defined class of sub-orthocentric tetrahedra. 
\begin{Le}\label{Mlem-11} A tetrahedron $\tet$ with nonobtuse facets is sub-orthocentric if and only if there exists a point $x$ whose orthogonal projections on the facets of $\tet$ are all sub-orthocentric.
\end{Le} 
{\bf Proof. } Let $\tet$ be a tetrahedron with vertices $opqr$ and let $x$ be given. For $\ell\in\{o,p,q,r\}$, write $x_\ell$ for the projection of $x$ on the triangular facet $\Delta_\ell$ of $\tet$ opposite $\ell$. Denote the orthocenter of $\Delta_\ell$ by $g_\ell$. We will distinguish two cases. The left picture in Figure~\ref{Mfigure213} depicts the first case, in which $x_i,x_j\in\{x_o,x_p,x_q,x_r\}$ lie on altitudes that meet at a vertex $\ell$. Then the second bullet of Proposition~\ref{Mpro-9} implies that $x$ lies on the tetrahedral altitude of $\tet$ from $\ell$. Thus, $x$ also projects on the altitude from $\ell$ on the third triangular facet meeting at $\ell$. Proposition~\ref{Mpro-11} now shows that if $\tet$ is not sub-orthocentric, then $x_\ell$ is not a sub-orthocentric point of $\Delta_\ell$. Conversely, if $\tet$ is sub-orthocentric, then $x_\ell$ is a sub-orthocentric point of $\Delta_\ell$ and moreover, by choosing $x$ close enough to, or even equal to $\ell$, all four points $\{x_i,x_j,x_k,x_\ell\}$ are orthocentric on their respective triangular facets. This proves the statement for the first case. The second case, depicted as unfolded view in the middle of Figure~\ref{Mfigure213}, concerns the complementary situation that each of the projections $x_o,x_p,x_q,x_r$ lies on a triangular altitude meeting a different vertex of $\tet$. Assume that $\tet$ is not sub-orthocentic. Without loss of generality, assume that $x_r\in\ol{og_r}\subset\Delta_r$ and that the projection $r'$ of $r$ on $\Delta_r$ lies right of the altitude from $o$. Then $x_o\in\Delta_o$ can only lie on the altitude from $p$ in $\Delta_o$, and consequently, $x_q$ must lie on the altitude from $r$ in $\Delta_q$. For this to be possible, $r'$ needs to lie in the interior of triangle $oqg_r$. This however establishes the position of $x$ as lying above the line $\ol{r'r}$. But then $x_p$ cannot lie between $q$ and $g_p$ in $\Delta_p$. Thus, if $\tet$ is not a sub-orthocentric tetrahedron, $x_o,x_p,x_q$, and $x_r$ cannot all be sub-orthocentric points.~\hfill $\Box$
\begin{Th}\label{Mth-6} Let $\four$ be a $4$-simplex. If $\four$ has four sub-orthocentric tetrahedral facets, also its fifth tetrahedral facet is sub-orthocentric.
\end{Th}
{\bf Proof. } Let $\four$ be a $4$-simplex with vertices $o,p,q,r,s$. Assume that each tetrahedral facet of $\four$ containing $o$ is sub-orthocentric. Then each triangular facet of $\four$ is nonobtuse, including those of the tetrahedral facet $\tet$ of $\four$ opposite $o$. Moreover, since $o$ projects on a sub-orthocentric point of each triangular facet $\Delta$ of $\tet$, so does $o'$, where $o'$ is the projection of $o$ into the hyperplane containing $\tet$. Lemma~\ref{Mlem-11} now shows that $\tet$ is sub-orthocentric.~\hfill $\Box$\\[2mm] 
Theorem~\ref{Mth-6}, which is illustrated in the right picture in Figure~\ref{Mfigure213} can be proved also using linear algebra only, but the geometric proof given here fits better in the context of this thesis.
\\[2mm]
Proposition~\ref{Mpro-11} and Theorem~\ref{Mth-6} will have interesting consequences for vertex Gramians of simplices with sub-orthocentric tetrahedral facets. This will be discussed in the next section.
\subsection{Simplices with sub-orthocentric tetrahedral facets}\label{Msect-5.2}
Consider now the vertex Gramian $G_o$ of the sub-orthocentric tetrahedron $\tet$ in Figure~\ref{Mfigure212}. With a minor abuse of notation, write $p,r,q$ for the {\em position vectors} with origin $o$ of its vertices with the corresponding labels. As all triangular facets of $\tet$ are nonobtuse, Lemma~\ref{Mpro-15} shows that $G_o\geq 0$ and that each diagonal entry is maximal in its row/column. Moreover, as $\tet$ is semi-orthocentric, one of the vectors $p,q,r$ is orthogonal to the difference of the other two. In Figure~\ref{Mfigure212}, $r\perp p-q$, or equivalently, $r^\top p=r^\top q$. Moreover, as $r$ projects between $g_r$ and $o$, we have that $r^\top p=r^\top q\leq p^\top q$. Thus, $G_o$ is a nonblocking matrix. Proposition~\ref{Mpro-11} shows that the same is true for {\em all} vertex Gramians of $\tet$, as each vertex of $\tet$ is incident to an edge orthogonal to its opposite edge, projecting on a sub-orthocentric point of its opposite facet.
\begin{Co}\label{Mlem-1} Let $\tet$ be a tetrahedron. Then the following are equivalent:\\[2mm]
$\bullet$ $\tet$ is sub-orthocentric;\\[2mm]
$\bullet$ $\tet$ has a nonblocking vertex Gramian;\\[2mm]
$\bullet$ all vertex Gramians of $\tet$ are nonblocking.
\end{Co}     
{\bf Proof. } This is an immediate consequence of Proposition~\ref{Mpro-11}.~\hfill $\Box$\\[2mm]
In fact, we have just shown that the property of being nonblocking is invariant under $\sim$ for all $A\in\RR^{3\times 3}_{\rm spd}$, as was already announced in Section~\ref{Msect-3.3}.\\[2mm]
Recall that a nonnegative $n\times n$ matrix $A\in\Rnspd$ is called {\em symmetric ultrametric} \index{symmetric ultrametric} \cite{VarNab} if each $3\times 3$ principal submatrix of $A$ has no unique minimal entry above its diagonal, and each diagonal element is maximal in its row. Writing from now on {\em ultrametric} instead of {\em symmetric ultrametric}, in the language of this chapter this reads as follows.
 \begin{Pro}\label{Mpro-13} Let $A=(a_{ij})\in\RR^{3\times 3}_{\rm spd}$. Then $A$ is ultrametric if and only if $A$ is nonblocking.
\end{Pro}
{\bf Proof. } If $\{a_{12},a_{13},a_{23}\}$ has no unique minimum, at least four of the six off-diagonal entries of $A$ are minimal and symmetrically placed. Thus, each column of $A$ contains an off-diagonal entry that is minimal in its row. Conversely, if this minimum is unique, it may without loss of generality be assumed to be the entry $a_{12}$. As $a_{21}=a_{12}$, this shows that the third column of $A$ is blocking.~\hfill $\Box$\\[2mm]
Corollary~\ref{Mlem-1} establishes ultrametricity of a $3\times 3$ matrix $A$ as a truly {\em geometric} property of its underlying sub-orthocentric tetrahedron $\SIS(A)$, and the $3\times 3$ ultrametric matrices as a simplicial matrix class as in Definition~\ref{Msmc}: either none, or all vertex Gramians of a tetrahedron are ultrametric. The matrix $G_0$ in (\ref{MG0}) shows that for $n\geq 4$, there exist nonblocking $A\in\Rnspd$ that are not ultrametric, and thus the validity of Proposition~\ref{Mpro-13} is restricted to $\RR^{3\times 3}_{\rm spd}$.
\begin{Th}\label{Mth-5} Let $n\geq 3$ and let $S$ be an $n$-simplex. Then equivalent are:\\[2mm] 
$\bullet$ $S$ has a vertex $v$ such that all tetrahedral facets $\tet$ with $v\in\tet$ are sub-orthocentric;\\[2mm]
$\bullet$ all tetrahedral facets of $S$ are sub-orthocentric;\\[2mm]
$\bullet$ $S$ has an ultrametric vertex Gramian;\\[2mm]
$\bullet$ all vertex Gramians of $S$ are ultrametric.
\end{Th}  
{\bf Proof. } Let $v\in S$ and assume that all tetrahedral facets $\tet$ with $v\in\tet$ are sub-orthocentric. Let $\tet$ be a tetrahedral face of $S$ with $v\not\in\tet$. The convex hull of $v$ and $\tet$ is a $4$-simplex $\four$ having four sub-orthocentric facets containing $v$. Now, Theorem~\ref{Mth-6} proves that also $\tet$ is sub-orthocentric. Thus, all tetrahedral facets of $S$ are sub-orthocentric. Due to Corollary~\ref{Mlem-1} and Proposition~\ref{Mpro-13}, the third and fourth bullet are equivalent to the first and the second.~\hfill $\Box$\\[2mm]
Theorem~\ref{Mth-5} shows that the set $\UU$ of all symmetric ultrametric $n\times n$ matrices is a simplicial matrix class. In \cite{MaMiSa}, it was proved that each strictly ultrametric matrix $A\in\Rnspd$ has a diagonally dominant Stieltjes matrix as inverse using concepts from discrete measure theory. In \cite{NabVar} a much simpler linear algebraic proof was given for the same result. We are now able to finish a geometric proof. We do not claim it to be simpler than the proof in \cite{NabVar}, especially not after all previously introduced concepts. But we do claim that the geometric ideas presented so far gives useful insights in inverse M-matrices that are complementary to the purely algebraic view, which is dominant in the literature.
\begin{Th}\label{Mco-8} Let $A\in\Rnspd$ be an ultrametric. Then $A\in\MM_\dd$.
\end{Th}
{\bf Proof. } Due to Theorem~\ref{Mth-5}, the underlying simplex $\SIS(A)$ of $A$ has only sub-orthocentric tetrahedral facets. Thus, it has nonobtuse triangular facets, and in particular, all $3\times 3$ principal submatrices of all vertex Gramians of $S$ are nonblocking. Hence, Corollary~\ref{Mco-9} proves the statement.~\hfill $\Box$
 \subsection{Tetrahedral facets of an orthogonal simplex}\label{Msect-5.3}
The path-simplices from Section~\ref{Msect-4} belong to the important class \cite{BrKoKrSo} of {\em orthogonal simplices}, \index{orthogonal simplices} $n$-simplices whose vertex-edge graph has a spanning tree of edges that are mutually orthogonal. They can also be characterized as simplices having the maximal amount of $r(n) = \binom{\npo}{2}-n$ right dihedral angles, taking into account a theorem by Fiedler~\cite{Fie} that states that each $n$-simplex has at least $n$ acute dihedral angles.
\begin{rem}{\rm There are two types of orthogonal {\em tetrahedra}. A {\em path tetrahedron} has a path of orthogonal edges and a {\em cube corner} has three orthogonal edges meeting at the same vertex}.
\end{rem}
Contrary to path-simplices, facets of orthogonal simplices are generally not orthogonal, as any cube corner shows. They are, however, all sub-orthocentric, as we will prove below. Reason is that each orthogonal simplex has at least two orthogonal facets. This follows from the fact that each tree has at least two leaves. Hence, the facet opposite a vertex that is a leaf of its spanning tree of orthogonal edges, obviously has a spanning tree of orthogonal edges, and is thus an orthogonal $(\nmo)$-simplex.
\begin{Th} Each tetrahedral facet of an orthogonal $n$-simplex $S$ is sub-or\-tho\-cen\-ter, and each vertex Gramian of $S$ is ultrametric.
\end{Th}
{\bf Proof. }  Let $\tet$ be an orthogonal tetrahedron. Then $\tet$ is a path-tetrahedron or a cube corner, and thus is $\tet$ sub-orthocentric. This proves the statement for $n=3$. As inductive hypotheses, assume the claim for all orthogonal $(\nmo)$-simplices. Let $S$ be an orthogonal $n$-simplex. Consider a vertex $u$ of $S$ that is a leaf of the spanning tree $T$ of orthogonal edges of $S$. Then the altitude of $S$ from $u$ meets a vertex (its parent in $T$) $v$ of the facet $F$ opposite $u$, which is an orthogonal $(\nmo)$-simplex. Now, let $\tet$ be a tetrahedral facet of $S$ with $v\in\tet$. If $\tet\subset F$ then $\tet$ is sub-orthocentric due to the inductive hypothesis. If $\tet\not\subset F$ then $u$ is also a vertex of $\tet$ and its altitude meets the sub-orthocentric point $v$ of its opposite triangular facet $\Delta$. As $\Delta$ is nonobtuse and $u$ projects on $v$, also $\tet$ is nonobtuse. Thus, $\tet$ is sub-orthocentric as a result of Proposition~\ref{Mpro-11}. This shows that {\em all} tetrahedral facets of $S$ containing $v$ are sub-orthocentric, and Theorem~\ref{Mth-5} now proves both the statements.~\hfill $\Box$\\[2mm]
An {\em $n$-cube corner simplex} is an orthogonal simplex whose spanning tree $T$ of orthogonal edges is the tree on $\npo$ vertices with $n$ leaves. If an orthogonal simplex $S$ is not an $n$-cube corner, then $T$ has a path of length three, or in other words, $S$ has a tetrahedral facet that is a path simplex. As orthocentric simplices have orthocentric facets, this proves the following.
\begin{Pro} Let $S$ be an orthogonal simplex. Then $S$ is orthocentric if and only if $S$ is an $n$-cube corner.\end{Pro} 
\section{Special $4$-simplicial facets}\label{Msect-6}
In the previous sections, we saw the following. Let $S$ be an $n$-simplex and $v$ a vertex of $S$. If all {\em triangular} facets $\Delta$ of an $n$-simplex $S$ with $v\in\Delta$ are equilateral, then {\em all} its triangular facets are equilateral and it followed that $S$ is nonobtuse. If all {\em tetrahedral} facets $\tet$ of $S$ with $v\in \tet$ are sub-orthocentric, then {\em all} tetrahedral facets of $S$ are sub-orthocentric and it followed that $S$ is nonobtuse. Furthermore, each nonobtuse triangle $\Delta$ is a facet of a sub-orthocentric tetrahedron $\tet$, and trivially, each (nonobtuse) line-segment is the edge of an equilateral triangle. A logical next step, to be investigated in this section, is to look for a class of $4$-simplicial facets of $S$ guaranteeing nonobtusity of $S$ in a similar manner.
\subsection{Some considerations in special $4$-simplex}
As a direct generalization of the equilateral triangle and the sub-orthocentric tetrahedron, we would like to identify a class $\FF$ of $4$-simplices such that:\\[2mm]
$\bullet$ each tetrahedral facet $\tet$ of a simplex $\four\in\FF$ is nonobtuse;\\[2mm]
$\bullet$ any nonobtuse tetrahedron is the facet of a $4$-simplex $\four\in\FF$; \\[2mm]
$\bullet$ if $S$ is a $5$-simplex with five of its facets in $\FF$, then also its sixth facet is in~$\FF$;\\[2mm]
$\bullet$ if $S$ is an $n$-simplex with $n \geq 5$, all whose $4$-facets are in $\FF$, then $S$ is nonobtuse.\\[2mm]
Note that the third bullet in this list implies, in the same way as in Proposition~\ref{Mpro-8} and Theorem~\ref{Mth-6}, that if $S$ is an $n$-simplex and $v$ a vertex of $S$, and all $4$-facets $\four$ of $S$ with $v\in \four$ are in $\FF$, then {\em all} $4$-facets of $S$ are in $\FF$. The second bullet implies that $\FF$ will be larger than simply the class of $4$-simplices with sub-orthocentric tetrahedral facets only, which has all the above properties. In matrix terms, the above boils down to: given $A\in\Rnspd$, we would like to identify conditions on each of the $4\times 4$ principal submatrices of $A$, such that the set of all $A\in\Rnspd$ satisfying those conditions is a simplicial matrix subclass $\CC$ of $\MM_\dd$. As $A\in\Rnspd$ has $\binom{n}{4}$ principal $4\times 4$ submatrices, any non-empty set of conditions on $4\times 4$ principal submatrices that guarantees $A$ to be in $\MM_\four$ can not be verified in less than $\mathcal{O}(n^4)$ operations.\\[2mm]
Since inverting $A$ to verify if $A^{-1}$ is an M-matrix can be done in $\mathcal{O}(n^3)$ operations, our investigations are mainly of theoretical interest. To start, we recall Theorem~\ref{Mth-8} with $k=3$, which may provide a tool towards the required class of $4$-simplices.
 \begin{Th}\label{Mth-11} Let $A\in\Rnspd$ be given. Suppose that all tetrahedral facets of the underlying simplex $S=\SIS(A)$ of $A$ are nonobtuse. Moreover, assume that for all $m\geq 4$, each $m\times m$ principal submatrix of each vertex Gramian of $S$ is nonblocking. Then $A\in\MM_\dd$.
\end{Th}
This theorem encourages to identify $4$-simplices with nonobtuse tetrahedral facets with nonblocking vertex Gramians only. More specifically, we would like that if four of its vertex Gramians are nonblocking, also its fifth is nonblocking. Additionally, this should guarantee that also all larger principal submatrices of all vertex Gramians are nonblocking.\\[2mm]
A complication in this approach is that we have seen in Section~\ref{Msect-3.3} that in general, not all vertex Gramians of a $4$-simplex with nonobtuse tetrahedral facets need to be nonblocking simultaneously. Thus in order to use the above theorem, a more restrictive property on the candidate $4$-simplices will be necessary. Moreover, as we saw in (\ref{Mexam-1}), even if all $4\times 4$ principal submatrices of a matrix $A\in\Rnspd$ are nonblocking, this does not guarantee that {\em larger} submatrices are also nonblocking. In the next section we will prove the negative but nonetheless interesting result that the additional restrictions can only be met by $4$-simplices having sub-orthocentric tetrahedral facets only. Thus, Theorem~\ref{Mth-11} cannot be used to find a class $\FF$ of $4$-simplices satisfying the second bullet in the list of requirements above. It only yields the $4$-simplices with sub-orthocentric tetrahedral facets discussed in Section~\ref{Msect-5} already.
\subsection{$4$-simplices with nonblocking vertex Gramians only}
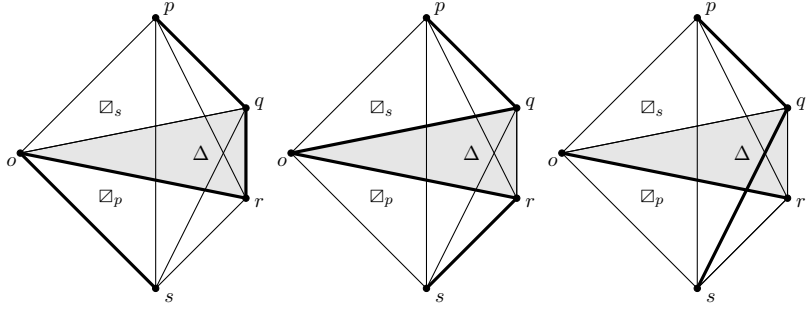
\begin{figure}[h]
\begin{center}  
\begin{tikzpicture}[scale=0.6, every node/.style={scale=0.7}]
\draw[fill=gray!20!white] (0,0)--(5,1)--(5,-1)--cycle;
\draw[very thick] (0,0)--(3,-3);
\draw (0,0)--(3,3);
\draw[very thick] (0,0)--(5,-1);
\draw (0,0)--(5,1);
\draw (3,-3)--(5,-1)--(5,1)--(3,3)--cycle;
\draw (5,-1)--(3,3);
\draw (5,1)--(3,-3);
\draw[very thick] (3,3)--(5,1);
\draw[very thick] (5,-1)--(5,1);
\draw[fill=black] (0,0) circle [radius=0.07];
\draw[fill=black] (3,-3) circle [radius=0.07];
\draw[fill=black] (3,3) circle [radius=0.07];
\draw[fill=black] (5,-1) circle [radius=0.07];
\draw[fill=black] (5,1) circle [radius=0.07];
\node at (-0.2,-0.2) {$o$};
\node at (3.3,-3.2) {$s$};
\node at (3.3, 3.2) {$p$};
\node at (5.3,-1.1) {$r$};
\node at (5.3, 1.1) {$q$};
\node at (4,0) {$\Delta$};
\node at (2,1) {$\boxslash_s$};
\node at (2,-1) {$\boxslash_p$}; \ 
\begin{scope}[shift={(6,0)}]
\draw[fill=gray!20!white] (0,0)--(5,1)--(5,-1)--cycle;
\draw (0,0)--(3,-3);
\draw[very thick] (5,-1)--(3,-3);
\draw (0,0)--(3,3);
\draw[very thick] (0,0)--(5,-1);
\draw (0,0)--(5,1);
\draw (3,-3)--(5,-1)--(5,1)--(3,3)--cycle;
\draw (5,-1)--(3,3);
\draw (5,1)--(3,-3);
\draw[very thick] (3,3)--(5,1);
\draw[very thick] (0,0)--(5,1);
\draw[fill=black] (0,0) circle [radius=0.07];
\draw[fill=black] (3,-3) circle [radius=0.07];
\draw[fill=black] (3,3) circle [radius=0.07];
\draw[fill=black] (5,-1) circle [radius=0.07];
\draw[fill=black] (5,1) circle [radius=0.07];
\node at (-0.2,-0.2) {$o$};
\node at (3.3,-3.2) {$s$};
\node at (3.3, 3.2) {$p$};
\node at (5.3,-1.1) {$r$};
\node at (5.3, 1.1) {$q$};
\node at (4,0) {$\Delta$};
\node at (2,1) {$\boxslash_s$};
\node at (2,-1) {$\boxslash_p$};
\end{scope} \ 
\begin{scope}[shift={(12,0)}]
\draw[fill=gray!20!white] (0,0)--(5,1)--(5,-1)--cycle;
\draw (0,0)--(3,-3);
\draw (5,-1)--(3,-3);
\draw (0,0)--(3,3);
\draw[very thick] (0,0)--(5,-1);
\draw (0,0)--(5,1);
\draw (3,-3)--(5,-1)--(5,1)--(3,3)--cycle;
\draw (5,-1)--(3,3);
\draw (5,1)--(3,-3);
\draw[very thick] (3,3)--(5,1);
\draw[very thick] (3,-3)--(5,1);
\draw[fill=black] (0,0) circle [radius=0.07];
\draw[fill=black] (3,-3) circle [radius=0.07];
\draw[fill=black] (3,3) circle [radius=0.07];
\draw[fill=black] (5,-1) circle [radius=0.07];
\draw[fill=black] (5,1) circle [radius=0.07];
\node at (-0.2,-0.2) {$o$};
\node at (3.3,-3.2) {$s$};
\node at (3.3, 3.2) {$p$};
\node at (5.3,-1.1) {$r$};
\node at (5.3, 1.1) {$q$};
\node at (4,0) {$\Delta$};
\node at (2,1) {$\boxslash_s$};
\node at (2,-1) {$\boxslash_p$};
\end{scope}
\end{tikzpicture}    
\end{center} 
\caption{\small{Possible relative configurations of blocking edges of $\tet_p$ and $\tet_s$.}}
\label{Mfigure214}
\end{figure}

First we will refine some observations on tetrahedra that will be of interest for what is to come. For this we will say that a vertex $x$ {\em causes a blocking column} in $G_y$ if the column of $G_y$ corresponding to the vertex $x$, or equivalently, to the position vector $x-y$, is blocking.
\begin{Le}\label{Mlem-12} Let $\tet$ be a tetrahedron with nonobtuse triangular facets. Denote the vertices of $\tet$ by $o,p,q,r$ and their respective vertex Gramians by $G_o,G_p,G_q,G_r$. Then each vertex Gramian of $\tet$ has at most one blocking column. 
In particular, the following are equivalent:\\[2mm]
$\bullet$ $p$ causes a blocking column in $G_o$;\\[2mm]
$\bullet$ $o$ causes a blocking column in $G_p$;\\[2mm]
$\bullet$ $q$ causes a blocking column in $G_r$;\\[2mm]
$\bullet$ $r$ causes a blocking column in $G_q$.
\end{Le}  
{\bf Proof. } Suppose a vertex Gramian of $\tet$ has a blocking column. Without loss of generality, assume it is its third column, as in Figure~\ref{Mfigure27}. Then the entries at position $(1,2)$ and $(2,1)$ are minimal in their respective columns. This proves the first statement. The equivalence of the third and fourth bullet is a slight extension of Proposition~\ref{Mpro-11} and Figure~\ref{Mfigure212}. In Figure~\ref{Mfigure212}, neither $q$ nor $r$ causes a blocking column in $G_r$ and $G_q$ as both altitudes from $q$ and $r$ meet a sub-orthocentric point on their opposite facets. However, shifting $r'$ slightly to its left will make $q'$ shift to its right, and vice versa. This proves the equivalence of the third and fourth bullet. Similarly, the first and second are equivalent to one another. They are also equivalent to the third and fourth, because in Corollary~\ref{Mlem-1} we proved that if one vertex Gramian has a blocking column, all of them have, combined with the fact that we have just argued that each vertex Gramian can have at most one blocking column.~\hfill $\Box$
\begin{rem}{\rm  An algebraic proof of the equivalence of the four bullets in the above lemma can be obtained by explicitly writing out all four vertex Gramians and verifying (in-)equalities. This is rather straightforward and provides solid proof, but it is not in the style of this thesis.}
\end{rem} 
The consequence of Lemma~\ref{Mlem-12} is, that to each tetrahedron $\tet$ with nonobtuse facets that is {\em not} sub-orthocentric, we can assign two edges that are its {\em blocking edges}. These edges, which are opposite one another in $\tet$ are the ones that cause the blocking columns in the vertex Gramians of $\tet$. Using this, we are now able to prove the following important negative result.
\begin{Th}\label{Mth-12} Let $\four$ be a $4$-simplex with nonobtuse tetrahedral facets. Then equivalent are:\\[2mm]
$\bullet$ each vertex Gramian of $\four$ is nonblocking;\\[2mm]
$\bullet$ each vertex Gramian is ultrametric.
\end{Th} 
{\bf Proof. } The implication from the second to the first bullet has already been proved. Indeed, if each vertex Gramian is ultrametric, then all tetrahedral facets of $\four$ are sub-orthocentric due to Theorem~\ref{Mth-5}. Then  Corollary~\ref{Mlem-1} shows that all $3\times 3$ submatrices of all vertex Gramians of $\four$ are nonblocking and Corollary~\ref{Mcor-x} proves that each vertex Gramian is nonblocking. To prove the converse implication, we will show if $\four$ has a tetrahedral facet that is {\em not} sub-orthocentric, which is equivalent to having a vertex Gramian that is not ultrametric, then {\em at least one} vertex Gramian of $\four$ has a blocking column. Let $o,p,q,r,s$ be the vertices of $\four$ and $G_o,\dots,G_s$ its corresponding vertex Gramians. Assume that the tetrahedron $\tet_s$ opposite $s$ with vertices $o,p,q,r$ is not sub-orthocentric. Theorem~\ref{Mth-6} shows that $\four$ must have another tetrahedral facet that is not sub-orthocentric. Without loss of generality, assume it is $\tet_p$ opposite $p$, with vertices $o,q,r,s$. Note that $\tet_p$ and $\tet_s$ share the triangular facet $\Delta$ with vertices $o,p,q$. See Figure~\ref{Mfigure214}, in which $\Delta$ is shaded. As $\tet_s$ is not sub-orthocentric, it has a pair of opposite blocking edges. We may without loss of generality assume that these are $\ol{or}$ and $\ol{pq}$. In Figure~\ref{Mfigure214}, these edges have been made bold. Now, also $\tet_p$ has a pair of opposite blocking edges, which, of course, we may not fix without loss of generality. Instead, we display all three possible options in Figure~\ref{Mfigure214}, again by drawing the pair of opposite blocking edges in bold. In the right picture, $\ol{or}$ is blocking in both $\tet_s$ and $\tet_p$.
We will now prove that each of the three configurations of blocking edges in Figure~\ref{Mfigure214} leads to a blocking column in a vertex Gramian of $\four$. The first claim is that for both the left and the right configuration, the vertex Gramian $G_o$ corresponding to $o$ has a blocking column. Indeed, considering $p,q,r,s$ as position vectors with respect to the origin $o$, we have that
\[ G_o = \left[\begin{array}{rrrr} p^\top p & p^\top q & p^\top r & p^\top s \\
q^\top p & q^\top q & q^\top r & q^\top s \\
r^\top p & r^\top q & r^\top r & r^\top s \\
s^\top p & s^\top q & s^\top r & s^\top s \end{array}\right] = \left[\begin{array}{rrr|r}  &  & & p^\top s \\
 & G_o^s & & q^\top s \\ 
 &  &  & r^\top s \\ \hline
s^\top p & s^\top q & s^\top r & s^\top s \end{array}\right] = \left[\begin{array}{r|rrr} p^\top p & p^\top q & p^\top r & p^\top s \\ \hline
q^\top p &  &  &  \\
r^\top p &  & G_o^p &  \\
s^\top p &  & & \end{array}\right]. \]
Since $\ol{or}$ is blocking in $\tet_s$, both $p^\top r$ and $q^\top r$ are not minimal in their row in the top-left $3\times 3$ principal submatrix $G_o^s$ of $G_o$. Thus, trivially, they can also not be minimal in their row in $G_o$, regardless of the configuration of blocking edges in $\tet_p$. Consider now first the right configuration in Figure~\ref{Mfigure214}. Then, as $\ol{or}$ is also blocking in the bottom-right $3\times 3$ principal submatrix $G_o^p$ of $G_o$, also $s^\top r$ (and again $q^\top r$) cannot be minimal in its row in $G_o^p$, nor in $G_o$. This proves that the third column of $G_o$ is blocking in the right configuration. In the left configuration, the last column of $G_o^p$ is blocking, showing that $q^\top s$ and $r^\top s$ are not minimal in their row in $G_o^p$, nor in $G_o$. Thus, either the last column of $G_o$ is blocking, or the entry $p^\top s$ is minimal in its row. But in the latter case, $p^\top s \leq p^\top q < q^\top r$, whereas also $q^\top r<q^\top s$ and $q^\top r < r^\top s$. This shows that $s^\top p = p^\top s < s^\top r$, hence $s^\top r$ cannot be minimal in its (fourth) row of $G_o$, and the third column of $G_o$ is blocking. Finally, to show that also the middle configuration leads to a blocking column, observe that the vertex Gramian $G_r$ corresponding to the vertex $r$ of $\four$ also has a blocking edge in $\tet_s$, and a blocking edge in $\tet_p$ that is {\em not} an edge of $\tet_s$. It is exactly that property that gave rise to the blocking column of the vertex Gramian $G_o$ in the left configuration. In other words, choose the position vectors from $r$ to $p,q,o,s$ {\em in that order} to form the Gramian $G_r$, and repeat the arguments for the left configuration.~\hfill $\Box$\\[2mm]
As mentioned above, this result can be called a negative result: only the simplices $S$ with sub-orthocentric tetrahedral facets have $4$-facets that have only nonblocking vertex Gramians, and thus no new class of $4$-simplices is found using the concepts of blocking column and nonblocking matrix. Note that the proof of Theorem~\ref{Mth-12} is mainly combinatorial in nature.
\subsection{The geometry of blocking columns}\label{Msect-6.3}
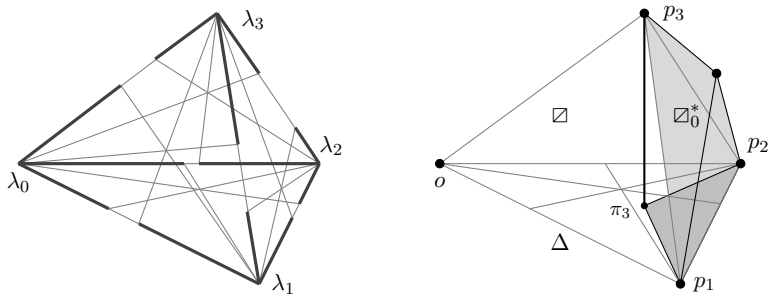
\begin{figure}[h]
\begin{center}  
\begin{tikzpicture}[scale=0.8, every node/.style={scale=0.8}]
\draw[gray] (0,2)--(5,2)--(4,0)--cycle;
\draw[gray] (0,2)--(4.667,1.333);
\draw[gray] (5,2)--(1.5,1.25);
\draw[gray] (4,0)--(2.75,2);
\draw[gray] (4,0)--(3.3,4.5)--(5,2);
\draw[gray] (3.3,4.5)--(4.55,1.1);
\draw[gray] (5,2)--(3.8,1.2);
\draw[gray] (4,0)--(4.6,2.6);
\draw[gray] (0,2)--(3.3,4.5);
\draw[gray] (3.3,4.5)--(3,2);
\draw[gray] (3.3,4.5)--(2,1);
\draw[gray] (4,0)--(1.7,3.3);
\draw[gray] (0,2)--(3.65,2.32);
\draw[gray] (0,2)--(4,3.5);
\draw[gray] (5,2)--(2.27,3.74);
\draw[very thick,gray!50!black] (0,2)--(1.7,3.3);
\draw[very thick,gray!50!black] (0,2)--(1.5,1.25);
\draw[very thick,gray!50!black] (0,2)--(2.75,2);
\draw[very thick,gray!50!black] (3.3,4.5)--(3.65,2.32);
\draw[very thick,gray!50!black] (3.3,4.5)--(4,3.5);
\draw[very thick,gray!50!black] (3.3,4.5)--(2.27,3.74);
\draw[very thick,gray!50!black] (5,2)--(4.6,2.6);
\draw[very thick,gray!50!black] (5,2)--(3,2);
\draw[very thick,gray!50!black] (5,2)--(4.667,1.333);
\draw[very thick,gray!50!black] (4,0)--(3.8,1.2);
\draw[very thick,gray!50!black] (4,0)--(2,1);
\draw[very thick,gray!50!black] (4,0)--(4.55,1.1);
\node at (0,1.7) {$\lambda_0$};
\node at (4.4,0) {$\lambda_1$};
\node at (5.2,2.3) {$\lambda_2$};
\node at (3.9,4.4) {$\lambda_3$};
\begin{scope}[shift={(7,0)}]
\draw[gray!30!white,fill=gray!30!white] (3.4,4.5)--(4,0)--(5,2)--(4.6,3.5)--cycle;
\draw[fill=gray!50!white] (3.4,1.3)--(5,2)--(4,0)--cycle;
\draw[gray] (0,2)--(5,2)--(4,0)--cycle;
\draw[gray] (0,2)--(4.667,1.333);
\draw[gray] (5,2)--(1.5,1.25);
\draw[gray] (4,0)--(2.75,2);
\draw[gray] (4,0)--(3.4,4.5)--(5,2);
\draw[gray] (0,2)--(3.4,4.5);
\draw[fill=black] (0,2) circle [radius=0.07];
\draw[fill=black] (5,2) circle [radius=0.07];
\draw[fill=black] (4,0) circle [radius=0.07];
\draw[fill=black] (4.6,3.5) circle [radius=0.07];
\draw (3.4,4.5)--(4.6,3.5)--(5,2);
\draw (4.6,3.5)--(4,0);
\draw[fill=black] (3.4,1.3) circle [radius=0.05];
\draw[fill=black] (3.4,4.5) circle [radius=0.07];
\draw[thick] (3.4,4.5)--(3.4,1.3);
\node at (0,1.7) {$o$};
\node at (4.4,0) {$p_1$};
\node at (5.3,2.3) {$p_2$};
\node at (3.9,4.5) {$p_3$};
\node[scale=0.9] at (3,1.2) {$\pi_3$};
\node at (4.1,2.8) {$\tet_0^\ast$};
\node at (2,2.8) {$\tet$};
\node at (2,0.7) {$\Delta$};
\end{scope}
\end{tikzpicture}    
\end{center} 
\caption{\small{Left: illustration of the sets $\lambda_0,\lambda_1,\lambda_2,\lambda_3$ of a generic (not sub-orthocentric) tetrahedron. Right: projecting $\tet_0^\ast$ on a triangular facet $\Delta$ of $\tet$.}}
\label{Mfigure216}
\end{figure}

In order to explain the shortcomings of the concept of blocking columns geometrically, consider the left picture in Figure~\ref{Mfigure215}. Here we see the triangular facet $\Delta$ of a tetrahedron $\tet$ with nonobtuse triangular facets only. If its fourth vertex $p_3$ does not project in the dark gray region, it does not cause a blocking column in the vertex Gramian $G_o$ of $\tet$. Thus, in order not to cause a blocking column in {\em any} of the vertex Gramians corresponding to $o,p_1,p_2$, it must project on the set of sub-orthocentric points of $\Delta$ displayed in the middle of Figure~\ref{Mfigure215}. 
\begin{figure}[h]
\begin{center}  
\begin{tikzpicture}[scale=0.8, every node/.style={scale=0.8}]
\draw[white!80!gray, fill=white!80!gray] (0,0)--(4,0)--(1,3)--cycle;
\draw[white!60!gray, fill=white!60!gray] (4,0)--(4,2)--(1,3)--cycle;
\draw[white!60!gray, fill=white!60!gray] (1,3)--(1,1)--(4,0)--cycle;
\draw (0,0)--(4,0)--(1,3)--cycle;
\draw (4,0)--(4,2)--(1,3)--cycle;
\draw[fill=black] (0,0) circle [radius=0.07];
\draw[fill=black] (4,0) circle [radius=0.07];
\draw[fill=black] (1,3) circle [radius=0.07];
\draw[dashed] (1,3)--(1,0);
\draw[dashed] (4,0)--(0.4,1.2);
\draw[fill=gray] (0.4,1.2) circle [radius=0.05];
\draw[fill=gray] (1,0) circle [radius=0.05];
\node at (4.1,-0.3) {$p_2$};
\node at (0.7,3) {$p_1$};
\node at (0,-0.3) {$o$};
\node at (1.5,1.5) {$\Delta$};
\node at (3,1.75) {$\Delta^\ast_0$};
\begin{scope}[shift={(5,0)}]
\draw[white!90!gray, fill=white!90!gray] (0,0)--(4,0)--(1,3)--cycle;
\draw (0,0)--(1,1);
\draw (1,1)--(1,3);
\draw (1,1)--(4,0);
\draw[fill=black] (0,0) circle [radius=0.07];
\draw[fill=black] (4,0) circle [radius=0.07];
\draw[fill=black] (1,3) circle [radius=0.07];
\node at (4,-0.3) {$p_2$};
\node at (0.7,3) {$p_1$};
\node at (0,-0.3) {$o$};
\node at (1.5,1.5) {$\Delta$};
\end{scope}
\begin{scope}[shift={(10,0)}]
\draw (0,0)--(4,0)--(1,3)--cycle;
\draw[very thick,gray] (0.4,1.2)--(0,0)--(1,0);
\draw[very thick,gray!30!white] (0.4,1.2)--(1,3)--(2,2);
\draw[very thick] (2,2)--(4,0)--(1,0);
\draw[gray,fill=gray] (0,0) circle [radius=0.07];
\draw[fill=black] (4,0) circle [radius=0.07];
\draw[gray!30!white,fill=gray!30!white] (1,3) circle [radius=0.07];
\draw[fill=gray] (0.4,1.2) circle [radius=0.05];
\draw[fill=gray] (1,0) circle [radius=0.05];
\draw[fill=gray] (2,2) circle [radius=0.05];
\node at (4,-0.3) {$p_2$};
\node at (0.7,3) {$p_1$};
\node at (0,-0.3) {$o$};
\node at (1,-0.3) {$\pi_1$};
\node at (2.2,2.2) {$\pi_0$};
\node at (0.1,1.2) {$\pi_2$};
\node[gray] at (0.5,0.3) {$\lambda_0$};
\node[gray!30!white] at (1.1,2.5) {$\lambda_1$};
\node at (3.3,0.3) {$\lambda_2$};
\end{scope}
\end{tikzpicture}    
\end{center} 
\caption{\small{The sub-orthocentric set of $\Delta$ is its subset of points that can be distinguished from points in its dual hull $\Delta^\ast$ by considering their projections on the edges of $\Delta$ only.}}
\label{Mfigure215}
\end{figure}
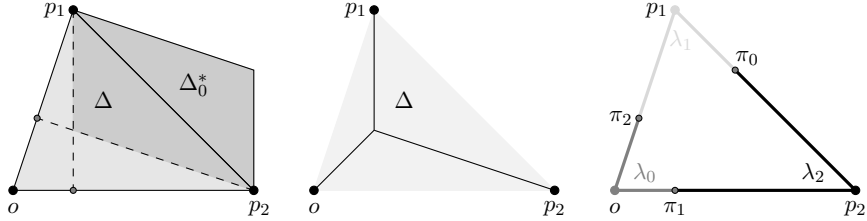\\[2mm]
The above can be rephrased in terms of the projections of the vertices of $\tet$ on its edges. Write $\pi_0,\pi_1$, and $\pi_2$ for the projections of $o$, $p_1$, and $p_2$ on their opposite edges, and define the closed sets $\lambda_0=\ol{\pi_1 o}\cup\ol{o\pi_2}$, $\lambda_1=\ol{\pi_2 p_1}\cup\ol{p_1\pi_0}$, and $\lambda_2=\ol{\pi_0 p_2}\cup\ol{p_2\pi_1}$. Note that $\lambda_0\cup\lambda_1\cup\lambda_2=\partial\Delta$, the boundary of $\Delta$, as depicted in the right picture in Figure~\ref{Mfigure215}, in which each set $\lambda_j$ has a different shade of gray. Then $p_3$ does not cause a blocking column in any of the vertex Gramians $G_0,G_1,G_2$ if and only if each of the sets $\lambda_0,\lambda_1,\lambda_2$ contains at least one of the projections of $p_3$ on its three opposite edges.\\[2mm]
Theorem~\ref{Mth-11} shows that a similar inspection of the locations of the projections of vertices on edges of a $4$-simplex $\four$ is too restrictive to find the required class $\FF$ of $4$-simplices. Indeed, if a tetrahedron $\tet$ is not sub-orthocentric, each pair of its vertices project on different points on their mutual opposite edge, as schematically depicted in the left of Figure~\ref{Mfigure216}.
\\[2mm]  
The requirement that the fifth vertex $p_4$ of a $4$-simplex with facet $\tet$ opposite $p_4$ causes no blocking column in the vertex Gramians corresponding to $o,p_1,p_2,p_3$ is equivalent to demanding that each of the four sets $\lambda_0,\lambda_1,\lambda_2,\lambda_3$, defined similarly as in the right picture of Figure~\ref{Mfigure215}, contains at least one of the six projections of $p_4$ on the respective edges of $\tet$.  See the left picture in Figure~\ref{Mfigure216}, in which each set $\lambda_j$ is in bold gray. Theorem~\ref{Mth-12} proves that not all five vertices of $\four$ can do so unless each tetrahedral facet is sub-orthocentric.
\subsection{Introducing a generalized concept of sub-orthocentricity}\label{Msect-6.4}
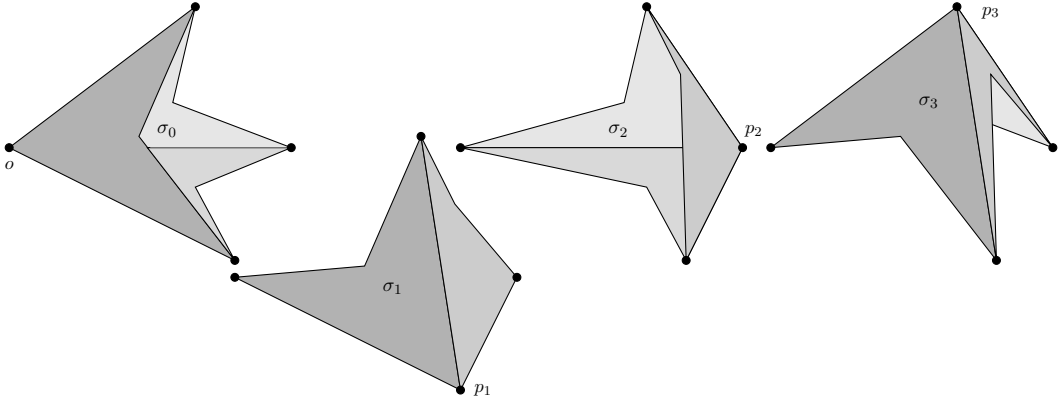
\begin{figure}[h]
\begin{center}  
\begin{tikzpicture}[scale=0.75, every node/.style={scale=0.75}]
\draw[fill=gray!20!white] (3.3,4.5)--(2.9,2.8)--(5,2)--(2,2)--cycle;
\draw[fill=gray!60!white] (0,2)--(4,0)--(2.3,2.2)--(3.3,4.5)--cycle;
\draw[fill=gray!30!white] (5,2)--(3.3,1.3)--(4,0)--(2.45,2)--cycle;
\draw[gray] (2.5,2)--(5,2);
\draw[fill=black] (0,2) circle [radius=0.07];
\draw[fill=black] (5,2) circle [radius=0.07];
\draw[fill=black] (4,0) circle [radius=0.07];
\draw[fill=black] (3.3,4.5) circle [radius=0.07];
\node[scale=0.85] at (0,1.7) {$o$};
\node[scale=0.9] at (2.8,2.3) {$\sigma_0$};
\begin{scope}[shift={(4,-2.3)}]
\draw[fill=gray!60!white] (0,2)--(2.3,2.2)--(3.3,4.5)--(4,0)--cycle;
\draw[fill=gray!40!white] (3.3,4.5)--(3.9,3.3)--(5,2)--(4,0)--cycle;
\draw[fill=black] (0,2) circle [radius=0.07];
\draw[fill=black] (5,2) circle [radius=0.07];
\draw[fill=black] (4,0) circle [radius=0.07];
\draw[fill=black] (3.3,4.5) circle [radius=0.07];
\node[scale=0.85] at (4.4,0) {$p_1$};
\node[scale=0.9] at (2.8,1.8) {$\sigma_1$};
\end{scope}
\begin{scope}[shift={(8,0)}]
\draw[fill=gray!20!white] (0,2)--(2.9,2.8)--(3.3,4.5)--(5,2)--cycle;
\draw[fill=gray!30!white] (4,0)--(3.3,1.3)--(0,2)--(5,2)--cycle;
\draw[fill=gray!40!white] (4,0)--(3.9,3.3)--(3.3,4.5)--(5,2)--cycle;
\draw[fill=black] (0,2) circle [radius=0.07];
\draw[fill=black] (5,2) circle [radius=0.07];
\draw[fill=black] (4,0) circle [radius=0.07];
\draw[fill=black] (3.3,4.5) circle [radius=0.07];
\node[scale=0.85] at (5.2,2.3) {$p_2$};
\node[scale=0.9] at (2.8,2.3) {$\sigma_2$};
\end{scope}
\begin{scope}[shift={(13.5,0)}]
\draw[fill=gray!20!white] (5,2)--(2.9,2.8)--(3.3,4.5)--(5,2)--cycle;
\draw[fill=gray!60!white] (0,2)--(2.3,2.2)--(4,0)--(3.3,4.5)--cycle;
\draw[fill=gray!40!white] (4,0)--(3.3,4.5)--(5,2)--(3.9,3.3)--cycle;
\draw[fill=black] (0,2) circle [radius=0.07];
\draw[fill=black] (5,2) circle [radius=0.07];
\draw[fill=black] (4,0) circle [radius=0.07];
\draw[fill=black] (3.3,4.5) circle [radius=0.07];
\node[scale=0.85] at (3.9,4.4) {$p_3$};
\node[scale=0.9] at (2.8,2.8) {$\sigma_3$};
\end{scope}
\end{tikzpicture}    
\end{center} 
\caption{\small{The closed subsets $\sigma_0,\sigma_1,\sigma_2,\sigma_3$ of $\partial\tet$ generalizing $\lambda_0,\lambda_1,\lambda_2$ in Figure~\ref{Mfigure214}.}}
\label{Mfigure217}
\end{figure}

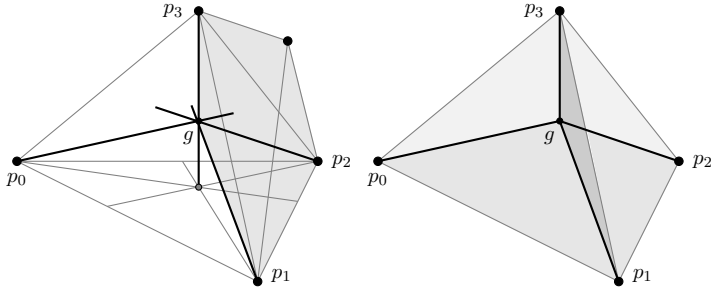
\begin{figure}[h]
\begin{center}  
\begin{tikzpicture}[scale=0.8, every node/.style={scale=0.8}]
\draw[white!80!gray, fill=white!80!gray] (4,0)--(5,2)--(4.5,4)--(3.02,4.5)--(3.02,2.67)--cycle;
\draw[gray] (0,2)--(5,2)--(4,0)--cycle;
\draw[gray] (0,2)--(4.667,1.333);
\draw[gray] (5,2)--(1.5,1.25);
\draw[gray] (4,0)--(3.02,4.5);
\draw[gray] (3.02,4.5)--(5,2);
\draw[gray] (0,2)--(3.02,4.5);
\draw[gray] (4,0)--(2.75,2);
\draw[gray] (5,2)--(4.5,4)--(3.02,4.5);
\draw[gray] (4.5,4)--(4,0);
\draw[fill=black] (0,2) circle [radius=0.07];
\draw[fill=black] (5,2) circle [radius=0.07];
\draw[fill=black] (4.5,4) circle [radius=0.07];
\draw[fill=black] (4,0) circle [radius=0.07];
\draw[thick] (5,2)--(2.3,2.9);
\draw[thick] (4,0)--(2.9,2.93);
\draw[thick] (0,2)--(3.6,2.8);
\draw[fill=black] (3.02,2.67) circle [radius=0.05];
\draw[fill=black] (3.02,4.5) circle [radius=0.07];
\draw[thick] (3.02,4.5)--(3.02,1.57);
\draw[fill=gray] (3.02,1.57) circle [radius=0.05];
\node[scale=0.9] at (0,1.7) {$p_0$};
\node[scale=0.9] at (4.4,0.1) {$p_1$};
\node[scale=0.9] at (5.4,2) {$p_2$};
\node[scale=0.9] at (2.6,4.5) {$p_3$};
\node[scale=0.9] at (2.85,2.4) {$g$};
\begin{scope}[shift={(6,0)}]
\draw[white!90!gray, fill=white!90!gray] (0,2)--(3.02,4.5)--(5,2)--(3.02,2.67)--cycle;
\draw[white!60!gray, fill=white!60!gray] (4,0)--(3.02,2.67)--(3.02,4.5)--cycle;
\draw[white!80!gray, fill=white!80!gray] (4,0)--(3.49,2.5)--(5,2)--cycle;
\draw[white!80!gray, fill=white!80!gray] (0,2)--(3.02,2.67)--(4,0)--cycle;
\draw[gray] (0,2)--(4,0)--(5,2);
\draw[gray] (4,0)--(3.02,4.5);
\draw[gray] (3.02,4.5)--(5,2);
\draw[gray] (0,2)--(3.02,4.5);
\draw[fill=black] (0,2) circle [radius=0.07];
\draw[fill=black] (5,2) circle [radius=0.07];
\draw[fill=black] (4,0) circle [radius=0.07];
\draw[thick] (5,2)--(3.02,2.67);
\draw[thick] (4,0)--(3.02,2.67);
\draw[thick] (0,2)--(3.02,2.67);
\draw[fill=black] (3.02,2.67) circle [radius=0.05];
\draw[fill=black] (3.02,4.5) circle [radius=0.07];
\draw[thick] (3.02,4.5)--(3.02,2.67);
\node[scale=0.9] at (0,1.7) {$p_0$};
\node[scale=0.9] at (4.4,0.1) {$p_1$};
\node[scale=0.9] at (5.4,2) {$p_2$};
\node[scale=0.9] at (2.6,4.5) {$p_3$};
\node[scale=0.9] at (2.85,2.4) {$g$};
\end{scope}
\end{tikzpicture}    
\end{center} 
\caption{\small{Left: the sets $\tet_0$ and $\tet_0^\ast$ of an orthocentric tetrahedron $\tet$. Right: the sub-orthocentric set $\tet_\ast$ of the orthocentric tetrahedron $\tet$.}}
\label{Mfigure219}
\end{figure}

Instead of considering the projections of the vertices of a $4$-simplex $\four$ on each of its six opposite {\em edges}, we will now consider its projections on its four opposite {\em triangular} facets. Given a nonobtuse tetrahedral facet $\tet$ of $\four$ as in Figure~\ref{Mfigure216} on the right, assume that interior of the part $\tet_0^\ast$ of $\tet^\ast$ opposite $o$ is nonempty. It is depicted in light gray. Let $\Delta_0$ be the triangular facet shared by $\tet$ and $\tet_0^\ast$. The projection of $\tet_0^\ast\setminus \Delta_0$ on the triangular facet $\Delta_3$ of $\tet$ opposite $p_3$ equals the darker gray interior of the triangle with vertices $p_1,p_2$, and $\pi_3$, where $\pi_3$ is the projection of $p_3$ on $\Delta_3$, together with the relative interior of the edge $\ol{p_1p_2}$.\\[2mm]
We want to guarantee that the vertex $p_4$ of $\four$ does not project in $\tet_0^\ast\setminus\Delta_0$, solely by inspecting the projections of $p_4$ onto the four facets of $\tet$. For this it is {\em sufficient} that one of these projections lies {\em outside} the projection of $\tet_0^\ast\setminus\Delta_0$ on those facets. This is equivalent with $p_4$ projecting in the (closed) part of the boundary $\partial\tet$ shaded (light, middle, or darker) gray in the left picture in Figure~\ref{Mfigure217}. We denote this part, which forms the generalization of the set $\lambda_0$ from Figures~\ref{Mfigure215} and~\ref{Mfigure216}, by $\sigma_0$. To guarantee that $p_4$ does not project in $\tet^\ast_j\setminus\Delta_j$ for any $j\in\{0,1,2,3\}$, each of the subsets $\sigma_0,\sigma_1,\sigma_2,\sigma_3$ of $\partial\tet$ depicted in Figure~\ref{Mfigure217} should contain at least one projection of $p_4$ on a triangular facet of $\tet$. Or equivalently, no set $\partial\tet\setminus \sigma_j$ should contain {\em all four} of the projections of $p_4$ on the triangular facets of $\tet$. 
\\[2mm] 
Consider in some more detail the set $\sigma_0\subset\partial\tet$, now in Figure~\ref{Mfigure218}. As in Figure~\ref{Mfigure216}, write $\pi_j$ for the projection of $p_j$ on its opposite facet $\Delta_j$. Let $T_j\subset\tet$ be the tetrahedron with vertices $p_1,p_2,p_3$, and $\pi_j$. Now, observe that if the projection of $p_4$ on $\tet$ lies in the interior $\inter(\tet_0)$ of the intersection $\tet_0=T_1\cap T_2\cap T_3$, or in the interior $\inter(\Delta_0)$ of $\Delta_0$, the set $\sigma_0$ will {\em not} contain a projection of $p_4$ on a triangular facet of $\tet$.
\begin{figure}[h]
\begin{center}  
\begin{tikzpicture}[scale=0.8, every node/.style={scale=0.8}]
\draw[fill=gray!20!white] (3.3,4.5)--(2.9,2.8)--(5,2)--(2,2)--cycle;
\draw[fill=gray!60!white] (0,2)--(4,0)--(2.3,2.2)--(3.3,4.5)--cycle;
\draw[fill=gray!30!white] (5,2)--(3.3,1.3)--(4,0)--(2.45,2)--cycle;
\draw[gray] (0,2)--(4,0)--(5,2);
\draw[gray] (4,0)--(3.3,4.5)--(5,2);
\draw[gray] (0,2)--(3.3,4.5);
\draw[gray] (2.5,2)--(5,2);
\draw[gray] (3.3,1.3)--(5,2)--(4,0)--cycle;
\draw[gray] (3.3,4.5)--(2.3,2.2)--(4,0);
\draw[gray] (3.3,4.5)--(2.9,2.8)--(5,2);
\draw[gray] (5,2)--(2.3,2.2);
\draw[gray] (4,0)--(2.9,2.8);
\draw[gray] (3.3,4.5)--(3.3,1.3);
\draw[fill=black] (2.3,2.2) circle [radius=0.05];
\draw[fill=black] (2.9,2.8) circle [radius=0.05];
\draw[fill=black] (3.3,1.3) circle [radius=0.05];
\draw[fill=black] (0,2) circle [radius=0.07];
\draw[fill=black] (5,2) circle [radius=0.07];
\draw[fill=black] (4,0) circle [radius=0.07];
\draw[fill=black] (3.3,4.5) circle [radius=0.07];
\node[scale=0.85] at (0,1.7) {$o$};
\node[scale=0.9] at (2,1.8) {$\sigma_0$};
\node at (0,1.7) {$o$}; 
\node at (4.4,0) {$p_1$};
\node at (5.3,2.3) {$p_2$};
\node at (3.9,4.5) {$p_3$};
\node[scale=0.9] at (3.1,1.5) {$\pi_3$};
\node[scale=0.9] at (2.1,2.4) {$\pi_2$};
\node[scale=0.9] at (2.8,2.5) {$\pi_1$};
\begin{scope}[shift={(10,0)}]
\draw[gray,fill=gray!40!white] (4,0)--(5,2)--(3.3,4.5)--cycle;
\draw[gray,fill=gray!20!white] (4,0)--(3.3,4.5)--(3.3,1.3)--cycle;
\draw[gray] (4,0)--(5,2)--(3.3,4.5)--cycle;
\draw[gray] (3.3,1.3)--(4,0);
\draw[gray] (3.3,1.3)--(5,2);
\draw[gray] (3.3,1.3)--(3.3,4.5);
\draw[fill=black] (3.3,4.5) circle [radius=0.07];
\draw[fill=black] (5,2) circle [radius=0.07];
\draw[fill=black] (4,0) circle [radius=0.07];
\draw[fill=black] (3.3,1.3) circle [radius=0.05];
\node at (4.1,2.5) {$T_3$};
\node at (4.4,0) {$p_1$};
\node at (5.3,2.3) {$p_2$};
\node at (3.9,4.5) {$p_3$};
\node[scale=0.9] at (3.1,1.5) {$\pi_3$};
\end{scope}
\begin{scope}[shift={(7,0)}]
\draw[gray,fill=gray!40!white] (4,0)--(5,2)--(3.3,4.5)--cycle;
\draw[gray,fill=gray!20!white] (4,0)--(3.3,4.5)--(2.1,2.4)--cycle;
\draw[gray] (4,0)--(5,2)--(3.3,4.5)--cycle;
\draw[gray] (2.1,2.4)--(4,0);
\draw[gray] (2.1,2.4)--(5,2);
\draw[gray] (2.1,2.4)--(3.3,4.5);
\draw[fill=black] (3.3,4.5) circle [radius=0.07];
\draw[fill=black] (5,2) circle [radius=0.07];
\draw[fill=black] (4,0) circle [radius=0.07];
\draw[fill=black] (2.1,2.4) circle [radius=0.05];
\node at (4.3,1.5) {$T_2$};
\node at (4.4,0) {$p_1$};
\node at (5.3,2.3) {$p_2$};
\node at (3.9,4.5) {$p_3$};
\node[scale=0.9] at (2,2.7) {$\pi_2$};
\end{scope}
\begin{scope}[shift={(3,0)}]
\draw[gray,fill=gray!40!white] (4,0)--(5,2)--(3.3,4.5)--cycle;
\draw[gray,fill=gray!20!white] (4,0)--(3.3,4.5)--(2.9,2.8)--cycle;
\draw[gray] (4,0)--(5,2)--(3.3,4.5)--cycle;
\draw[gray] (2.9,2.8)--(4,0);
\draw[gray] (2.9,2.8)--(5,2);
\draw[gray] (2.9,2.8)--(3.3,4.5);
\draw[fill=black] (3.3,4.5) circle [radius=0.07];
\draw[fill=black] (5,2) circle [radius=0.07];
\draw[fill=black] (4,0) circle [radius=0.07];
\draw[fill=black] (2.9,2.8) circle [radius=0.05];
\node at (4.3,1.5) {$T_1$};
\node at (4.4,0) {$p_1$};
\node at (5.3,2.3) {$p_2$};
\node at (3.9,4.5) {$p_3$};
\node[scale=0.9] at (2.8,2.5) {$\pi_1$};
\end{scope}
\end{tikzpicture}    
\end{center} 
\caption{\small{The tetrahedra $T_1,T_2,T_3$ that intersect as the tetrahedron $\tet_0\subset\tet$.}}
\label{Mfigure218}
\end{figure}
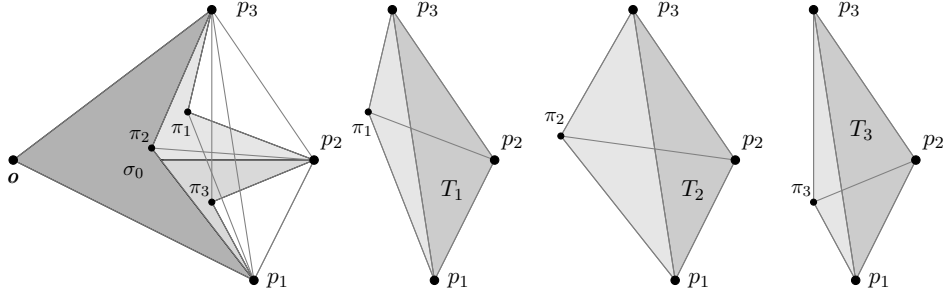\\[2mm] 
Thus, we have identified a set points in $\tet$ that can be distinguished from those in $\tet_0^\ast\setminus\Delta_0$ by solely inspecting the four projections of $p_4$ on the four facets of $\tet$.
\begin{rem}{\rm The intersection $\tet_0$ of the tetrahedra $T_1,T_2,T_3$ is itself a tetrahedron, because $T_1,T_2,T_3$ share the same triangular base facet $F$. This will be proved in Lemma~\ref{Mlem-10}.}
\end{rem}
In a similar manner as the tetrahedron $\tet_0\subset\tet$, we can define tetrahedral subsets $\tet_1,\tet_2,\tet_3$ of $\tet$, with $\tet_j$ sharing the facet $\Delta_j$ with $\tet$. Note that $\tet_0,\tet_1,\tet_2,\tet_3$ have disjoint interior, as each three of them project as disjoint subsets on a facet of $\tet$.\\[2mm]
Now, let $\tet$ be a nonobtuse tetrahedron. Consider the four subsets $\tet_0^\ast,\tet_1^\ast,\tet_2^\ast,\tet_3^\ast$ of $\tet^\ast$ as defined in (\ref{MDandE}). Observe that $\tet_j^\ast$ has empty interior if and only if $\Delta_j$ makes a right angle with one of the other facets of $\tet$, which is the case if and only if $\tet_j$ has empty interior. Define the subset $\tet_\ast$ of $\tet$ by
\be\label{Meq-100} \tet_\ast = \tet\setminus \bigcup (\inter(\tet_j)\cup\inter(\Delta_j)),\ee
where the (first) union is over all $j\in\{0,1,2,3\}$ for which $\tet_j^\ast$ has nonempty interior. Thus, for each $\tet_j^\ast$ with nonempty interior, we remove from $\tet$ both $\inter(\tet_j)$ and $\inter(\Delta_j)$. What remains is a closed subset of $\tet$ that generalizes the concept of sub-orthocentric points of a triangle.   
\begin{Def}[sub-orthocentric points]{\rm Let $\tet$ be a nonobtuse tetrahedron. We call the elements of the closed nonempty set $\tet_\ast\subset\tet$ the {\em sub-orthocentric points}
\index{sub-orthocentric points} of $\tet$}.
\end{Def} 
As a simple but extremal example, consider the {\em orthocentric} tetrahedron $\tet$ displayed in the left of Figure~\ref{Mfigure218}, together with the sets $\tet_0^\ast$ and $\tet_0$. Each of the tetrahedra $\tet_j$ is simply the convex hull of the orthocenter and the triangular facet $\Delta_j$ of $\tet$. Removing the interiors of these tetrahedra $\tet_j$ and of their facets $\Delta_j$ results in the closed, two-dimensional sub-orthocentric set $\tet_\ast$ displayed in the right picture of Figure~\ref{Mfigure219}.
\\[2mm]
For a {\em sub-orthocentric} tetrahedron $\tet$, the set $\tet_\ast$ is rather more complicated to visualize. It helps to move the vertex $p_3$ in the left in Figure~\ref{Mfigure218} over a small distance towards $o$, while $o, p_1$, and $p_2$ remain where they are. What results is a sub-orthocentric tetrahedron, and moreover, the altitudes from both $p_1$ and $p_2$ remain the same, whereas the ones from $o$ and $p_3$ have slightly shifted, as visible in Figure~\ref{Mfigure220}. It is easy to verify that the tetrahedra $T_1$ and $T_2$ are now both contained in $T_3$, thus, $\tet_0=T_1\cap T_2$, which is the tetrahedron with facet $\Delta_0$ and top $s_2$, the intersection of the altitudes from $p_1$ and $p_2$. Similarly, $\tet_3$ is the tetrehadron with facet $\Delta_3$ and apex $s_2$, hence, $\tet_0$ and $\tet_3$ share the triangular facet with vertices $s_2,p_1$, and $p_2$. In the left of Figure~\ref{Mfigure220}, $\tet_3$ is depicted in dark gray, and on top of it in lighter gray we see~$\tet_0$.
\begin{figure}[h]
\begin{center}  
\begin{tikzpicture}[scale=0.8, every node/.style={scale=0.8}]
\draw[white!80!gray, fill=white!80!gray] (4,0)--(5,2)--(2.6,4.25)--(2.98,2.78)--cycle;
\draw[white!60!gray, fill=white!60!gray] (4,0)--(5,2)--(3.0,2.78)--(0,2)--cycle;
\draw[gray] (0,2)--(5,2)--(4,0)--cycle;
\draw[gray] (0,2)--(4.667,1.333);
\draw[gray] (5,2)--(1.5,1.25);
\draw[gray] (4,0)--(2.6,4.25);
\draw[gray] (2.6,4.25)--(5,2);
\draw[gray] (0,2)--(2.6,4.25);
\draw[gray] (4,0)--(2.75,2);
\draw[fill=black] (0,2) circle [radius=0.07];
\draw[fill=black] (5,2) circle [radius=0.07];
\draw[fill=black] (4,0) circle [radius=0.07];
\draw[thick] (5,2)--(2.3,3.03);
\draw[thick] (4,0)--(2.9,2.93);
\draw[thick] (0,2)--(3.05,3.51);
\draw (2.95,2.78)--(2.6,4.25);
\draw (2.95,2.78)--(0,2);
\draw[fill=black] (2.95,2.78) circle [radius=0.05];
\draw[fill=black] (2.6,3.3) circle [radius=0.05];
\draw[fill=black] (2.6,4.25) circle [radius=0.07];
\draw[thick] (2.6,4.25)--(2.6,1.62);
\draw[fill=gray] (3.0,1.57) circle [radius=0.05];
\node[scale=0.9] at (0,1.7) {$p_0$};
\node[scale=0.9] at (4.4,0.1) {$p_1$};
\node[scale=0.9] at (5.4,2) {$p_2$};
\node[scale=0.9] at (2.6,4.5) {$p_3$};
\node[scale=0.9] at (2.4,3.5) {$s_1$};
\node[scale=0.9] at (3.3,2.8) {$s_2$};
\begin{scope}[shift={(6.5,0)}]
\draw[white!40!gray, fill=white!40!gray] (2.6,3.3)--(3.01,2.78)--(5,2)--cycle;
\draw[white!80!gray, fill=white!80!gray] (2.6,3.3)--(3.01,2.78)--(0,2)--cycle;
\draw[gray] (0,2)--(5,2)--(4,0)--cycle;
\draw[gray] (0,2)--(4.667,1.333);
\draw[gray] (5,2)--(1.5,1.25);
\draw[gray] (4,0)--(2.6,4.25);
\draw[gray] (2.6,4.25)--(5,2);
\draw[gray] (0,2)--(2.6,4.25);
\draw[gray] (4,0)--(2.75,2);
\draw (0,2)--(3.01,2.78);
\draw (2.6,3.3)--(5,2);
\draw[white!80!gray, fill=white!80!gray] (2.6,3.3)--(3.01,2.78)--(2.6,4.25)--cycle;
\draw (3.01,2.78)--(2.6,4.25);
\draw[fill=black] (0,2) circle [radius=0.07];
\draw[fill=black] (5,2) circle [radius=0.07];
\draw[fill=black] (4,0) circle [radius=0.07];
\draw[white!40!gray, fill=white!40!gray] (2.6,3.3)--(3.01,2.78)--(4,0)--cycle;
\draw (2.6,3.3)--(4,0);
\draw (3.01,2.78)--(2.6,3.3);
\draw[thick] (5,2)--(3.01,2.78);
\draw[thick] (4,0)--(3.01,2.78);
\draw[thick] (0,2)--(2.6,3.3);
\draw[fill=black] (3.01,2.78) circle [radius=0.05];
\draw[fill=black] (2.6,3.3) circle [radius=0.05];
\draw[fill=black] (2.6,4.25) circle [radius=0.07];
\draw[thick] (2.6,4.25)--(2.6,3.3);
\draw[fill=gray] (3.02,1.57) circle [radius=0.05];
\node[scale=0.9] at (0,1.7) {$p_0$};
\node[scale=0.9] at (4.4,0.1) {$p_1$};
\node[scale=0.9] at (5.4,2) {$p_2$};
\node[scale=0.9] at (2.6,4.5) {$p_3$};
\node[scale=0.9] at (2.25,3.45) {$s_1$};
\node[scale=0.9] at (3.4,2.35) {$s_2$};
\end{scope}
\end{tikzpicture}    
\end{center} 
\caption{\small{Left: the sets $\tet_0$ and $\tet_3$ of a sub-orthocentric tetrahedron $\tet$. Right: The triangles $\nabla_0,\dots,\nabla_3$ such that $\tet_\ast = \cup\,\conv(\nabla_i,\nabla_j)$.}}
\label{Mfigure220}
\end{figure}
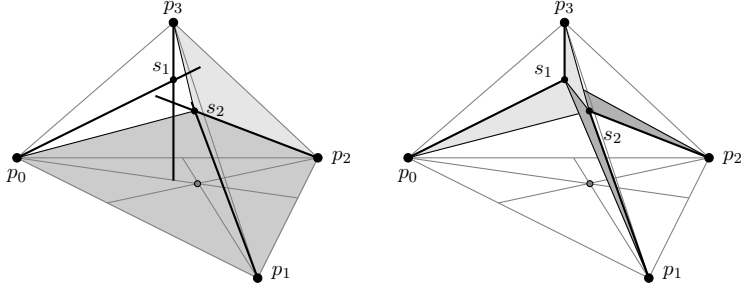 \\[2mm]
Likewise, the intersection $s_1$ of the altitudes from $o$ and $p_3$ is the apex of the tetrahedra $\tet_1$ and $\tet_2$ with respective base facets $\Delta_1$ and $\Delta_2$. Interestingly, we see that the union of the closed tetrahedra $\tet_0,\tet_1,\tet_2,\tet_3$ does not equal $\tet$. To study the sub-orthocentric set $\tet_\ast$, let $\nabla_j$ be the triangle with vertices $s_1, s_2$, and $p_j$ and note that $\nabla_0$ and $\nabla_3$ lie in the same plane, as do $\nabla_1$ and $\nabla_2$.  Closer inspection of the right picture in Figure~\ref{Mfigure220} reveals that
\be\label{Meq-70} \tet_\ast = \bigcup_{i,j\in I_0^3} \conv(\nabla_i,\nabla_j). \ee
Note that the sets $\conv(\nabla_0,\nabla_3)$ and $\conv(\nabla_1,\nabla_2)$ are triangles in $\tet_\ast$ containing $\tet_0\cap\tet_3$ and $\tet_1\cap\tet_2$, respectively, whereas the remaining four sets are closed nondegenerate tetrahedra. If $\tet$ is orthocentric, then $s_1=s_2=g$. Hence $\nabla_j$ is the line segment between $p_j$ and $g$, after which (\ref{Meq-70}) reduces to precisely the set depicted in the right of Figure~\ref{Mfigure218}.
\begin{rem}{\rm In the extreme case that $p_3$ is located orthogonally above the vertex $p_0$ of $\Delta$, then $s_1=p_0$ and $s_2$ is the orthocenter of $\Delta_3$, and (\ref{Meq-70}) shows (correctly) that $\tet_\ast$ is $\tet$ with only $\inter(\tet_0)\cup\inter(\Delta_0)$ removed. Interestingly, a {\em path-tetrahedron} $\tet$ is the only tetrahedron for which
\be \tet^\ast = \tet = \tet_\ast, \ee
because it has the property that each triangular facet makes a right angle with another facet. } 
\end{rem}
Armed with the new concept of the sub-orthocentric set of a nonobtuse tetrahedron, we can now also define a sub-orthocentric $4$-simplex in analogy to the three-dimensional case.
\begin{Def}[sub-orthocentric $4$-simplex] \index{sub-orthocentric $4$-simplex} {\rm A $4$-simplex $\four$ with nonobtuse tetrahedral facets is called {\em sub-orthocentric} if each vertex of $\four$ projects on a sub-orthocentric point of its opposite tetrahedral facet.}
\end{Def}
As the set $\tet_\ast\subset\tet$, a sub-orthocentric $4$-simplex is trivially nonobtuse due to Theorem~\ref{Mth-1}. Observe also that all {\em edges} of $\tet$ are in $\tet_\ast$, because only points from $\inter(\tet)$ and the interior of its triangular facets are removed to get $\tet_\ast$.\\[2mm]
The following result is now of great interest.
\begin{Th}\label{Mth-14} A $4$-simplex whose facets are all sub-or\-tho\-cen\-ter, is sub-or\-tho\-cen\-ter.
\end{Th}
{\bf Proof. } Let $v$ be a vertex a $4$-simplex $\four$ with sub-orthocentric tetrahedral facets, and write $\tet$ for the tetrahedral facet of $\four$ opposite $v$. Then $\tet$ is nonobtuse. Write $\pi$ for the projection of $v$ on $\tet$. By assumtion, $v$ projects in the sub-orthocentric set of each of the four triangular facets of $\tet$, hence so does $\pi$. Taking into account the proof of Lemma~\ref{Mlem-11} and using the right picture in Figure~\ref{Mfigure220}, this implies that $\pi$ lies on the boldface part of an altitude of $\tet$, or on $\ol{s_1s_2}$. But this implies that $\pi\in\tet_\ast$. Hence, $\four$ is sub-orthocentric.~\hfill $\Box$\\[2mm]
Thus, the property of sub-orthocentricity, contrary to that of nonobtusity, is transferable from the facets of a $4$-simplex to the $4$-simplex itself. This motivates the final section of this chapter.
\section{sub-orthocentric simplices}\label{Msect-7}
Let $S$ be a nonobtuse $n$-simplex with vertices $v_0,\dots,v_n$ and opposite facets $F_0,\dots,F_n$. Let $j\in I_0^n$ be given. According to Theorem~\ref{Mth-1}, $v_j$ projects onto a point $\pi_j$ in $F_j$. This induces a face-to-face triangulation $\TT_j$ of $S$ into $n$ subsimplices
\be\label{Meq-37} S_{j\ell} = \conv\{\pi_j,F_\ell\}, \hdrie \ell \in I_0^n\setminus\{j\}. \ee
Note that in case $\pi_j\in\partial F_\ell$ some, but at most $\nmo$ of the simplices $S_{j\ell}$ may be degenerate. Now, for each given $\ell\in I_0^n$, consider the intersection
\be\label{Meq-38} S_\ell = \bigcap_{j\not=\ell} S_{j\ell} = \bigcap_{j\not=\ell}\conv\{\pi_j,F_\ell\} \supset F_\ell. \ee
For $n=3$, the intersections $S_0,S_1,S_2,S_3$ are precisely the tetrahedra defined in Section~\ref{Msect-6.4}. As announced there, we will prove that each intersection $S_\ell$ in (\ref{Meq-38}) is an $n$-simplex.
\begin{Le}\label{Mlem-10} Let $T_1,\dots,T_p$ be $n$-simplices with common facet $F_0$, and assume that they all lie at the same side of $F_0$. Then also $T=\cap_{j=1}^p T_j$ is an $n$-simplex with facet $F_0$.
\end{Le}
{\bf Proof. } It is sufficient to prove the statement for $p=2$. As in Section~\ref{Msect-3.1}, let $T_1$ be an $n$-simplex with facets $F_0,\dots, F_n$. For each $j\in I_0^n$, let $H_j$ be the affine hyperplane with $F_j\subset H_j$, and $H_j^+$ the closed half-space separated by $H_j$ with $T_1\subset H_j^+$. Given an $n$-simplex $T_2\subset H_0^+$ having $F_0$ as a facet,
\be T_2\cap T_1 = T_2 \cap \bigcap_{j=0}^n H_j^+ = T_2 \cap \bigcap_{j=1}^n H_j^+.\ee
It suffices to prove that $T_2\cap H_1^+$ is an $n$-simplex with facet $F_0$. For this, note that $H_1$ contains the facet $G=F_1\cap F_0$ of $F_0$. Thus, $H_1$ bisects $T_2$ in its edge opposite $G$ into two $n$-simplices, with the part in $H_1^+$ being an $n$-simplex with facet $F_0$.~\hfill $\Box$\\[2mm]
Since the sets $S_{j\ell}$ defined in (\ref{Meq-37}) have $F_\ell$ as common facet, Lemma~\ref{Mlem-1} shows that $S_0,\dots,S_n$ in (\ref{Meq-38}) are simplices. Moreover, the same argument as in Section~\ref{Msect-6.4} can be used to prove that they have disjoint interior. The next result shows when they have empty interior. 
\begin{Pro}\label{Mprop-13} Let $S_0^\ast,\dots,S_n^\ast$ be the parts of the dual hull $S^\ast$ of a nonobtuse $n$-simplex $S$ as in (\ref{MDandE}), and let $S_0,\dots,S_n$ the parts of $S$ as in (\ref{Meq-38}). Then 
\be \inter\left(S_\ell^\ast\right) =\emptyset \hdrie \Leftrightarrow \hdrie \inter\left(S_\ell\right)=\emptyset. \ee
\end{Pro} 
 {\bf Proof. } The interior of the set $S_\ell^\ast$ is empty if and only if the facet $F_\ell$ of $S$ makes a right angle with one of the remaining facets, which holds true if and only if the projection $\pi_j$ of a vertex $v_j$ of $F_\ell$ onto its opposite facet $F_j$ lands in $F_\ell\cap F_j$, which happens if and only if $S_{j\ell}$ in (\ref{Meq-37}), and hence $S_\ell$ in (\ref{Meq-38}), degenerates. This proves the statement.~\hfill $\Box$\\[2mm]
 We are now ready to define the sub-orthocentric points of an arbitrary $n$-simplex $S$. In words, for each set $S_\ell$ with non-empty interior, we remove from $S$ both $\inter(S_\ell)$ and $\inter(F_\ell)$.  
 \begin{Def}[sub-orthocentric set] \index{sub-orthocentric set} {\rm Let $S$ be a nonobtuse $n$-simplex. The {\em sub-orthocentric} set $S_\ast\subset S$ of $S$ is the set
 \be S_\ast = S \setminus \bigcup (\inter(S_\ell)\cup \inter(F_\ell)), \ee
 where the union ranges only over the sets $S_\ell$ from (\ref{Meq-38}) with {\em nonempty} interior}.
 \end{Def}
To motivate this definition, recall that it was our goal to remove from $S$ precisely those points that cannot be distinguished from points in $\inter(S_\ell^\ast)$ solely by looking at their respective projections on the facets $F_j,j\not=\ell$, of $S$, as discussed in detail for $n\in\{2,3\}$ in Section~\ref{Msect-6}. But if $\inter(S_\ell^\ast)=\emptyset$, then also $\inter(S_\ell)=\emptyset$ due to Proposition~\ref{Mprop-13}, and nothing needs to be removed from $S$, in particular not the nonempty set $\inter(F_\ell)$. In case $\inter(S_j)$ is not empty, also $\inter(F_\ell)$ needs to be removed. Note that $S_\ast$ contains all $(n-2)$-facets of $S$.
\begin{Def}[sub-orthocentric simplex] {\rm We call an $n$-simplex $S$ with nonobtuse facets {\em sub-orthocentric} \index{sub-orthocentric simplex} if each vertex of $S$ projects in the sub-orthocentric set of its opposite facet.}
\end{Def}
Observe again that trivially, a sub-orthocentric simplex is nonobtuse.\\[2mm]
The following conjecture we proved for $n=4$ in Theorem~\ref{Mth-14}. If it holds also for all larger dimensions, the sub-orthocentric $n$-simplices are the desired generalizations of equilateral triangles and sub-orthocentric tetrahedra.
\begin{Con} \label{Mconj-2} A simplex with only sub-orthocentric facets is sub-orthocentric.
\end{Con} 
{\bf Outline: } Assume that $S$ is an $n$-simplex with sub-orthocentric $(\nmo)$-facets. This implies in particular that $S$ has {\em nonobtuse} facets. It moreover implies that each vertex $v$ of $S$ projects in the sub-orthocentric set of the facets of its opposite facet, hence in all $(n-2)$-facets of $S$ opposite $v$. From this, we need to show that $v$ projects in the sub-orthocentric set of its opposite facet. It is sufficient to show that the intersection of the cylinders of the sub-orthocentric sets of the facets $F$ of $S$ is contained in the suborthocentric set of $S$ itself.\hfill $\Box$ \\[2mm]
The immediate consequence of this conjecture is the following final result. 
\begin{Th} Let $S$ be an $n$-simplex whose $k$-facets are all sub-orthocentric. Then, assuming the validity of Conjecture~\ref{Mconj-2} we have that $S$ is nonobtuse.
\end{Th} 
{\bf Proof. } A finite repetitive application of Conjecture~\ref{Mconj-2} shows that all $\ell$-facets with $k\leq\ell\leq n$ are sub-orthocentric. All sub-orthocentric simplices are nonobtuse.~\hfill $\Box$\\[2mm]
Note that Conjecture~\ref{Mconj-2} is the more interesting, if the set of $n$-simplices with sub-orthocentric facets is even {\em strictly} contained in the set of sub-orthocentric $n$-simplices. Since the sub-orthocentric set $\tet_\ast$ of a nonobtuse tetrahedron $\tet$ is in general much larger than the subset of points in $\tet$ that project on the sub-orthocentric sets of each of the triangular facets of $\tet$, as discussed in Section~\ref{Msect-6.4}, this is most probably the~case.   

\subsection*{Acknowledgments} 
Jan Brandts and Apo Cihangir acknowledge the support by Research Project 613.001.019 of the Netherlands Organisation for Scientific Research (NWO). They also thank Michal K\v{r}\'{\i}\v{z}ek for his valuable comments on earlier versions of this paper.


\begin{thebibliography}{10}  

\bibitem{BePl}
A.~Berman and B.~Plemmons (1994).
\newblock Nonnegative matrices in the Mathematical Sciences.
\newblock {\em Society for Industrial and Applied Mathematics}.

\bibitem{BeSh}
A.~Berman, N.~Shaked-Monderer (2003). 
\newblock Completely positive matrices. 
\newblock {\em World Scientific Publishing} Co., Inc., River Edge, NJ.

\bibitem{BoVa}  
S.P.~Boyd and L.~Vandenberghe (2004). 
\newblock Convex Optimization. 
\newblock {\em Cambridge University Press}.

\bibitem{Bra}
D. Braess (2001).
\newblock Finite elements: theory, fast solvers, and applications in solid mechanics.
\newblock {\em University Press, Cambridge. Second Edition}. 

\bibitem{BrKr1}
J.H. Brandts and M. K\v{r}\'{\i}\v{z}ek (2003).
\newblock Gradient superconvergence on uniform simplicial partitions of polytopes.
\newblock {\em IMA Journal of Numerical Analysis}, 23:1--17. 

\bibitem{BrKoKr1}
J.H. Brandts, S. Korotov, and M. K\v{r}\'{\i}\v{z}ek(2008).
\newblock The discrete maximum principle for linear simplicial finite element approximations of a reaction-diffusion problem.
\newblock {\em Linear Algebra and its Applications}, 429(10):2344--2357.

\bibitem{BrKoKr}
J.H. Brandts, S. Korotov, and M. K\v{r}\'{\i}\v{z}ek (2007).
\newblock Dissection of the path-simplex in~$\mathbb{R}^n$ into~$n$ path-subsimplices.
\newblock {\em Linear Algebra and its Applications}, 421(2-3):382--393.

\bibitem{BrKoKrSo}
J.H. Brandts, S. Korotov, M. K\v{r}\'{\i}\v{z}ek, and J. \v{S}olc (2009).
\newblock On acute and nonobtuse simplicial partitions. 
\newblock {\em Society for Industrial and Applied Mathematics}, 51(2):317--335.

\bibitem{Cou1}
N.A.~Court (1934).
\newblock Notes on the orthocentric tetrahedron.
\newblock {\em American Mathematical Monthly,} 41:499--502.

\bibitem{Cou2}
N.A.~Court (1953).
\newblock The semi-orthocentric tetrahedron.
\newblock {\em American Mathematical Monthly,} 60:306--310.

\bibitem{Cox}
H.S.M.~Coxeter (1973).
\newblock Regular polytopes, third edition.
\newblock {\em Dover books on mathematics,} Dover Publications Inc. New York.

\bibitem{Cox2}
H.S.M.~Coxeter, (1991).
\newblock Orthogonal trees.
\newblock {\em Proc. 7th ACM Symp. Computational Geometry,} 89--97.

\bibitem{DeMaMa}
C.~Dellacherie, S.~Martinez, J.~San Martin (2014). 
\newblock Inverse M-Matrices and Ultrametric Matrices. 
\newblock {\em Springer Lecture Notes in Maths,} 2118.

\bibitem{Ege}
E.~Egerv\'ary (1940).
\newblock On orthocentric simplices.
\newblock {\em Acta Universitatis Szegediensis. Acta Scientiarum Mathematicarum,} IX:218--226.

\bibitem{Fie}
M. Fiedler (1957).
\newblock \"{U}ber qualitative Winkeleigenschaften der Simplexe.
\newblock {\em Czechoslovak Mathematical Journal,} 7(82):463--478.

\bibitem{Fie2}
M.~Fiedler (1998).
\newblock Ultrametric sets in Euclidean point spaces.
\newblock {\em Electronic Journal of Linear Algebra,} 3:23--30.
  
\bibitem{Fie3} 
M.~Fiedler (2001).
\newblock Special ultrametric matrices and graphs.
\newblock {\em SIAM Journal on Matrix Analysis and Applications,} 22:106--113.
 
\bibitem{Fie4}
M.~Fiedler (2011).
\newblock Matrices and Graphs in Geometry.
\newblock {\em Encyclopedia of Mathematics and Its Applications,} 139. Cambridge University Press, Cambridge.

\bibitem{FiPt}
M.~Fiedler and V.~Pt\'ak (1962).
\newblock On matrices with non-positive off-diagonal elements and positive principal minors. 
\newblock {\em Czechoslovak Mathematical Journal,} 12(87):382--400.

\bibitem{HaWe}
H.~Havlicek and G.~Wei\ss\,(2003).
\newblock Altitudes of a tetrahedron and traceless quadratic forms.
\newblock {\em The American Mathematical Monthly,} 110(8):679--693.

\bibitem{Joh}
C.R.~Johnson (1982).
\newblock {\em Inverse M-matrices.} 
\newblock {\em Linear Algebra and its Applications,} 47:195--216.
 
\bibitem{JoSm}
C.R.~Johnson and R.L.~Smith (2011).
\newblock Inverse M-matrices, II.
\newblock {\em Linear Algebra and its Applications,} 435:953--983.

\bibitem{KaZi} 
G. Kalai and G.M. Ziegler (Eds) (1997). 
\newblock Lectures on 0/1-polytopes. Polytopes—combinatorics and computation .
\newblock {\em DMV Seminar,} Band 29, Birkh\"auser Verlag, Basel, Boston, Berlin.

\bibitem{Mar}
T.L.~Markham (1972).
\newblock Nonnegative matrices whose inverses are M-matrices.
\newblock {\em Proceedings of the American Mathematical Society,} 36:326--330.

\bibitem{MaMiSa}
S.~Martinez, G.~Michon, and J.~San Martin (1994).
\newblock Inverses of ultrametric matrices are of Stieltjes type. 
\newblock {\em SIAM Journal on Matrix Analysis and Applications,} 15(1):98--106.

\bibitem{McNeSchTs}
J.J.~McDonald, M.~Neumann, H.~Schneider, and M.J.~Tsatsomeros (1995).
\newblock Inverse M-matrix inequalities and generalized ultrametric matrices.
\newblock {\em Linear Algebra and its Applications,} 220:321--341.

\bibitem{NabVar}
R.~Nabben and R.~Varga (1994).
\newblock A linear algebra proof that the inverse of a strictly ultrametric matrix is a strictly diagonally dominant Stieltjes matrix.
\newblock {\em SIAM Journal on Matrix Analysis and Applications,} 15(1):107--113.

\bibitem{NabVar3}
R.~Nabben and R.~Varga (1995).
\newblock Generalized ultrametric matrices - a class of inverse M-matrices.
\newblock {\em Linear Algebra and its Applications,} 220:365--390.
 
\bibitem{NabVar2}
R.~Nabben and R.~Varga (1995).
\newblock On classes of inverse Z-matrices.
\newblock {\em Linear Algebra and its Applications,} 223:521--552. 

\bibitem{Schl}
L.~Schl\"afli (1950).
\newblock Gesammelte Mathematische Abhandlungen I und II.
\newblock {\em Verlag Birkhäuser,} Basel.
   
\bibitem{VarNab}
R.S.~Varga and R.~Nabben (1994).
\newblock On symmetric ultrametric matrices.
\newblock {\em Numerical Linear Algebra,} (L. Reichel, A. Ruttan, and R.S. Varga, eds.), pp. 193--199.

\bibitem{Wil} 
R.A.~Willoughby (1977).
\newblock The inverse M-matrix problem.
\newblock {\em Linear Algebra and its Applications,} 18(1):75--94.


\end{thebibliography}
\end{document}